

\ifx\shlhetal\undefinedcontrolsequence\let\shlhetal\relax\fi

\def\fmtname{AmS-TeX}

\def\fmtversion{2.1}
\catcode`\@=11
\ifx\amstexloaded@\relax\catcode`\@=\active
  \endinput\else\let\amstexloaded@\relax\fi
\newlinechar=`\^^J
\def\W@{\immediate\write\sixt@@n}
\def\CR@{\W@{^^J\fmtname - Version \fmtversion^^J%
COPYRIGHT 1985, 1990, 1991 - AMERICAN MATHEMATICAL SOCIETY^^J%
Use of this macro package is not restricted provided^^J%
each use is acknowledged upon publication.^^J}}
\CR@ \everyjob{\CR@}
\message{Loading definitions for}
\message{misc utility macros,}
\toksdef\toks@@=2
\long\def\rightappend@#1\to#2{\toks@{\\{#1}}\toks@@
 =\expandafter{#2}\xdef#2{\the\toks@@\the\toks@}\toks@{}\toks@@{}}
\def\alloclist@{}
\newif\ifalloc@
\def\showallocations{{\def\\{\immediate\write\m@ne}\alloclist@}\alloc@true}
\def\alloc@#1#2#3#4#5{\global\advance\count1#1by\@ne
 \ch@ck#1#4#2\allocationnumber=\count1#1
 \global#3#5=\allocationnumber
 \edef\next@{\string#5=\string#2\the\allocationnumber}%
 \expandafter\rightappend@\next@\to\alloclist@}
\newcount\count@@
\newcount\count@@@
\def\FN@{\futurelet\next}
\def\DN@{\def\next@}
\def\DNii@{\def\nextii@}
\def\RIfM@{\relax\ifmmode}
\def\RIfMIfI@{\relax\ifmmode\ifinner}
\def\setboxz@h{\setbox\z@\hbox}
\def\wdz@{\wd\z@}
\def\boxz@{\box\z@}
\def\setbox@ne{\setbox\@ne}
\def\wd@ne{\wd\@ne}
\def\iterate{\body\expandafter\iterate\else\fi}
\def\err@#1{\errmessage{AmS-TeX error: #1}}
\newhelp\defaulthelp@{Sorry, I already gave what help I could...^^J
Maybe you should try asking a human?^^J
An error might have occurred before I noticed any problems.^^J
``If all else fails, read the instructions.''}
\def\Err@{\errhelp\defaulthelp@\err@}
\def\eat@#1{}
\def\in@#1#2{\def\in@@##1#1##2##3\in@@{\ifx\in@##2\in@false\else\in@true\fi}%
 \in@@#2#1\in@\in@@}
\newif\ifin@
\def\space@.{\futurelet\space@\relax}
\space@. %
\newhelp\athelp@
{Only certain combinations beginning with @ make sense to me.^^J
Perhaps you wanted \string\@\space for a printed @?^^J
I've ignored the character or group after @.}
{\catcode`\~=\active 
 \lccode`\~=`\@ \lowercase{\gdef~{\FN@\at@}}}
\def\at@{\let\next@\at@@
 \ifcat\noexpand\next a\else\ifcat\noexpand\next0\else
 \ifcat\noexpand\next\relax\else
   \let\next\at@@@\fi\fi\fi
 \next@}
\def\at@@#1{\expandafter
 \ifx\csname\space @\string#1\endcsname\relax
  \expandafter\at@@@ \else
  \csname\space @\string#1\expandafter\endcsname\fi}
\def\at@@@#1{\errhelp\athelp@ \err@{\Invalid@@ @}}
\def\atdef@#1{\expandafter\def\csname\space @\string#1\endcsname}
\newhelp\defahelp@{If you typed \string\define\space cs instead of
\string\define\string\cs\space^^J
I've substituted an inaccessible control sequence so that your^^J
definition will be completed without mixing me up too badly.^^J
If you typed \string\define{\string\cs} the inaccessible control sequence^^J
was defined to be \string\cs, and the rest of your^^J
definition appears as input.}
\newhelp\defbhelp@{I've ignored your definition, because it might^^J
conflict with other uses that are important to me.}
\def\define{\FN@\define@}
\def\define@{\ifcat\noexpand\next\relax
 \expandafter\define@@\else\errhelp\defahelp@                               
 \err@{\string\define\space must be followed by a control
 sequence}\expandafter\def\expandafter\nextii@\fi}                          
\def\undefined@@@@@@@@@@{}
\def\preloaded@@@@@@@@@@{}
\def\next@@@@@@@@@@{}
\def\define@@#1{\ifx#1\relax\errhelp\defbhelp@                              
 \err@{\string#1\space is already defined}\DN@{\DNii@}\else
 \expandafter\ifx\csname\expandafter\eat@\string                            
 #1@@@@@@@@@@\endcsname\undefined@@@@@@@@@@\errhelp\defbhelp@
 \err@{\string#1\space can't be defined}\DN@{\DNii@}\else
 \expandafter\ifx\csname\expandafter\eat@\string#1\endcsname\relax          
 \global\let#1\undefined\DN@{\def#1}\else\errhelp\defbhelp@
 \err@{\string#1\space is already defined}\DN@{\DNii@}\fi
 \fi\fi\next@}

\def\predefine#1#2{\let#1#2}
\def\undefine#1{\let#1\undefined}
\message{page layout,}
\newdimen\captionwidth@
\captionwidth@\hsize
\advance\captionwidth@-1.5in
\def\pagewidth#1{\hsize#1\relax
 \captionwidth@\hsize\advance\captionwidth@-1.5in}
\def\pageheight#1{\vsize#1\relax}
\def\hcorrection#1{\advance\hoffset#1\relax}
\def\vcorrection#1{\advance\voffset#1\relax}
\message{accents/punctuation,}

\let\graveaccent\`
\let\acuteaccent\'
\let\tildeaccent\~
\let\hataccent\^
\let\underscore\_
\let\B\=
\let\D\.
\let\ic@\/
\def\/{\unskip\ic@}
\def\textfonti{\the\textfont\@ne}
\def\t#1#2{{\edef\next@{\the\font}\textfonti\accent"7F \next@#1#2}}
\def~{\unskip\nobreak\ \ignorespaces}
\def\.{.\spacefactor\@m}
\atdef@;{\leavevmode\null;}
\atdef@:{\leavevmode\null:}
\atdef@?{\leavevmode\null?}
\edef\@{\string @}
\def\relaxnext@{\let\next\relax}
\atdef@-{\relaxnext@\leavevmode
 \DN@{\ifx\next-\DN@-{\FN@\nextii@}\else
  \DN@{\leavevmode\hbox{-}}\fi\next@}%
 \DNii@{\ifx\next-\DN@-{\leavevmode\hbox{---}}\else
  \DN@{\leavevmode\hbox{--}}\fi\next@}%
 \FN@\next@}
\def\srdr@{\kern.16667em}
\def\drsr@{\kern.02778em}
\def\sldl@{\drsr@}
\def\dlsl@{\srdr@}
\atdef@"{\unskip\relaxnext@
 \DN@{\ifx\next\space@\DN@. {\FN@\nextii@}\else
  \DN@.{\FN@\nextii@}\fi\next@.}%
 \DNii@{\ifx\next`\DN@`{\FN@\nextiii@}\else
  \ifx\next\lq\DN@\lq{\FN@\nextiii@}\else
  \DN@####1{\FN@\nextiv@}\fi\fi\next@}%
 \def\nextiii@{\ifx\next`\DN@`{\sldl@``}\else\ifx\next\lq
  \DN@\lq{\sldl@``}\else\DN@{\dlsl@`}\fi\fi\next@}%
 \def\nextiv@{\ifx\next'\DN@'{\srdr@''}\else
  \ifx\next\rq\DN@\rq{\srdr@''}\else\DN@{\drsr@'}\fi\fi\next@}%
 \FN@\next@}

\def\textfontii{\the\textfont\tw@}
\def\lbrace@{\delimiter"4266308 }
\def\rbrace@{\delimiter"5267309 }
\def\{{\RIfM@\lbrace@\else{\textfontii f}\spacefactor\@m\fi}
\def\}{\RIfM@\rbrace@\else
 \let\@sf\empty\ifhmode\edef\@sf{\spacefactor\the\spacefactor}\fi
 {\textfontii g}\@sf\relax\fi}
\let\lbrace\{
\let\rbrace\}
\def\AmSTeX{{\textfontii A\kern-.1667em%
  \lower.5ex\hbox{M}\kern-.125emS}-\TeX}
\message{line and page breaks,}
\def\vmodeerr@#1{\Err@{\string#1\space not allowed between paragraphs}}
\def\mathmodeerr@#1{\Err@{\string#1\space not allowed in math mode}}
\def\linebreak{\RIfM@\mathmodeerr@\linebreak\else
 \ifhmode\unskip\unkern\break\else\vmodeerr@\linebreak\fi\fi}

\newskip\saveskip@
\def\allowlinebreak{\RIfM@\mathmodeerr@\allowlinebreak\else
 \ifhmode\saveskip@\lastskip\unskip
 \allowbreak\ifdim\saveskip@>\z@\hskip\saveskip@\fi
 \else\vmodeerr@\allowlinebreak\fi\fi}
\def\nolinebreak{\RIfM@\mathmodeerr@\nolinebreak\else
 \ifhmode\saveskip@\lastskip\unskip
 \nobreak\ifdim\saveskip@>\z@\hskip\saveskip@\fi
 \else\vmodeerr@\nolinebreak\fi\fi}
\def\newline{\relaxnext@
 \DN@{\RIfM@\expandafter\mathmodeerr@\expandafter\newline\else
  \ifhmode\ifx\next\par\else
  \expandafter\unskip\expandafter\null\expandafter\hfill\expandafter\break\fi
  \else
  \expandafter\vmodeerr@\expandafter\newline\fi\fi}%
 \FN@\next@}
\def\dmatherr@#1{\Err@{\string#1\space not allowed in display math mode}}
\def\nondmatherr@#1{\Err@{\string#1\space not allowed in non-display math
 mode}}
\def\onlydmatherr@#1{\Err@{\string#1\space allowed only in display math mode}}
\def\nonmatherr@#1{\Err@{\string#1\space allowed only in math mode}}
\def\mathbreak{\RIfMIfI@\break\else
 \dmatherr@\mathbreak\fi\else\nonmatherr@\mathbreak\fi}
\def\nomathbreak{\RIfMIfI@\nobreak\else
 \dmatherr@\nomathbreak\fi\else\nonmatherr@\nomathbreak\fi}
\def\allowmathbreak{\RIfMIfI@\allowbreak\else
 \dmatherr@\allowmathbreak\fi\else\nonmatherr@\allowmathbreak\fi}
\def\pagebreak{\RIfM@
 \ifinner\nondmatherr@\pagebreak\else\postdisplaypenalty-\@M\fi
 \else\ifvmode\removelastskip\break\else\vadjust{\break}\fi\fi}
\def\nopagebreak{\RIfM@
 \ifinner\nondmatherr@\nopagebreak\else\postdisplaypenalty\@M\fi
 \else\ifvmode\nobreak\else\vadjust{\nobreak}\fi\fi}
\def\nonvmodeerr@#1{\Err@{\string#1\space not allowed within a paragraph
 or in math}}
\def\vnonvmode@#1#2{\relaxnext@\DNii@{\ifx\next\par\DN@{#1}\else
 \DN@{#2}\fi\next@}%
 \ifvmode\DN@{#1}\else
 \DN@{\FN@\nextii@}\fi\next@}
\def\newpage{\vnonvmode@{\vfill\break}{\nonvmodeerr@\newpage}}
\def\smallpagebreak{\vnonvmode@\smallbreak{\nonvmodeerr@\smallpagebreak}}
\def\medpagebreak{\vnonvmode@\medbreak{\nonvmodeerr@\medpagebreak}}
\def\bigpagebreak{\vnonvmode@\bigbreak{\nonvmodeerr@\bigpagebreak}}
\def\NoBlackBoxes{\global\overfullrule\z@}
\def\BlackBoxes{\global\overfullrule5\p@}
\def\Invalid@#1{\def#1{\Err@{\Invalid@@\string#1}}}
\def\Invalid@@{Invalid use of }
\message{figures,}
\Invalid@\caption
\Invalid@\captionwidth
\newdimen\smallcaptionwidth@
\def\topspace{\mid@false\ins@}
\def\midspace{\mid@true\ins@}
\newif\ifmid@
\def\captionfont@{}
\def\ins@#1{\relaxnext@\allowbreak
 \smallcaptionwidth@\captionwidth@\gdef\thespace@{#1}%
 \DN@{\ifx\next\space@\DN@. {\FN@\nextii@}\else
  \DN@.{\FN@\nextii@}\fi\next@.}%
 \DNii@{\ifx\next\caption\DN@\caption{\FN@\nextiii@}%
  \else\let\next@\nextiv@\fi\next@}%
 \def\nextiv@{\vnonvmode@
  {\ifmid@\expandafter\midinsert\else\expandafter\topinsert\fi
   \vbox to\thespace@{}\endinsert}
  {\ifmid@\nonvmodeerr@\midspace\else\nonvmodeerr@\topspace\fi}}%
 \def\nextiii@{\ifx\next\captionwidth\expandafter\nextv@
  \else\expandafter\nextvi@\fi}%
 \def\nextv@\captionwidth##1##2{\smallcaptionwidth@##1\relax\nextvi@{##2}}%
 \def\nextvi@##1{\def\thecaption@{\captionfont@##1}%
  \DN@{\ifx\next\space@\DN@. {\FN@\nextvii@}\else
   \DN@.{\FN@\nextvii@}\fi\next@.}%
  \FN@\next@}%
 \def\nextvii@{\vnonvmode@
  {\ifmid@\expandafter\midinsert\else
  \expandafter\topinsert\fi\vbox to\thespace@{}\nobreak\smallskip
  \setboxz@h{\noindent\ignorespaces\thecaption@\unskip}%
  \ifdim\wdz@>\smallcaptionwidth@\centerline{\vbox{\hsize\smallcaptionwidth@
   \noindent\ignorespaces\thecaption@\unskip}}%
  \else\centerline{\boxz@}\fi\endinsert}
  {\ifmid@\nonvmodeerr@\midspace
  \else\nonvmodeerr@\topspace\fi}}%
 \FN@\next@}
\message{comments,}
\def\newcodes@{\catcode`\\12\catcode`\{12\catcode`\}12\catcode`\#12%
 \catcode`\%12\relax}
\def\oldcodes@{\catcode`\\0\catcode`\{1\catcode`\}2\catcode`\#6%
 \catcode`\%14\relax}
\def\comment{\newcodes@\endlinechar=10 \comment@}
{\lccode`\0=`\\
\lowercase{\gdef\comment@#1^^J{\comment@@#10endcomment\comment@@@}%
\gdef\comment@@#10endcomment{\FN@\comment@@@}%
\gdef\comment@@@#1\comment@@@{\ifx\next\comment@@@\let\next\comment@
 \else\def\next{\oldcodes@\endlinechar=`\^^M\relax}%
 \fi\next}}}
\def\pr@m@s{\ifx'\next\DN@##1{\prim@s}\else\let\next@\egroup\fi\next@}
\def\prime{{\null\prime@\null}}
\mathchardef\prime@="0230
\let\dsize\displaystyle

\let\ssize\scriptstyle

\message{math spacing,}
\def\,{\RIfM@\mskip\thinmuskip\relax\else\kern.16667em\fi}
\def\!{\RIfM@\mskip-\thinmuskip\relax\else\kern-.16667em\fi}
\let\thinspace\,
\let\negthinspace\!
\def\medspace{\RIfM@\mskip\medmuskip\relax\else\kern.222222em\fi}
\def\negmedspace{\RIfM@\mskip-\medmuskip\relax\else\kern-.222222em\fi}
\def\thickspace{\RIfM@\mskip\thickmuskip\relax\else\kern.27777em\fi}
\let\;\thickspace
\def\negthickspace{\RIfM@\mskip-\thickmuskip\relax\else
 \kern-.27777em\fi}
\atdef@,{\RIfM@\mskip.1\thinmuskip\else\leavevmode\null,\fi}
\atdef@!{\RIfM@\mskip-.1\thinmuskip\else\leavevmode\null!\fi}
\atdef@.{\RIfM@&&\else\leavevmode.\spacefactor3000 \fi}
\def\and{\DOTSB\;\mathchar"3026 \;}

\message{fractions,}
\def\frac#1#2{{#1\over#2}}

\newdimen\ex@
\ex@.2326ex
\Invalid@\thickness
\def\thickfrac{\relaxnext@
 \DN@{\ifx\next\thickness\let\next@\nextii@\else
 \DN@{\nextii@\thickness1}\fi\next@}%
 \DNii@\thickness##1##2##3{{##2\above##1\ex@##3}}%
 \FN@\next@}

\def\thickfracwithdelims#1#2{\relaxnext@\def\ldelim@{#1}\def\rdelim@{#2}%
 \DN@{\ifx\next\thickness\let\next@\nextii@\else
 \DN@{\nextii@\thickness1}\fi\next@}%
 \DNii@\thickness##1##2##3{{##2\abovewithdelims
 \ldelim@\rdelim@##1\ex@##3}}%
 \FN@\next@}

\def\:{\nobreak\hskip.1111em\mathpunct{}\nonscript\mkern-\thinmuskip{:}\hskip
 .3333emplus.0555em\relax}
\def\snug{\unskip\kern-\mathsurround}
\message{smash commands,}
\def\topsmash{\top@true\bot@false\smash@}
\def\botsmash{\top@false\bot@true\smash@}
\newif\iftop@
\newif\ifbot@
\def\smash{\top@true\bot@true\smash@}
\def\smash@{\RIfM@\expandafter\mathpalette\expandafter\mathsm@sh\else
 \expandafter\makesm@sh\fi}
\def\finsm@sh{\iftop@\ht\z@\z@\fi\ifbot@\dp\z@\z@\fi\leavevmode\boxz@}
\message{large operator symbols,}
\def\LimitsOnSums{\global\let\slimits@\displaylimits}
\def\NoLimitsOnSums{\global\let\slimits@\nolimits}
\LimitsOnSums
\mathchardef\coprod@="1360       \def\coprod{\DOTSB\coprod@\slimits@}
\mathchardef\bigvee@="1357       \def\bigvee{\DOTSB\bigvee@\slimits@}
\mathchardef\bigwedge@="1356     \def\bigwedge{\DOTSB\bigwedge@\slimits@}
\mathchardef\biguplus@="1355     \def\biguplus{\DOTSB\biguplus@\slimits@}
\mathchardef\bigcap@="1354       \def\bigcap{\DOTSB\bigcap@\slimits@}
\mathchardef\bigcup@="1353       \def\bigcup{\DOTSB\bigcup@\slimits@}
\mathchardef\prod@="1351         \def\prod{\DOTSB\prod@\slimits@}
\mathchardef\sum@="1350          \def\sum{\DOTSB\sum@\slimits@}
\mathchardef\bigotimes@="134E    \def\bigotimes{\DOTSB\bigotimes@\slimits@}
\mathchardef\bigoplus@="134C     \def\bigoplus{\DOTSB\bigoplus@\slimits@}
\mathchardef\bigodot@="134A      \def\bigodot{\DOTSB\bigodot@\slimits@}
\mathchardef\bigsqcup@="1346     \def\bigsqcup{\DOTSB\bigsqcup@\slimits@}
\message{integrals,}
\def\LimitsOnInts{\global\let\ilimits@\displaylimits}
\def\NoLimitsOnInts{\global\let\ilimits@\nolimits}
\NoLimitsOnInts
\def\int{\DOTSI\intop\ilimits@}
\def\oint{\DOTSI\ointop\ilimits@}
\def\intic@{\mathchoice{\hskip.5em}{\hskip.4em}{\hskip.4em}{\hskip.4em}}
\def\negintic@{\mathchoice
 {\hskip-.5em}{\hskip-.4em}{\hskip-.4em}{\hskip-.4em}}
\def\intkern@{\mathchoice{\!\!\!}{\!\!}{\!\!}{\!\!}}
\def\intdots@{\mathchoice{\plaincdots@}
 {{\cdotp}\mkern1.5mu{\cdotp}\mkern1.5mu{\cdotp}}
 {{\cdotp}\mkern1mu{\cdotp}\mkern1mu{\cdotp}}
 {{\cdotp}\mkern1mu{\cdotp}\mkern1mu{\cdotp}}}
\newcount\intno@
\def\iint{\DOTSI\intno@\tw@\FN@\ints@}
\def\iiint{\DOTSI\intno@\thr@@\FN@\ints@}
\def\iiiint{\DOTSI\intno@4 \FN@\ints@}
\def\idotsint{\DOTSI\intno@\z@\FN@\ints@}
\def\ints@{\findlimits@\ints@@}
\newif\iflimtoken@
\newif\iflimits@
\def\findlimits@{\limtoken@true\ifx\next\limits\limits@true
 \else\ifx\next\nolimits\limits@false\else
 \limtoken@false\ifx\ilimits@\nolimits\limits@false\else
 \ifinner\limits@false\else\limits@true\fi\fi\fi\fi}
\def\multint@{\int\ifnum\intno@=\z@\intdots@                                
 \else\intkern@\fi                                                          
 \ifnum\intno@>\tw@\int\intkern@\fi                                         
 \ifnum\intno@>\thr@@\int\intkern@\fi                                       
 \int}                                                                      
\def\multintlimits@{\intop\ifnum\intno@=\z@\intdots@\else\intkern@\fi
 \ifnum\intno@>\tw@\intop\intkern@\fi
 \ifnum\intno@>\thr@@\intop\intkern@\fi\intop}
\def\ints@@{\iflimtoken@                                                    
 \def\ints@@@{\iflimits@\negintic@\mathop{\intic@\multintlimits@}\limits    
  \else\multint@\nolimits\fi                                                
  \eat@}                                                                    
 \else                                                                      
 \def\ints@@@{\iflimits@\negintic@
  \mathop{\intic@\multintlimits@}\limits\else
  \multint@\nolimits\fi}\fi\ints@@@}
\def\LimitsOnNames{\global\let\nlimits@\displaylimits}
\def\NoLimitsOnNames{\global\let\nlimits@\nolimits@}
\LimitsOnNames
\def\nolimits@{\relaxnext@
 \DN@{\ifx\next\limits\DN@\limits{\nolimits}\else
  \let\next@\nolimits\fi\next@}%
 \FN@\next@}
\message{operator names,}
\def\newmcodes@{\mathcode`\'"27\mathcode`\*"2A\mathcode`\."613A%
 \mathcode`\-"2D\mathcode`\/"2F\mathcode`\:"603A }
\def\operatorname#1{\mathop{\newmcodes@\kern\z@\fam\z@#1}\nolimits@}
\def\operatornamewithlimits#1{\mathop{\newmcodes@\kern\z@\fam\z@#1}\nlimits@}
\def\qopname@#1{\mathop{\fam\z@#1}\nolimits@}
\def\qopnamewl@#1{\mathop{\fam\z@#1}\nlimits@}
\def\arccos{\qopname@{arccos}}
\def\arcsin{\qopname@{arcsin}}
\def\arctan{\qopname@{arctan}}
\def\arg{\qopname@{arg}}
\def\cos{\qopname@{cos}}
\def\cosh{\qopname@{cosh}}
\def\cot{\qopname@{cot}}
\def\coth{\qopname@{coth}}
\def\csc{\qopname@{csc}}
\def\deg{\qopname@{deg}}
\def\det{\qopnamewl@{det}}
\def\dim{\qopname@{dim}}
\def\exp{\qopname@{exp}}
\def\gcd{\qopnamewl@{gcd}}
\def\hom{\qopname@{hom}}
\def\inf{\qopnamewl@{inf}}
\def\injlim{\qopnamewl@{inj\,lim}}
\def\ker{\qopname@{ker}}
\def\lg{\qopname@{lg}}
\def\lim{\qopnamewl@{lim}}
\def\liminf{\qopnamewl@{lim\,inf}}
\def\limsup{\qopnamewl@{lim\,sup}}
\def\ln{\qopname@{ln}}
\def\log{\qopname@{log}}
\def\max{\qopnamewl@{max}}
\def\min{\qopnamewl@{min}}
\def\Pr{\qopnamewl@{Pr}}
\def\projlim{\qopnamewl@{proj\,lim}}
\def\sec{\qopname@{sec}}
\def\sin{\qopname@{sin}}
\def\sinh{\qopname@{sinh}}
\def\sup{\qopnamewl@{sup}}
\def\tan{\qopname@{tan}}
\def\tanh{\qopname@{tanh}}
\def\varinjlim{\mathop{\vtop{\ialign{##\crcr
 \hfil\rm lim\hfil\crcr\noalign{\nointerlineskip}\rightarrowfill\crcr
 \noalign{\nointerlineskip\kern-\ex@}\crcr}}}}
\def\varprojlim{\mathop{\vtop{\ialign{##\crcr
 \hfil\rm lim\hfil\crcr\noalign{\nointerlineskip}\leftarrowfill\crcr
 \noalign{\nointerlineskip\kern-\ex@}\crcr}}}}
\def\varliminf{\mathop{\underline{\vrule height\z@ depth.2exwidth\z@
 \hbox{\rm lim}}}}

\newdimen\buffer@
\buffer@\fontdimen13 \tenex
\newdimen\buffer
\buffer\buffer@

\def\ResetBuffer{\fontdimen13 \tenex\buffer@\global\buffer\buffer@}
\def\shave#1{\mathop{\hbox{$\m@th\fontdimen13 \tenex\z@                     
 \displaystyle{#1}$}}\fontdimen13 \tenex\buffer}

\message{multilevel sub/superscripts,}
\Invalid@\\
\def\Let@{\relax\iffalse{\fi\let\\=\cr\iffalse}\fi}
\Invalid@\vspace
\def\vspace@{\def\vspace##1{\crcr\noalign{\vskip##1\relax}}}
\def\multilimits@{\bgroup\vspace@\Let@
 \baselineskip\fontdimen10 \scriptfont\tw@
 \advance\baselineskip\fontdimen12 \scriptfont\tw@
 \lineskip\thr@@\fontdimen8 \scriptfont\thr@@
 \lineskiplimit\lineskip
 \vbox\bgroup\ialign\bgroup\hfil$\m@th\scriptstyle{##}$\hfil\crcr}
\def\Sb{_\multilimits@}
\def\endSb{\crcr\egroup\egroup\egroup}
\def\Sp{^\multilimits@}

\def\spreadlines#1{\RIfMIfI@\onlydmatherr@\spreadlines\else
 \openup#1\relax\fi\else\onlydmatherr@\spreadlines\fi}
\def\Mathstrut@{\copy\Mathstrutbox@}
\newbox\Mathstrutbox@
\setbox\Mathstrutbox@\null
\setboxz@h{$\m@th($}
\ht\Mathstrutbox@\ht\z@
\dp\Mathstrutbox@\dp\z@
\message{matrices,}
\newdimen\spreadmlines@
\def\spreadmatrixlines#1{\RIfMIfI@
 \onlydmatherr@\spreadmatrixlines\else
 \spreadmlines@#1\relax\fi\else\onlydmatherr@\spreadmatrixlines\fi}
\def\matrix{\null\,\vcenter\bgroup\Let@\vspace@
 \normalbaselines\openup\spreadmlines@\ialign
 \bgroup\hfil$\m@th##$\hfil&&\quad\hfil$\m@th##$\hfil\crcr
 \Mathstrut@\crcr\noalign{\kern-\baselineskip}}
\def\endmatrix{\crcr\Mathstrut@\crcr\noalign{\kern-\baselineskip}\egroup
 \egroup\,}
\def\format{\crcr\egroup\iffalse{\fi\ifnum`}=0 \fi\format@}
\newtoks\hashtoks@
\hashtoks@{#}
\def\format@#1\\{\def\preamble@{#1}%
 \def\l{$\m@th\the\hashtoks@$\hfil}%
 \def\c{\hfil$\m@th\the\hashtoks@$\hfil}%
 \def\r{\hfil$\m@th\the\hashtoks@$}%
 \edef\preamble@@{\preamble@}\ifnum`{=0 \fi\iffalse}\fi
 \ialign\bgroup\span\preamble@@\crcr}
\def\smallmatrix{\null\,\vcenter\bgroup\vspace@\Let@
 \baselineskip9\ex@\lineskip\ex@
 \ialign\bgroup\hfil$\m@th\scriptstyle{##}$\hfil&&\thickspace\hfil
 $\m@th\scriptstyle{##}$\hfil\crcr}
\def\endsmallmatrix{\crcr\egroup\egroup\,}

\newmuskip\dotsspace@
\dotsspace@1.5mu
\def\strip@#1 {#1}
\def\spacehdots#1\for#2{\multispan{#2}\xleaders
 \hbox{$\m@th\mkern\strip@#1 \dotsspace@.\mkern\strip@#1 \dotsspace@$}\hfill}
\def\hdotsfor#1{\spacehdots\@ne\for{#1}}
\def\multispan@#1{\omit\mscount#1\unskip\loop\ifnum\mscount>\@ne\sp@n\repeat}
\def\spaceinnerhdots#1\for#2\after#3{\multispan@{\strip@#2 }#3\xleaders
 \hbox{$\m@th\mkern\strip@#1 \dotsspace@.\mkern\strip@#1 \dotsspace@$}\hfill}
\def\innerhdotsfor#1\after#2{\spaceinnerhdots\@ne\for#1\after{#2}}
\def\cases{\bgroup\spreadmlines@\jot\left\{\,\matrix\format\l&\quad\l\\}
\def\endcases{\endmatrix\right.\egroup}
\message{multiline displays,}
\newif\ifinany@
\newif\ifinalign@
\newif\ifingather@
\def\strut@{\copy\strutbox@}
\newbox\strutbox@
\setbox\strutbox@\hbox{\vrule height8\p@ depth3\p@ width\z@}
\def\topaligned{\null\,\vtop\aligned@}
\def\botaligned{\null\,\vbox\aligned@}
\def\aligned{\null\,\vcenter\aligned@}
\def\aligned@{\bgroup\vspace@\Let@
 \ifinany@\else\openup\jot\fi\ialign
 \bgroup\hfil\strut@$\m@th\displaystyle{##}$&
 $\m@th\displaystyle{{}##}$\hfil\crcr}
\def\endaligned{\crcr\egroup\egroup}

\def\alignedat#1{\null\,\vcenter\bgroup\doat@{#1}\vspace@\Let@
 \ifinany@\else\openup\jot\fi\ialign\bgroup\span\preamble@@\crcr}
\newcount\atcount@
\def\doat@#1{\toks@{\hfil\strut@$\m@th
 \displaystyle{\the\hashtoks@}$&$\m@th\displaystyle
 {{}\the\hashtoks@}$\hfil}
 \atcount@#1\relax\advance\atcount@\m@ne                                    
 \loop\ifnum\atcount@>\z@\toks@=\expandafter{\the\toks@&\hfil$\m@th
 \displaystyle{\the\hashtoks@}$&$\m@th
 \displaystyle{{}\the\hashtoks@}$\hfil}\advance
  \atcount@\m@ne\repeat                                                     
 \xdef\preamble@{\the\toks@}\xdef\preamble@@{\preamble@}}

\def\gathered{\null\,\vcenter\bgroup\vspace@\Let@
 \ifinany@\else\openup\jot\fi\ialign
 \bgroup\hfil\strut@$\m@th\displaystyle{##}$\hfil\crcr}
\def\endgathered{\crcr\egroup\egroup}
\newif\iftagsleft@
\def\TagsOnLeft{\global\tagsleft@true}
\def\TagsOnRight{\global\tagsleft@false}
\TagsOnLeft
\newif\ifmathtags@
\def\TagsAsMath{\global\mathtags@true}
\def\TagsAsText{\global\mathtags@false}
\TagsAsText
\def\tagform@#1{\hbox{\rm(\ignorespaces#1\unskip)}}
\def\thetag{\leavevmode\tagform@}
\def\tag#1$${\iftagsleft@\leqno\else\eqno\fi                                
 \maketag@#1\maketag@                                                       
 $$}                                                                        
\def\maketag@{\FN@\maketag@@}
\def\maketag@@{\ifx\next"\expandafter\maketag@@@\else\expandafter\maketag@@@@
 \fi}
\def\maketag@@@"#1"#2\maketag@{\hbox{\rm#1}}                                
\def\maketag@@@@#1\maketag@{\ifmathtags@\tagform@{$\m@th#1$}\else
 \tagform@{#1}\fi}
\interdisplaylinepenalty\@M
\def\allowdisplaybreaks{\RIfMIfI@
 \onlydmatherr@\allowdisplaybreaks\else
 \interdisplaylinepenalty\z@\fi\else\onlydmatherr@\allowdisplaybreaks\fi}
\Invalid@\allowdisplaybreak
\Invalid@\displaybreak
\Invalid@\intertext
\def\allowdisplaybreak@{\def\allowdisplaybreak{\crcr\noalign{\allowbreak}}}
\def\displaybreak@{\def\displaybreak{\crcr\noalign{\break}}}
\def\intertext@{\def\intertext##1{\crcr\noalign{%
 \penalty\postdisplaypenalty \vskip\belowdisplayskip
 \vbox{\normalbaselines\noindent##1}%
 \penalty\predisplaypenalty \vskip\abovedisplayskip}}}
\newskip\centering@
\centering@\z@ plus\@m\p@
\def\align{\relax\ifingather@\DN@{\csname align (in
  \string\gather)\endcsname}\else
 \ifmmode\ifinner\DN@{\onlydmatherr@\align}\else
  \let\next@\align@\fi
 \else\DN@{\onlydmatherr@\align}\fi\fi\next@}
\newhelp\andhelp@
{An extra & here is so disastrous that you should probably exit^^J
and fix things up.}
\newif\iftag@
\newcount\and@
\def\align@{\inalign@true\inany@true
 \vspace@\allowdisplaybreak@\displaybreak@\intertext@
 \def\tag{\global\tag@true\ifnum\and@=\z@\DN@{&&}\else
          \DN@{&}\fi\next@}%
 \iftagsleft@\DN@{\csname align \endcsname}\else
  \DN@{\csname align \space\endcsname}\fi\next@}
\def\Tag@{\iftag@\else\errhelp\andhelp@\err@{Extra & on this line}\fi}
\newdimen\lwidth@
\newdimen\rwidth@
\newdimen\maxlwidth@
\newdimen\maxrwidth@
\newdimen\totwidth@
\def\measure@#1\endalign{\lwidth@\z@\rwidth@\z@\maxlwidth@\z@\maxrwidth@\z@
 \global\and@\z@                                                            
 \setbox@ne\vbox                                                            
  {\everycr{\noalign{\global\tag@false\global\and@\z@}}\Let@                
  \halign{\setboxz@h{$\m@th\displaystyle{\@lign##}$}
   \global\lwidth@\wdz@                                                     
   \ifdim\lwidth@>\maxlwidth@\global\maxlwidth@\lwidth@\fi                  
   \global\advance\and@\@ne                                                 
   &\setboxz@h{$\m@th\displaystyle{{}\@lign##}$}\global\rwidth@\wdz@        
   \ifdim\rwidth@>\maxrwidth@\global\maxrwidth@\rwidth@\fi                  
   \global\advance\and@\@ne                                                
   &\Tag@
   \eat@{##}\crcr#1\crcr}}
 \totwidth@\maxlwidth@\advance\totwidth@\maxrwidth@}                       
\def\displ@y@{\global\dt@ptrue\openup\jot
 \everycr{\noalign{\global\tag@false\global\and@\z@\ifdt@p\global\dt@pfalse
 \vskip-\lineskiplimit\vskip\normallineskiplimit\else
 \penalty\interdisplaylinepenalty\fi}}}
\def\black@#1{\noalign{\ifdim#1>\displaywidth
 \dimen@\prevdepth\nointerlineskip                                          
 \vskip-\ht\strutbox@\vskip-\dp\strutbox@                                   
 \vbox{\noindent\hbox to#1{\strut@\hfill}}
 \prevdepth\dimen@                                                          
 \fi}}
\expandafter\def\csname align \space\endcsname#1\endalign
 {\measure@#1\endalign\global\and@\z@                                       
 \ifingather@\everycr{\noalign{\global\and@\z@}}\else\displ@y@\fi           
 \Let@\tabskip\centering@                                                   
 \halign to\displaywidth
  {\hfil\strut@\setboxz@h{$\m@th\displaystyle{\@lign##}$}
  \global\lwidth@\wdz@\boxz@\global\advance\and@\@ne                        
  \tabskip\z@skip                                                           
  &\setboxz@h{$\m@th\displaystyle{{}\@lign##}$}
  \global\rwidth@\wdz@\boxz@\hfill\global\advance\and@\@ne                  
  \tabskip\centering@                                                       
  &\setboxz@h{\@lign\strut@\maketag@##\maketag@}
  \dimen@\displaywidth\advance\dimen@-\totwidth@
  \divide\dimen@\tw@\advance\dimen@\maxrwidth@\advance\dimen@-\rwidth@     
  \ifdim\dimen@<\tw@\wdz@\llap{\vtop{\normalbaselines\null\boxz@}}
  \else\llap{\boxz@}\fi                                                    
  \tabskip\z@skip                                                          
  \crcr#1\crcr                                                             
  \black@\totwidth@}}                                                      
\newdimen\lineht@
\expandafter\def\csname align \endcsname#1\endalign{\measure@#1\endalign
 \global\and@\z@
 \ifdim\totwidth@>\displaywidth\let\displaywidth@\totwidth@\else
  \let\displaywidth@\displaywidth\fi                                        
 \ifingather@\everycr{\noalign{\global\and@\z@}}\else\displ@y@\fi
 \Let@\tabskip\centering@\halign to\displaywidth
  {\hfil\strut@\setboxz@h{$\m@th\displaystyle{\@lign##}$}%
  \global\lwidth@\wdz@\global\lineht@\ht\z@                                 
  \boxz@\global\advance\and@\@ne
  \tabskip\z@skip&\setboxz@h{$\m@th\displaystyle{{}\@lign##}$}%
  \global\rwidth@\wdz@\ifdim\ht\z@>\lineht@\global\lineht@\ht\z@\fi         
  \boxz@\hfil\global\advance\and@\@ne
  \tabskip\centering@&\kern-\displaywidth@                                  
  \setboxz@h{\@lign\strut@\maketag@##\maketag@}%
  \dimen@\displaywidth\advance\dimen@-\totwidth@
  \divide\dimen@\tw@\advance\dimen@\maxlwidth@\advance\dimen@-\lwidth@
  \ifdim\dimen@<\tw@\wdz@
   \rlap{\vbox{\normalbaselines\boxz@\vbox to\lineht@{}}}\else
   \rlap{\boxz@}\fi
  \tabskip\displaywidth@\crcr#1\crcr\black@\totwidth@}}
\expandafter\def\csname align (in \string\gather)\endcsname
  #1\endalign{\vcenter{\align@#1\endalign}}
\Invalid@\endalign
\newif\ifxat@
\def\alignat{\RIfMIfI@\DN@{\onlydmatherr@\alignat}\else
 \DN@{\csname alignat \endcsname}\fi\else
 \DN@{\onlydmatherr@\alignat}\fi\next@}
\newif\ifmeasuring@
\newbox\savealignat@
\expandafter\def\csname alignat \endcsname#1#2\endalignat                   
 {\inany@true\xat@false
 \def\tag{\global\tag@true\count@#1\relax\multiply\count@\tw@
  \xdef\tag@{}\loop\ifnum\count@>\and@\xdef\tag@{&\tag@}\advance\count@\m@ne
  \repeat\tag@}%
 \vspace@\allowdisplaybreak@\displaybreak@\intertext@
 \displ@y@\measuring@true                                                   
 \setbox\savealignat@\hbox{$\m@th\displaystyle\Let@
  \attag@{#1}
  \vbox{\halign{\span\preamble@@\crcr#2\crcr}}$}%
 \measuring@false                                                           
 \Let@\attag@{#1}
 \tabskip\centering@\halign to\displaywidth
  {\span\preamble@@\crcr#2\crcr                                             
  \black@{\wd\savealignat@}}}                                               
\Invalid@\endalignat
\def\xalignat{\RIfMIfI@
 \DN@{\onlydmatherr@\xalignat}\else
 \DN@{\csname xalignat \endcsname}\fi\else
 \DN@{\onlydmatherr@\xalignat}\fi\next@}
\expandafter\def\csname xalignat \endcsname#1#2\endxalignat
 {\inany@true\xat@true
 \def\tag{\global\tag@true\def\tag@{}\count@#1\relax\multiply\count@\tw@
  \loop\ifnum\count@>\and@\xdef\tag@{&\tag@}\advance\count@\m@ne\repeat\tag@}%
 \vspace@\allowdisplaybreak@\displaybreak@\intertext@
 \displ@y@\measuring@true\setbox\savealignat@\hbox{$\m@th\displaystyle\Let@
 \attag@{#1}\vbox{\halign{\span\preamble@@\crcr#2\crcr}}$}%
 \measuring@false\Let@
 \attag@{#1}\tabskip\centering@\halign to\displaywidth
 {\span\preamble@@\crcr#2\crcr\black@{\wd\savealignat@}}}
\def\attag@#1{\let\Maketag@\maketag@\let\TAG@\Tag@                          
 \let\Tag@=0\let\maketag@=0
 \ifmeasuring@\def\llap@##1{\setboxz@h{##1}\hbox to\tw@\wdz@{}}%
  \def\rlap@##1{\setboxz@h{##1}\hbox to\tw@\wdz@{}}\else
  \let\llap@\llap\let\rlap@\rlap\fi                                         
 \toks@{\hfil\strut@$\m@th\displaystyle{\@lign\the\hashtoks@}$\tabskip\z@skip
  \global\advance\and@\@ne&$\m@th\displaystyle{{}\@lign\the\hashtoks@}$\hfil
  \ifxat@\tabskip\centering@\fi\global\advance\and@\@ne}
 \iftagsleft@
  \toks@@{\tabskip\centering@&\Tag@\kern-\displaywidth
   \rlap@{\@lign\maketag@\the\hashtoks@\maketag@}%
   \global\advance\and@\@ne\tabskip\displaywidth}\else
  \toks@@{\tabskip\centering@&\Tag@\llap@{\@lign\maketag@
   \the\hashtoks@\maketag@}\global\advance\and@\@ne\tabskip\z@skip}\fi      
 \atcount@#1\relax\advance\atcount@\m@ne
 \loop\ifnum\atcount@>\z@
 \toks@=\expandafter{\the\toks@&\hfil$\m@th\displaystyle{\@lign
  \the\hashtoks@}$\global\advance\and@\@ne
  \tabskip\z@skip&$\m@th\displaystyle{{}\@lign\the\hashtoks@}$\hfil\ifxat@
  \tabskip\centering@\fi\global\advance\and@\@ne}\advance\atcount@\m@ne
 \repeat                                                                    
 \xdef\preamble@{\the\toks@\the\toks@@}
 \xdef\preamble@@{\preamble@}
 \let\maketag@\Maketag@\let\Tag@\TAG@}                                      
\Invalid@\endxalignat
\def\xxalignat{\RIfMIfI@
 \DN@{\onlydmatherr@\xxalignat}\else\DN@{\csname xxalignat
  \endcsname}\fi\else
 \DN@{\onlydmatherr@\xxalignat}\fi\next@}
\expandafter\def\csname xxalignat \endcsname#1#2\endxxalignat{\inany@true
 \vspace@\allowdisplaybreak@\displaybreak@\intertext@
 \displ@y\setbox\savealignat@\hbox{$\m@th\displaystyle\Let@
 \xxattag@{#1}\vbox{\halign{\span\preamble@@\crcr#2\crcr}}$}%
 \Let@\xxattag@{#1}\tabskip\z@skip\halign to\displaywidth
 {\span\preamble@@\crcr#2\crcr\black@{\wd\savealignat@}}}
\def\xxattag@#1{\toks@{\tabskip\z@skip\hfil\strut@
 $\m@th\displaystyle{\the\hashtoks@}$&%
 $\m@th\displaystyle{{}\the\hashtoks@}$\hfil\tabskip\centering@&}%
 \atcount@#1\relax\advance\atcount@\m@ne\loop\ifnum\atcount@>\z@
 \toks@=\expandafter{\the\toks@&\hfil$\m@th\displaystyle{\the\hashtoks@}$%
  \tabskip\z@skip&$\m@th\displaystyle{{}\the\hashtoks@}$\hfil
  \tabskip\centering@}\advance\atcount@\m@ne\repeat
 \xdef\preamble@{\the\toks@\tabskip\z@skip}\xdef\preamble@@{\preamble@}}
\Invalid@\endxxalignat
\newdimen\gwidth@
\newdimen\gmaxwidth@
\def\gmeasure@#1\endgather{\gwidth@\z@\gmaxwidth@\z@\setbox@ne\vbox{\Let@
 \halign{\setboxz@h{$\m@th\displaystyle{##}$}\global\gwidth@\wdz@
 \ifdim\gwidth@>\gmaxwidth@\global\gmaxwidth@\gwidth@\fi
 &\eat@{##}\crcr#1\crcr}}}
\def\gather{\RIfMIfI@\DN@{\onlydmatherr@\gather}\else
 \ingather@true\inany@true\def\tag{&}%
 \vspace@\allowdisplaybreak@\displaybreak@\intertext@
 \displ@y\Let@
 \iftagsleft@\DN@{\csname gather \endcsname}\else
  \DN@{\csname gather \space\endcsname}\fi\fi
 \else\DN@{\onlydmatherr@\gather}\fi\next@}
\expandafter\def\csname gather \space\endcsname#1\endgather
 {\gmeasure@#1\endgather\tabskip\centering@
 \halign to\displaywidth{\hfil\strut@\setboxz@h{$\m@th\displaystyle{##}$}%
 \global\gwidth@\wdz@\boxz@\hfil&
 \setboxz@h{\strut@{\maketag@##\maketag@}}%
 \dimen@\displaywidth\advance\dimen@-\gwidth@
 \ifdim\dimen@>\tw@\wdz@\llap{\boxz@}\else
 \llap{\vtop{\normalbaselines\null\boxz@}}\fi
 \tabskip\z@skip\crcr#1\crcr\black@\gmaxwidth@}}
\newdimen\glineht@
\expandafter\def\csname gather \endcsname#1\endgather{\gmeasure@#1\endgather
 \ifdim\gmaxwidth@>\displaywidth\let\gdisplaywidth@\gmaxwidth@\else
 \let\gdisplaywidth@\displaywidth\fi\tabskip\centering@\halign to\displaywidth
 {\hfil\strut@\setboxz@h{$\m@th\displaystyle{##}$}%
 \global\gwidth@\wdz@\global\glineht@\ht\z@\boxz@\hfil&\kern-\gdisplaywidth@
 \setboxz@h{\strut@{\maketag@##\maketag@}}%
 \dimen@\displaywidth\advance\dimen@-\gwidth@
 \ifdim\dimen@>\tw@\wdz@\rlap{\boxz@}\else
 \rlap{\vbox{\normalbaselines\boxz@\vbox to\glineht@{}}}\fi
 \tabskip\gdisplaywidth@\crcr#1\crcr\black@\gmaxwidth@}}
\newif\ifctagsplit@
\def\CenteredTagsOnSplits{\global\ctagsplit@true}
\def\TopOrBottomTagsOnSplits{\global\ctagsplit@false}
\TopOrBottomTagsOnSplits
\def\split{\relax\ifinany@\let\next@\insplit@\else
 \ifmmode\ifinner\def\next@{\onlydmatherr@\split}\else
 \let\next@\outsplit@\fi\else
 \def\next@{\onlydmatherr@\split}\fi\fi\next@}
\def\insplit@{\global\setbox\z@\vbox\bgroup\vspace@\Let@\ialign\bgroup
 \hfil\strut@$\m@th\displaystyle{##}$&$\m@th\displaystyle{{}##}$\hfill\crcr}
\def\endsplit{\crcr\egroup\egroup\iftagsleft@\expandafter\lendsplit@\else
 \expandafter\rendsplit@\fi}
\def\rendsplit@{\global\setbox9 \vbox
 {\unvcopy\z@\global\setbox8 \lastbox\unskip}
 \setbox@ne\hbox{\unhcopy8 \unskip\global\setbox\tw@\lastbox
 \unskip\global\setbox\thr@@\lastbox}
 \global\setbox7 \hbox{\unhbox\tw@\unskip}
 \ifinalign@\ifctagsplit@                                                   
  \gdef\split@{\hbox to\wd\thr@@{}&
   \vcenter{\vbox{\moveleft\wd\thr@@\boxz@}}}
 \else\gdef\split@{&\vbox{\moveleft\wd\thr@@\box9}\crcr
  \box\thr@@&\box7}\fi                                                      
 \else                                                                      
  \ifctagsplit@\gdef\split@{\vcenter{\boxz@}}\else
  \gdef\split@{\box9\crcr\hbox{\box\thr@@\box7}}\fi
 \fi
 \split@}                                                                   
\def\lendsplit@{\global\setbox9\vtop{\unvcopy\z@}
 \setbox@ne\vbox{\unvcopy\z@\global\setbox8\lastbox}
 \setbox@ne\hbox{\unhcopy8\unskip\setbox\tw@\lastbox
  \unskip\global\setbox\thr@@\lastbox}
 \ifinalign@\ifctagsplit@                                                   
  \gdef\split@{\hbox to\wd\thr@@{}&
  \vcenter{\vbox{\moveleft\wd\thr@@\box9}}}
  \else                                                                     
  \gdef\split@{\hbox to\wd\thr@@{}&\vbox{\moveleft\wd\thr@@\box9}}\fi
 \else
  \ifctagsplit@\gdef\split@{\vcenter{\box9}}\else
  \gdef\split@{\box9}\fi
 \fi\split@}
\def\outsplit@#1$${\align\insplit@#1\endalign$$}
\newdimen\multlinegap@
\multlinegap@1em
\newdimen\multlinetaggap@
\multlinetaggap@1em
\def\MultlineGap#1{\global\multlinegap@#1\relax}
\def\multlinegap#1{\RIfMIfI@\onlydmatherr@\multlinegap\else
 \multlinegap@#1\relax\fi\else\onlydmatherr@\multlinegap\fi}
\def\nomultlinegap{\multlinegap{\z@}}
\def\multline{\RIfMIfI@
 \DN@{\onlydmatherr@\multline}\else
 \DN@{\multline@}\fi\else
 \DN@{\onlydmatherr@\multline}\fi\next@}
\newif\iftagin@
\def\tagin@#1{\tagin@false\in@\tag{#1}\ifin@\tagin@true\fi}
\def\multline@#1$${\inany@true\vspace@\allowdisplaybreak@\displaybreak@
 \tagin@{#1}\iftagsleft@\DN@{\multline@l#1$$}\else
 \DN@{\multline@r#1$$}\fi\next@}
\newdimen\mwidth@
\def\rmmeasure@#1\endmultline{%
 \def\shoveleft##1{##1}\def\shoveright##1{##1}
 \setbox@ne\vbox{\Let@\halign{\setboxz@h
  {$\m@th\@lign\displaystyle{}##$}\global\mwidth@\wdz@
  \crcr#1\crcr}}}
\newdimen\mlineht@
\newif\ifzerocr@
\newif\ifonecr@
\def\lmmeasure@#1\endmultline{\global\zerocr@true\global\onecr@false
 \everycr{\noalign{\ifonecr@\global\onecr@false\fi
  \ifzerocr@\global\zerocr@false\global\onecr@true\fi}}
  \def\shoveleft##1{##1}\def\shoveright##1{##1}%
 \setbox@ne\vbox{\Let@\halign{\setboxz@h
  {$\m@th\@lign\displaystyle{}##$}\ifonecr@\global\mwidth@\wdz@
  \global\mlineht@\ht\z@\fi\crcr#1\crcr}}}
\newbox\mtagbox@
\newdimen\ltwidth@
\newdimen\rtwidth@
\def\multline@l#1$${\iftagin@\DN@{\lmultline@@#1$$}\else
 \DN@{\setbox\mtagbox@\null\ltwidth@\z@\rtwidth@\z@
  \lmultline@@@#1$$}\fi\next@}
\def\lmultline@@#1\endmultline\tag#2$${%
 \setbox\mtagbox@\hbox{\maketag@#2\maketag@}
 \lmmeasure@#1\endmultline\dimen@\mwidth@\advance\dimen@\wd\mtagbox@
 \advance\dimen@\multlinetaggap@                                            
 \ifdim\dimen@>\displaywidth\ltwidth@\z@\else\ltwidth@\wd\mtagbox@\fi       
 \lmultline@@@#1\endmultline$$}
\def\lmultline@@@{\displ@y
 \def\shoveright##1{##1\hfilneg\hskip\multlinegap@}%
 \def\shoveleft##1{\setboxz@h{$\m@th\displaystyle{}##1$}%
  \setbox@ne\hbox{$\m@th\displaystyle##1$}%
  \hfilneg
  \iftagin@
   \ifdim\ltwidth@>\z@\hskip\ltwidth@\hskip\multlinetaggap@\fi
  \else\hskip\multlinegap@\fi\hskip.5\wd@ne\hskip-.5\wdz@##1}
  \halign\bgroup\Let@\hbox to\displaywidth
   {\strut@$\m@th\displaystyle\hfil{}##\hfil$}\crcr
   \hfilneg                                                                 
   \iftagin@                                                                
    \ifdim\ltwidth@>\z@                                                     
     \box\mtagbox@\hskip\multlinetaggap@                                    
    \else
     \rlap{\vbox{\normalbaselines\hbox{\strut@\box\mtagbox@}%
     \vbox to\mlineht@{}}}\fi                                               
   \else\hskip\multlinegap@\fi}                                             
\def\multline@r#1$${\iftagin@\DN@{\rmultline@@#1$$}\else
 \DN@{\setbox\mtagbox@\null\ltwidth@\z@\rtwidth@\z@
  \rmultline@@@#1$$}\fi\next@}
\def\rmultline@@#1\endmultline\tag#2$${\ltwidth@\z@
 \setbox\mtagbox@\hbox{\maketag@#2\maketag@}%
 \rmmeasure@#1\endmultline\dimen@\mwidth@\advance\dimen@\wd\mtagbox@
 \advance\dimen@\multlinetaggap@
 \ifdim\dimen@>\displaywidth\rtwidth@\z@\else\rtwidth@\wd\mtagbox@\fi
 \rmultline@@@#1\endmultline$$}
\def\rmultline@@@{\displ@y
 \def\shoveright##1{##1\hfilneg\iftagin@\ifdim\rtwidth@>\z@
  \hskip\rtwidth@\hskip\multlinetaggap@\fi\else\hskip\multlinegap@\fi}%
 \def\shoveleft##1{\setboxz@h{$\m@th\displaystyle{}##1$}%
  \setbox@ne\hbox{$\m@th\displaystyle##1$}%
  \hfilneg\hskip\multlinegap@\hskip.5\wd@ne\hskip-.5\wdz@##1}%
 \halign\bgroup\Let@\hbox to\displaywidth
  {\strut@$\m@th\displaystyle\hfil{}##\hfil$}\crcr
 \hfilneg\hskip\multlinegap@}
\def\endmultline{\iftagsleft@\expandafter\lendmultline@\else
 \expandafter\rendmultline@\fi}
\def\lendmultline@{\hfilneg\hskip\multlinegap@\crcr\egroup}
\def\rendmultline@{\iftagin@                                                
 \ifdim\rtwidth@>\z@                                                        
  \hskip\multlinetaggap@\box\mtagbox@                                       
 \else\llap{\vtop{\normalbaselines\null\hbox{\strut@\box\mtagbox@}}}\fi     
 \else\hskip\multlinegap@\fi                                                
 \hfilneg\crcr\egroup}
\def\bmod{\mskip-\medmuskip\mkern5mu\mathbin{\fam\z@ mod}\penalty900
 \mkern5mu\mskip-\medmuskip}
\def\pmod#1{\allowbreak\ifinner\mkern8mu\else\mkern18mu\fi
 ({\fam\z@ mod}\,\,#1)}
\def\pod#1{\allowbreak\ifinner\mkern8mu\else\mkern18mu\fi(#1)}
\def\mod#1{\allowbreak\ifinner\mkern12mu\else\mkern18mu\fi{\fam\z@ mod}\,\,#1}
\message{continued fractions,}
\newcount\cfraccount@
\def\cfrac{\bgroup\bgroup\advance\cfraccount@\@ne\strut
 \iffalse{\fi\def\\{\over\displaystyle}\iffalse}\fi}
\def\lcfrac{\bgroup\bgroup\advance\cfraccount@\@ne\strut
 \iffalse{\fi\def\\{\hfill\over\displaystyle}\iffalse}\fi}
\def\rcfrac{\bgroup\bgroup\advance\cfraccount@\@ne\strut\hfill
 \iffalse{\fi\def\\{\over\displaystyle}\iffalse}\fi}
\def\gloop@#1\repeat{\gdef\body{#1}\iterate}
\def\endcfrac{\gloop@\ifnum\cfraccount@>\z@\global\advance\cfraccount@\m@ne
 \egroup\hskip-\nulldelimiterspace\egroup\repeat}
\message{compound symbols,}
\def\binrel@#1{\setboxz@h{\thinmuskip0mu
  \medmuskip\m@ne mu\thickmuskip\@ne mu$#1\m@th$}%
 \setbox@ne\hbox{\thinmuskip0mu\medmuskip\m@ne mu\thickmuskip
  \@ne mu${}#1{}\m@th$}%
 \setbox\tw@\hbox{\hskip\wd@ne\hskip-\wdz@}}
\def\overset#1\to#2{\binrel@{#2}\ifdim\wd\tw@<\z@
 \mathbin{\mathop{\kern\z@#2}\limits^{#1}}\else\ifdim\wd\tw@>\z@
 \mathrel{\mathop{\kern\z@#2}\limits^{#1}}\else
 {\mathop{\kern\z@#2}\limits^{#1}}{}\fi\fi}
\def\underset#1\to#2{\binrel@{#2}\ifdim\wd\tw@<\z@
 \mathbin{\mathop{\kern\z@#2}\limits_{#1}}\else\ifdim\wd\tw@>\z@
 \mathrel{\mathop{\kern\z@#2}\limits_{#1}}\else
 {\mathop{\kern\z@#2}\limits_{#1}}{}\fi\fi}
\def\oversetbrace#1\to#2{\overbrace{#2}^{#1}}
\def\undersetbrace#1\to#2{\underbrace{#2}_{#1}}
\def\sideset#1\and#2\to#3{%
 \setbox@ne\hbox{$\dsize{\vphantom{#3}}#1{#3}\m@th$}%
 \setbox\tw@\hbox{$\dsize{#3}#2\m@th$}%
 \hskip\wd@ne\hskip-\wd\tw@\mathop{\hskip\wd\tw@\hskip-\wd@ne
  {\vphantom{#3}}#1{#3}#2}}
\def\rightarrowfill@#1{\setboxz@h{$#1-\m@th$}\ht\z@\z@
  $#1\m@th\copy\z@\mkern-6mu\cleaders
  \hbox{$#1\mkern-2mu\box\z@\mkern-2mu$}\hfill
  \mkern-6mu\mathord\rightarrow$}
\def\leftarrowfill@#1{\setboxz@h{$#1-\m@th$}\ht\z@\z@
  $#1\m@th\mathord\leftarrow\mkern-6mu\cleaders
  \hbox{$#1\mkern-2mu\copy\z@\mkern-2mu$}\hfill
  \mkern-6mu\box\z@$}
\def\leftrightarrowfill@#1{\setboxz@h{$#1-\m@th$}\ht\z@\z@
  $#1\m@th\mathord\leftarrow\mkern-6mu\cleaders
  \hbox{$#1\mkern-2mu\box\z@\mkern-2mu$}\hfill
  \mkern-6mu\mathord\rightarrow$}
\def\overrightarrow{\mathpalette\overrightarrow@}
\def\overrightarrow@#1#2{\vbox{\ialign{##\crcr\rightarrowfill@#1\crcr
 \noalign{\kern-\ex@\nointerlineskip}$\m@th\hfil#1#2\hfil$\crcr}}}

\def\overleftarrow{\mathpalette\overleftarrow@}
\def\overleftarrow@#1#2{\vbox{\ialign{##\crcr\leftarrowfill@#1\crcr
 \noalign{\kern-\ex@\nointerlineskip}$\m@th\hfil#1#2\hfil$\crcr}}}
\def\overleftrightarrow{\mathpalette\overleftrightarrow@}
\def\overleftrightarrow@#1#2{\vbox{\ialign{##\crcr\leftrightarrowfill@#1\crcr
 \noalign{\kern-\ex@\nointerlineskip}$\m@th\hfil#1#2\hfil$\crcr}}}
\def\underrightarrow{\mathpalette\underrightarrow@}
\def\underrightarrow@#1#2{\vtop{\ialign{##\crcr$\m@th\hfil#1#2\hfil$\crcr
 \noalign{\nointerlineskip}\rightarrowfill@#1\crcr}}}

\def\underleftarrow{\mathpalette\underleftarrow@}
\def\underleftarrow@#1#2{\vtop{\ialign{##\crcr$\m@th\hfil#1#2\hfil$\crcr
 \noalign{\nointerlineskip}\leftarrowfill@#1\crcr}}}
\def\underleftrightarrow{\mathpalette\underleftrightarrow@}
\def\underleftrightarrow@#1#2{\vtop{\ialign{##\crcr$\m@th\hfil#1#2\hfil$\crcr
 \noalign{\nointerlineskip}\leftrightarrowfill@#1\crcr}}}
\message{various kinds of dots,}
\let\DOTSI\relax
\let\DOTSB\relax

\newif\ifmath@
{\uccode`7=`\\ \uccode`8=`m \uccode`9=`a \uccode`0=`t \uccode`!=`h
 \uppercase{\gdef\math@#1#2#3#4#5#6\math@{\global\math@false\ifx 7#1\ifx 8#2%
 \ifx 9#3\ifx 0#4\ifx !#5\xdef\meaning@{#6}\global\math@true\fi\fi\fi\fi\fi}}}
\newif\ifmathch@
{\uccode`7=`c \uccode`8=`h \uccode`9=`\"
 \uppercase{\gdef\mathch@#1#2#3#4#5#6\mathch@{\global\mathch@false
  \ifx 7#1\ifx 8#2\ifx 9#5\global\mathch@true\xdef\meaning@{9#6}\fi\fi\fi}}}
\newcount\classnum@
\def\getmathch@#1.#2\getmathch@{\classnum@#1 \divide\classnum@4096
 \ifcase\number\classnum@\or\or\gdef\thedots@{\dotsb@}\or
 \gdef\thedots@{\dotsb@}\fi}
\newif\ifmathbin@
{\uccode`4=`b \uccode`5=`i \uccode`6=`n
 \uppercase{\gdef\mathbin@#1#2#3{\relaxnext@
  \DNii@##1\mathbin@{\ifx\space@\next\global\mathbin@true\fi}%
 \global\mathbin@false\DN@##1\mathbin@{}%
 \ifx 4#1\ifx 5#2\ifx 6#3\DN@{\FN@\nextii@}\fi\fi\fi\next@}}}
\newif\ifmathrel@
{\uccode`4=`r \uccode`5=`e \uccode`6=`l
 \uppercase{\gdef\mathrel@#1#2#3{\relaxnext@
  \DNii@##1\mathrel@{\ifx\space@\next\global\mathrel@true\fi}%
 \global\mathrel@false\DN@##1\mathrel@{}%
 \ifx 4#1\ifx 5#2\ifx 6#3\DN@{\FN@\nextii@}\fi\fi\fi\next@}}}
\newif\ifmacro@
{\uccode`5=`m \uccode`6=`a \uccode`7=`c
 \uppercase{\gdef\macro@#1#2#3#4\macro@{\global\macro@false
  \ifx 5#1\ifx 6#2\ifx 7#3\global\macro@true
  \xdef\meaning@{\macro@@#4\macro@@}\fi\fi\fi}}}
\def\macro@@#1->#2\macro@@{#2}
\newif\ifDOTS@
\newcount\DOTSCASE@
{\uccode`6=`\\ \uccode`7=`D \uccode`8=`O \uccode`9=`T \uccode`0=`S
 \uppercase{\gdef\DOTS@#1#2#3#4#5{\global\DOTS@false\DN@##1\DOTS@{}%
  \ifx 6#1\ifx 7#2\ifx 8#3\ifx 9#4\ifx 0#5\let\next@\DOTS@@\fi\fi\fi\fi\fi
  \next@}}}
{\uccode`3=`B \uccode`4=`I \uccode`5=`X
 \uppercase{\gdef\DOTS@@#1{\relaxnext@
  \DNii@##1\DOTS@{\ifx\space@\next\global\DOTS@true\fi}%
  \DN@{\FN@\nextii@}%
  \ifx 3#1\global\DOTSCASE@\z@\else
  \ifx 4#1\global\DOTSCASE@\@ne\else
  \ifx 5#1\global\DOTSCASE@\tw@\else\DN@##1\DOTS@{}%
  \fi\fi\fi\next@}}}
\newif\ifnot@
{\uccode`5=`\\ \uccode`6=`n \uccode`7=`o \uccode`8=`t
 \uppercase{\gdef\not@#1#2#3#4{\relaxnext@
  \DNii@##1\not@{\ifx\space@\next\global\not@true\fi}%
 \global\not@false\DN@##1\not@{}%
 \ifx 5#1\ifx 6#2\ifx 7#3\ifx 8#4\DN@{\FN@\nextii@}\fi\fi\fi
 \fi\next@}}}
\newif\ifkeybin@
\def\keybin@{\keybin@true
 \ifx\next+\else\ifx\next=\else\ifx\next<\else\ifx\next>\else\ifx\next-\else
 \ifx\next*\else\ifx\next:\else\keybin@false\fi\fi\fi\fi\fi\fi\fi}
\def\dots{\RIfM@\expandafter\mdots@\else\expandafter\tdots@\fi}
\def\tdots@{\unskip\relaxnext@
 \DN@{$\m@th\mathinner{\ldotp\ldotp\ldotp}\,
   \ifx\next,\,$\else\ifx\next.\,$\else\ifx\next;\,$\else\ifx\next:\,$\else
   \ifx\next?\,$\else\ifx\next!\,$\else$ \fi\fi\fi\fi\fi\fi}%
 \ \FN@\next@}
\def\mdots@{\FN@\mdots@@}
\def\mdots@@{\gdef\thedots@{\dotso@}
 \ifx\next\boldkey\gdef\thedots@\boldkey{\boldkeydots@}\else                
 \ifx\next\boldsymbol\gdef\thedots@\boldsymbol{\boldsymboldots@}\else       
 \ifx,\next\gdef\thedots@{\dotsc}
 \else\ifx\not\next\gdef\thedots@{\dotsb@}
 \else\keybin@
 \ifkeybin@\gdef\thedots@{\dotsb@}
 \else\xdef\meaning@{\meaning\next..........}\xdef\meaning@@{\meaning@}
  \expandafter\math@\meaning@\math@
  \ifmath@
   \expandafter\mathch@\meaning@\mathch@
   \ifmathch@\expandafter\getmathch@\meaning@\getmathch@\fi                 
  \else\expandafter\macro@\meaning@@\macro@                                 
  \ifmacro@                                                                
   \expandafter\not@\meaning@\not@\ifnot@\gdef\thedots@{\dotsb@}
  \else\expandafter\DOTS@\meaning@\DOTS@
  \ifDOTS@
   \ifcase\number\DOTSCASE@\gdef\thedots@{\dotsb@}%
    \or\gdef\thedots@{\dotsi}\else\fi                                      
  \else\expandafter\math@\meaning@\math@                                   
  \ifmath@\expandafter\mathbin@\meaning@\mathbin@
  \ifmathbin@\gdef\thedots@{\dotsb@}
  \else\expandafter\mathrel@\meaning@\mathrel@
  \ifmathrel@\gdef\thedots@{\dotsb@}
  \fi\fi\fi\fi\fi\fi\fi\fi\fi\fi\fi\fi
 \thedots@}
\def\plainldots@{\mathinner{\ldotp\ldotp\ldotp}}
\def\plaincdots@{\mathinner{\cdotp\cdotp\cdotp}}
\def\dotsi{\!\plaincdots@}
\let\dotsb@\plaincdots@
\newif\ifextra@
\newif\ifrightdelim@
\def\rightdelim@{\global\rightdelim@true                                    
 \ifx\next)\else                                                            
 \ifx\next]\else
 \ifx\next\rbrack\else
 \ifx\next\}\else
 \ifx\next\rbrace\else
 \ifx\next\rangle\else
 \ifx\next\rceil\else
 \ifx\next\rfloor\else
 \ifx\next\rgroup\else
 \ifx\next\rmoustache\else
 \ifx\next\right\else
 \ifx\next\bigr\else
 \ifx\next\biggr\else
 \ifx\next\Bigr\else                                                        
 \ifx\next\Biggr\else\global\rightdelim@false
 \fi\fi\fi\fi\fi\fi\fi\fi\fi\fi\fi\fi\fi\fi\fi}
\def\extra@{%
 \global\extra@false\rightdelim@\ifrightdelim@\global\extra@true            
 \else\ifx\next$\global\extra@true                                          
 \else\xdef\meaning@{\meaning\next..........}
 \expandafter\macro@\meaning@\macro@\ifmacro@                               
 \expandafter\DOTS@\meaning@\DOTS@
 \ifDOTS@
 \ifnum\DOTSCASE@=\tw@\global\extra@true                                    
 \fi\fi\fi\fi\fi}
\newif\ifbold@
\def\dotso@{\relaxnext@
 \ifbold@
  \let\next\delayed@
  \DNii@{\extra@\plainldots@\ifextra@\,\fi}%
 \else
  \DNii@{\DN@{\extra@\plainldots@\ifextra@\,\fi}\FN@\next@}%
 \fi
 \nextii@}
\def\extrap@#1{%
 \ifx\next,\DN@{#1\,}\else
 \ifx\next;\DN@{#1\,}\else
 \ifx\next.\DN@{#1\,}\else\extra@
 \ifextra@\DN@{#1\,}\else
 \let\next@#1\fi\fi\fi\fi\next@}
\def\ldots{\DN@{\extrap@\plainldots@}%
 \FN@\next@}
\def\cdots{\DN@{\extrap@\plaincdots@}%
 \FN@\next@}

\def\dotsc{\relaxnext@
 \DN@{\ifx\next;\plainldots@\,\else
  \ifx\next.\plainldots@\,\else\extra@\plainldots@
  \ifextra@\,\fi\fi\fi}%
 \FN@\next@}
\def\cdot{\mathchar"2201 }

\message{special superscripts,}
\def\dddot#1{{\mathop{#1}\limits^{\vbox to-1.4\ex@{\kern-\tw@\ex@
 \hbox{\rm...}\vss}}}}
\def\ddddot#1{{\mathop{#1}\limits^{\vbox to-1.4\ex@{\kern-\tw@\ex@
 \hbox{\rm....}\vss}}}}
\def\sphat{^{\mathchoice{}{}%
 {\,\,\botsmash{\hbox{\lower4\ex@\hbox{$\m@th\widehat{\null}$}}}}%
 {\,\botsmash{\hbox{\lower3\ex@\hbox{$\m@th\hat{\null}$}}}}}}

\def\spacute{^{\!\botsmash{\hbox{\lower\@ne ex\hbox{\'{}}}}}}
\def\spgrave{^{\mathchoice{}{}{}{\!}%
 \botsmash{\hbox{\lower\@ne ex\hbox{\`{}}}}}}
\def\spdot{^{\hbox{\raise\ex@\hbox{\rm.}}}}
\def\spddot{^{\hbox{\raise\ex@\hbox{\rm..}}}}
\def\spdddot{^{\hbox{\raise\ex@\hbox{\rm...}}}}
\def\spddddot{^{\hbox{\raise\ex@\hbox{\rm....}}}}
\def\spbreve{^{\!\botsmash{\hbox{\lower4\ex@\hbox{\u{}}}}}}

\message{\string\text,}
\def\textonlyfont@#1#2{\def#1{\RIfM@
 \Err@{Use \string#1\space only in text}\else#2\fi}}
\textonlyfont@\rm\tenrm
\textonlyfont@\it\tenit
\textonlyfont@\sl\tensl
\textonlyfont@\bf\tenbf
\def\oldnos#1{\RIfM@{\mathcode`\,="013B \fam\@ne#1}\else
 \leavevmode\hbox{$\m@th\mathcode`\,="013B \fam\@ne#1$}\fi}
\def\text{\RIfM@\expandafter\text@\else\expandafter\text@@\fi}
\def\text@@#1{\leavevmode\hbox{#1}}
\def\mathhexbox@#1#2#3{\text{$\m@th\mathchar"#1#2#3$}}
\def\dag{{\mathhexbox@279}}
\def\ddag{{\mathhexbox@27A}}
\def\S{{\mathhexbox@278}}
\def\P{{\mathhexbox@27B}}
\newif\iffirstchoice@
\firstchoice@true
\def\text@#1{\mathchoice
 {\hbox{\everymath{\displaystyle}\def\textfonti{\the\textfont\@ne}%
  \def\textfontii{\the\textfont\tw@}\textdef@@ T#1}}
 {\hbox{\firstchoice@false
  \everymath{\textstyle}\def\textfonti{\the\textfont\@ne}%
  \def\textfontii{\the\textfont\tw@}\textdef@@ T#1}}
 {\hbox{\firstchoice@false
  \everymath{\scriptstyle}\def\textfonti{\the\scriptfont\@ne}%
  \def\textfontii{\the\scriptfont\tw@}\textdef@@ S\rm#1}}
 {\hbox{\firstchoice@false
  \everymath{\scriptscriptstyle}\def\textfonti
  {\the\scriptscriptfont\@ne}%
  \def\textfontii{\the\scriptscriptfont\tw@}\textdef@@ s\rm#1}}}
\def\textdef@@#1{\textdef@#1\rm\textdef@#1\bf\textdef@#1\sl\textdef@#1\it}
\def\rmfam{0}
\def\textdef@#1#2{%
 \DN@{\csname\expandafter\eat@\string#2fam\endcsname}%
 \if S#1\edef#2{\the\scriptfont\next@\relax}%
 \else\if s#1\edef#2{\the\scriptscriptfont\next@\relax}%
 \else\edef#2{\the\textfont\next@\relax}\fi\fi}
\scriptfont\itfam\tenit \scriptscriptfont\itfam\tenit
\scriptfont\slfam\tensl \scriptscriptfont\slfam\tensl
\newif\iftopfolded@
\newif\ifbotfolded@
\def\topfoldedtext{\topfolded@true\botfolded@false\foldedtext@}
\def\botfoldedtext{\botfolded@true\topfolded@false\foldedtext@}
\def\foldedtext{\topfolded@false\botfolded@false\foldedtext@}
\Invalid@\foldedwidth
\def\foldedtext@{\relaxnext@
 \DN@{\ifx\next\foldedwidth\let\next@\nextii@\else
  \DN@{\nextii@\foldedwidth{.3\hsize}}\fi\next@}%
 \DNii@\foldedwidth##1##2{\setbox\z@\vbox
  {\normalbaselines\hsize##1\relax
  \tolerance1600 \noindent\ignorespaces##2}\ifbotfolded@\boxz@\else
  \iftopfolded@\vtop{\unvbox\z@}\else\vcenter{\boxz@}\fi\fi}%
 \FN@\next@}
\message{math font commands,}
\def\bold{\RIfM@\expandafter\bold@\else
 \expandafter\nonmatherr@\expandafter\bold\fi}
\def\bold@#1{{\bold@@{#1}}}
\def\bold@@#1{\fam\bffam\relax#1}
\def\slanted{\RIfM@\expandafter\slanted@\else
 \expandafter\nonmatherr@\expandafter\slanted\fi}
\def\slanted@#1{{\slanted@@{#1}}}
\def\slanted@@#1{\fam\slfam\relax#1}
\def\roman{\RIfM@\expandafter\roman@\else
 \expandafter\nonmatherr@\expandafter\roman\fi}
\def\roman@#1{{\roman@@{#1}}}
\def\roman@@#1{\fam\rmfam\relax#1}
\def\italic{\RIfM@\expandafter\italic@\else
 \expandafter\nonmatherr@\expandafter\italic\fi}
\def\italic@#1{{\italic@@{#1}}}
\def\italic@@#1{\fam\itfam\relax#1}
\def\Cal{\RIfM@\expandafter\Cal@\else
 \expandafter\nonmatherr@\expandafter\Cal\fi}
\def\Cal@#1{{\Cal@@{#1}}}
\def\Cal@@#1{\noaccents@\fam\tw@#1}
\mathchardef\Gamma="0000
\mathchardef\Delta="0001
\mathchardef\Theta="0002
\mathchardef\Lambda="0003
\mathchardef\Xi="0004
\mathchardef\Pi="0005
\mathchardef\Sigma="0006
\mathchardef\Upsilon="0007
\mathchardef\Phi="0008
\mathchardef\Psi="0009
\mathchardef\Omega="000A
\mathchardef\varGamma="0100
\mathchardef\varDelta="0101
\mathchardef\varTheta="0102
\mathchardef\varLambda="0103
\mathchardef\varXi="0104
\mathchardef\varPi="0105
\mathchardef\varSigma="0106
\mathchardef\varUpsilon="0107
\mathchardef\varPhi="0108
\mathchardef\varPsi="0109
\mathchardef\varOmega="010A
\let\alloc@@\alloc@
\def\hexnumber@#1{\ifcase#1 0\or 1\or 2\or 3\or 4\or 5\or 6\or 7\or 8\or
 9\or A\or B\or C\or D\or E\or F\fi}
\def\loadmsam{%
 \font@\tenmsa=msam10
 \font@\sevenmsa=msam7
 \font@\fivemsa=msam5
 \alloc@@8\fam\chardef\sixt@@n\msafam
 \textfont\msafam=\tenmsa
 \scriptfont\msafam=\sevenmsa
 \scriptscriptfont\msafam=\fivemsa
 \edef\next{\hexnumber@\msafam}%
 \mathchardef\dabar@"0\next39
 \edef\dashrightarrow{\mathrel{\dabar@\dabar@\mathchar"0\next4B}}%
 \edef\dashleftarrow{\mathrel{\mathchar"0\next4C\dabar@\dabar@}}%
 \let\dasharrow\dashrightarrow
 \edef\ulcorner{\delimiter"4\next70\next70 }%
 \edef\urcorner{\delimiter"5\next71\next71 }%
 \edef\llcorner{\delimiter"4\next78\next78 }%
 \edef\lrcorner{\delimiter"5\next79\next79 }%
 \edef\yen{{\noexpand\mathhexbox@\next55}}%
 \edef\checkmark{{\noexpand\mathhexbox@\next58}}%
 \edef\circledR{{\noexpand\mathhexbox@\next72}}%
 \edef\maltese{{\noexpand\mathhexbox@\next7A}}%
 \global\let\loadmsam\empty}%
\def\loadmsbm{%
 \font@\tenmsb=msbm10 \font@\sevenmsb=msbm7 \font@\fivemsb=msbm5
 \alloc@@8\fam\chardef\sixt@@n\msbfam
 \textfont\msbfam=\tenmsb
 \scriptfont\msbfam=\sevenmsb \scriptscriptfont\msbfam=\fivemsb
 \global\let\loadmsbm\empty
 }
\def\widehat#1{\ifx\undefined\msbfam \DN@{362}%
  \else \setboxz@h{$\m@th#1$}%
    \edef\next@{\ifdim\wdz@>\tw@ em%
        \hexnumber@\msbfam 5B%
      \else 362\fi}\fi
  \mathaccent"0\next@{#1}}
\def\widetilde#1{\ifx\undefined\msbfam \DN@{365}%
  \else \setboxz@h{$\m@th#1$}%
    \edef\next@{\ifdim\wdz@>\tw@ em%
        \hexnumber@\msbfam 5D%
      \else 365\fi}\fi
  \mathaccent"0\next@{#1}}
\message{\string\newsymbol,}
\def\newsymbol#1#2#3#4#5{\define#1{}%
  \count@#2\relax \advance\count@\m@ne 
 \ifcase\count@
   \ifx\undefined\msafam\loadmsam\fi \let\next@\msafam
 \or \ifx\undefined\msbfam\loadmsbm\fi \let\next@\msbfam
 \else  \Err@{\Invalid@@\string\newsymbol}\let\next@\tw@\fi
 \mathchardef#1="#3\hexnumber@\next@#4#5\space}

\def\Bbb{\RIfM@\expandafter\Bbb@\else
 \expandafter\nonmatherr@\expandafter\Bbb\fi}
\def\Bbb@#1{{\Bbb@@{#1}}}
\def\Bbb@@#1{\noaccents@\fam\msbfam\relax#1}
\message{bold Greek and bold symbols,}
\def\loadbold{%
 \font@\tencmmib=cmmib10 \font@\sevencmmib=cmmib7 \font@\fivecmmib=cmmib5
 \skewchar\tencmmib'177 \skewchar\sevencmmib'177 \skewchar\fivecmmib'177
 \alloc@@8\fam\chardef\sixt@@n\cmmibfam
 \textfont\cmmibfam\tencmmib
 \scriptfont\cmmibfam\sevencmmib \scriptscriptfont\cmmibfam\fivecmmib
 \font@\tencmbsy=cmbsy10 \font@\sevencmbsy=cmbsy7 \font@\fivecmbsy=cmbsy5
 \skewchar\tencmbsy'60 \skewchar\sevencmbsy'60 \skewchar\fivecmbsy'60
 \alloc@@8\fam\chardef\sixt@@n\cmbsyfam
 \textfont\cmbsyfam\tencmbsy
 \scriptfont\cmbsyfam\sevencmbsy \scriptscriptfont\cmbsyfam\fivecmbsy
 \let\loadbold\empty
}
\def\boldnotloaded#1{\Err@{\ifcase#1\or First\else Second\fi
       bold symbol font not loaded}}
\def\mathchari@#1#2#3{\ifx\undefined\cmmibfam
    \boldnotloaded@\@ne
  \else\mathchar"#1\hexnumber@\cmmibfam#2#3\space \fi}
\def\mathcharii@#1#2#3{\ifx\undefined\cmbsyfam
    \boldnotloaded\tw@
  \else \mathchar"#1\hexnumber@\cmbsyfam#2#3\space\fi}
\edef\bffam@{\hexnumber@\bffam}
\def\boldkey#1{\ifcat\noexpand#1A%
  \ifx\undefined\cmmibfam \boldnotloaded\@ne
  \else {\fam\cmmibfam#1}\fi
 \else
 \ifx#1!\mathchar"5\bffam@21 \else
 \ifx#1(\mathchar"4\bffam@28 \else\ifx#1)\mathchar"5\bffam@29 \else
 \ifx#1+\mathchar"2\bffam@2B \else\ifx#1:\mathchar"3\bffam@3A \else
 \ifx#1;\mathchar"6\bffam@3B \else\ifx#1=\mathchar"3\bffam@3D \else
 \ifx#1?\mathchar"5\bffam@3F \else\ifx#1[\mathchar"4\bffam@5B \else
 \ifx#1]\mathchar"5\bffam@5D \else
 \ifx#1,\mathchari@63B \else
 \ifx#1-\mathcharii@200 \else
 \ifx#1.\mathchari@03A \else
 \ifx#1/\mathchari@03D \else
 \ifx#1<\mathchari@33C \else
 \ifx#1>\mathchari@33E \else
 \ifx#1*\mathcharii@203 \else
 \ifx#1|\mathcharii@06A \else
 \ifx#10\bold0\else\ifx#11\bold1\else\ifx#12\bold2\else\ifx#13\bold3\else
 \ifx#14\bold4\else\ifx#15\bold5\else\ifx#16\bold6\else\ifx#17\bold7\else
 \ifx#18\bold8\else\ifx#19\bold9\else
  \Err@{\string\boldkey\space can't be used with #1}%
 \fi\fi\fi\fi\fi\fi\fi\fi\fi\fi\fi\fi\fi\fi\fi
 \fi\fi\fi\fi\fi\fi\fi\fi\fi\fi\fi\fi\fi\fi}
\def\boldsymbol#1{%
 \DN@{\Err@{You can't use \string\boldsymbol\space with \string#1}#1}%
 \ifcat\noexpand#1A%
   \let\next@\relax
   \ifx\undefined\cmmibfam \boldnotloaded\@ne
   \else {\fam\cmmibfam#1}\fi
 \else
  \xdef\meaning@{\meaning#1.........}%
  \expandafter\math@\meaning@\math@
  \ifmath@
   \expandafter\mathch@\meaning@\mathch@
   \ifmathch@
    \expandafter\boldsymbol@@\meaning@\boldsymbol@@
   \fi
  \else
   \expandafter\macro@\meaning@\macro@
   \expandafter\delim@\meaning@\delim@
   \ifdelim@
    \expandafter\delim@@\meaning@\delim@@
   \else
    \boldsymbol@{#1}%
   \fi
  \fi
 \fi
 \next@}
\def\mathhexboxii@#1#2{\ifx\undefined\cmbsyfam
    \boldnotloaded\tw@
  \else \mathhexbox@{\hexnumber@\cmbsyfam}{#1}{#2}\fi}
\def\boldsymbol@#1{\let\next@\relax\let\next#1%
 \ifx\next\cdot\mathcharii@201 \else
 \ifx\next\prime{{\null\mathcharii@030 \null}}\else
 \ifx\next\lbrack\mathchar"4\bffam@5B \else
 \ifx\next\rbrack\mathchar"5\bffam@5D \else
 \ifx\next\{\mathcharii@466 \else
 \ifx\next\lbrace\mathcharii@466 \else
 \ifx\next\}\mathcharii@567 \else
 \ifx\next\rbrace\mathcharii@567 \else
 \ifx\next\surd{{\mathcharii@170}}\else
 \ifx\next\S{{\mathhexboxii@78}}\else
 \ifx\next\P{{\mathhexboxii@7B}}\else
 \ifx\next\dag{{\mathhexboxii@79}}\else
 \ifx\next\ddag{{\mathhexboxii@7A}}\else
 \DN@{\Err@{You can't use \string\boldsymbol\space with \string#1}#1}%
 \fi\fi\fi\fi\fi\fi\fi\fi\fi\fi\fi\fi\fi}
\def\boldsymbol@@#1.#2\boldsymbol@@{\classnum@#1 \count@@@\classnum@        
 \divide\classnum@4096 \count@\classnum@                                    
 \multiply\count@4096 \advance\count@@@-\count@ \count@@\count@@@           
 \divide\count@@@\@cclvi \count@\count@@                                    
 \multiply\count@@@\@cclvi \advance\count@@-\count@@@                       
 \divide\count@@@\@cclvi                                                    
 \multiply\classnum@4096 \advance\classnum@\count@@                         
 \ifnum\count@@@=\z@                                                        
  \count@"\bffam@ \multiply\count@\@cclvi
  \advance\classnum@\count@
  \DN@{\mathchar\number\classnum@}%
 \else
  \ifnum\count@@@=\@ne                                                      
   \ifx\undefined\cmmibfam \DN@{\boldnotloaded\@ne}%
   \else \count@\cmmibfam \multiply\count@\@cclvi
     \advance\classnum@\count@
     \DN@{\mathchar\number\classnum@}\fi
  \else
   \ifnum\count@@@=\tw@                                                    
     \ifx\undefined\cmbsyfam
       \DN@{\boldnotloaded\tw@}%
     \else
       \count@\cmbsyfam \multiply\count@\@cclvi
       \advance\classnum@\count@
       \DN@{\mathchar\number\classnum@}%
     \fi
  \fi
 \fi
\fi}
\newif\ifdelim@
\newcount\delimcount@
{\uccode`6=`\\ \uccode`7=`d \uccode`8=`e \uccode`9=`l
 \uppercase{\gdef\delim@#1#2#3#4#5\delim@
  {\delim@false\ifx 6#1\ifx 7#2\ifx 8#3\ifx 9#4\delim@true
   \xdef\meaning@{#5}\fi\fi\fi\fi}}}
\def\delim@@#1"#2#3#4#5#6\delim@@{\if#32%
\let\next@\relax
 \ifx\undefined\cmbsyfam \boldnotloaded\@ne
 \else \mathcharii@#2#4#5\space \fi\fi}
\def\vert{\delimiter"026A30C }
\def\Vert{\delimiter"026B30D }
\let\|\Vert
\def\backslash{\delimiter"026E30F }
\def\boldkeydots@#1{\bold@true\let\next=#1\let\delayed@=#1\mdots@@
 \boldkey#1\bold@false}  
\def\boldsymboldots@#1{\bold@true\let\next#1\let\delayed@#1\mdots@@
 \boldsymbol#1\bold@false}
\message{Euler fonts,}

\def\frak{\mathfont@\frak}

\def\loadmathfont#1{%
   \expandafter\font@\csname ten#1\endcsname=#110
   \expandafter\font@\csname seven#1\endcsname=#17
   \expandafter\font@\csname five#1\endcsname=#15
   \edef\next{\noexpand\alloc@@8\fam\chardef\sixt@@n
     \expandafter\noexpand\csname#1fam\endcsname}%
   \next
   \textfont\csname#1fam\endcsname \csname ten#1\endcsname
   \scriptfont\csname#1fam\endcsname \csname seven#1\endcsname
   \scriptscriptfont\csname#1fam\endcsname \csname five#1\endcsname
   \expandafter\def\csname #1\expandafter\endcsname\expandafter{%
      \expandafter\mathfont@\csname#1\endcsname}%
 \expandafter\gdef\csname load#1\endcsname{}%
}
\def\mathfont@#1{\RIfM@\expandafter\mathfont@@\expandafter#1\else
  \expandafter\nonmatherr@\expandafter#1\fi}
\def\mathfont@@#1#2{{\mathfont@@@#1{#2}}}
\def\mathfont@@@#1#2{\noaccents@
   \fam\csname\expandafter\eat@\string#1fam\endcsname
   \relax#2}
\message{math accents,}
\def\accentclass@{7}
\def\noaccents@{\def\accentclass@{0}}
\def\makeacc@#1#2{\def#1{\mathaccent"\accentclass@#2 }}
\makeacc@\hat{05E}
\makeacc@\check{014}
\makeacc@\tilde{07E}
\makeacc@\acute{013}
\makeacc@\grave{012}
\makeacc@\dot{05F}
\makeacc@\ddot{07F}
\makeacc@\breve{015}
\makeacc@\bar{016}

\newcount\skewcharcount@
\newcount\familycount@
\def\theskewchar@{\familycount@\@ne
 \global\skewcharcount@\the\skewchar\textfont\@ne                           
 \ifnum\fam>\m@ne\ifnum\fam<16
  \global\familycount@\the\fam\relax
  \global\skewcharcount@\the\skewchar\textfont\the\fam\relax\fi\fi          
 \ifnum\skewcharcount@>\m@ne
  \ifnum\skewcharcount@<128
  \multiply\familycount@256
  \global\advance\skewcharcount@\familycount@
  \global\advance\skewcharcount@28672
  \mathchar\skewcharcount@\else
  \global\skewcharcount@\m@ne\fi\else
 \global\skewcharcount@\m@ne\fi}                                            
\newcount\pointcount@
\def\getpoints@#1.#2\getpoints@{\pointcount@#1 }
\newdimen\accentdimen@
\newcount\accentmu@
\def\dimentomu@{\multiply\accentdimen@ 100
 \expandafter\getpoints@\the\accentdimen@\getpoints@
 \multiply\pointcount@18
 \divide\pointcount@\@m
 \global\accentmu@\pointcount@}
\def\Makeacc@#1#2{\def#1{\RIfM@\DN@{\mathaccent@
 {"\accentclass@#2 }}\else\DN@{\nonmatherr@{#1}}\fi\next@}}
\def\unbracefonts@{\let\Cal@\Cal@@\let\roman@\roman@@\let\bold@\bold@@
 \let\slanted@\slanted@@}
\def\mathaccent@#1#2{\ifnum\fam=\m@ne\xdef\thefam@{1}\else
 \xdef\thefam@{\the\fam}\fi                                                 
 \accentdimen@\z@                                                           
 \setboxz@h{\unbracefonts@$\m@th\fam\thefam@\relax#2$}
 \ifdim\accentdimen@=\z@\DN@{\mathaccent#1{#2}}
  \setbox@ne\hbox{\unbracefonts@$\m@th\fam\thefam@\relax#2\theskewchar@$}
  \setbox\tw@\hbox{$\m@th\ifnum\skewcharcount@=\m@ne\else
   \mathchar\skewcharcount@\fi$}
  \global\accentdimen@\wd@ne\global\advance\accentdimen@-\wdz@
  \global\advance\accentdimen@-\wd\tw@                                     
  \global\multiply\accentdimen@\tw@
  \dimentomu@\global\advance\accentmu@\@ne                                 
 \else\DN@{{\mathaccent#1{#2\mkern\accentmu@ mu}%
    \mkern-\accentmu@ mu}{}}\fi                                             
 \next@}\Makeacc@\Hat{05E}
\Makeacc@\Check{014}
\Makeacc@\Tilde{07E}
\Makeacc@\Acute{013}
\Makeacc@\Grave{012}
\Makeacc@\Dot{05F}
\Makeacc@\Ddot{07F}
\Makeacc@\Breve{015}
\Makeacc@\Bar{016}
\def\Vec{\RIfM@\DN@{\mathaccent@{"017E }}\else
 \DN@{\nonmatherr@\Vec}\fi\next@}
\def\accentedsymbol#1#2{\csname newbox\expandafter\endcsname
  \csname\expandafter\eat@\string#1@box\endcsname
 \expandafter\setbox\csname\expandafter\eat@
  \string#1@box\endcsname\hbox{$\m@th#2$}\define
  #1{\copy\csname\expandafter\eat@\string#1@box\endcsname{}}}
\message{roots,}
\def\sqrt#1{\radical"270370 {#1}}
\let\underline@\underline
\let\overline@\overline
\def\underline#1{\underline@{#1}}
\def\overline#1{\overline@{#1}}
\Invalid@\leftroot
\Invalid@\uproot
\newcount\uproot@
\newcount\leftroot@
\def\root{\relaxnext@
  \DN@{\ifx\next\uproot\let\next@\nextii@\else
   \ifx\next\leftroot\let\next@\nextiii@\else
   \let\next@\plainroot@\fi\fi\next@}%
  \DNii@\uproot##1{\uproot@##1\relax\FN@\nextiv@}%
  \def\nextiv@{\ifx\next\space@\DN@. {\FN@\nextv@}\else
   \DN@.{\FN@\nextv@}\fi\next@.}%
  \def\nextv@{\ifx\next\leftroot\let\next@\nextvi@\else
   \let\next@\plainroot@\fi\next@}%
  \def\nextvi@\leftroot##1{\leftroot@##1\relax\plainroot@}%
   \def\nextiii@\leftroot##1{\leftroot@##1\relax\FN@\nextvii@}%
  \def\nextvii@{\ifx\next\space@
   \DN@. {\FN@\nextviii@}\else
   \DN@.{\FN@\nextviii@}\fi\next@.}%
  \def\nextviii@{\ifx\next\uproot\let\next@\nextix@\else
   \let\next@\plainroot@\fi\next@}%
  \def\nextix@\uproot##1{\uproot@##1\relax\plainroot@}%
  \bgroup\uproot@\z@\leftroot@\z@\FN@\next@}
\def\plainroot@#1\of#2{\setbox\rootbox\hbox{$\m@th\scriptscriptstyle{#1}$}%
 \mathchoice{\r@@t\displaystyle{#2}}{\r@@t\textstyle{#2}}
 {\r@@t\scriptstyle{#2}}{\r@@t\scriptscriptstyle{#2}}\egroup}
\def\r@@t#1#2{\setboxz@h{$\m@th#1\sqrt{#2}$}%
 \dimen@\ht\z@\advance\dimen@-\dp\z@
 \setbox@ne\hbox{$\m@th#1\mskip\uproot@ mu$}\advance\dimen@ 1.667\wd@ne
 \mkern-\leftroot@ mu\mkern5mu\raise.6\dimen@\copy\rootbox
 \mkern-10mu\mkern\leftroot@ mu\boxz@}
\def\boxed#1{\setboxz@h{$\m@th\displaystyle{#1}$}\dimen@.4\ex@
 \advance\dimen@3\ex@\advance\dimen@\dp\z@
 \hbox{\lower\dimen@\hbox{%
 \vbox{\hrule height.4\ex@
 \hbox{\vrule width.4\ex@\hskip3\ex@\vbox{\vskip3\ex@\boxz@\vskip3\ex@}%
 \hskip3\ex@\vrule width.4\ex@}\hrule height.4\ex@}%
 }}}
\message{commutative diagrams,}
\let\ampersand@\relax
\newdimen\minaw@
\minaw@11.11128\ex@
\newdimen\minCDaw@
\minCDaw@2.5pc
\def\minCDarrowwidth#1{\RIfMIfI@\onlydmatherr@\minCDarrowwidth
 \else\minCDaw@#1\relax\fi\else\onlydmatherr@\minCDarrowwidth\fi}
\newif\ifCD@
\def\CD{\bgroup\vspace@\relax\let\ampersand@&\iffalse}\fi
 \CD@true\vcenter\bgroup\Let@\tabskip\z@skip\baselineskip20\ex@
 \lineskip3\ex@\lineskiplimit3\ex@\halign\bgroup
 &\hfill$\m@th##$\hfill\crcr}
\def\endCD{\crcr\egroup\egroup\egroup}
\newdimen\bigaw@
\atdef@>#1>#2>{\ampersand@                                                  
 \setboxz@h{$\m@th\ssize\;{#1}\;\;$}
 \setbox@ne\hbox{$\m@th\ssize\;{#2}\;\;$}
 \setbox\tw@\hbox{$\m@th#2$}
 \ifCD@\global\bigaw@\minCDaw@\else\global\bigaw@\minaw@\fi                 
 \ifdim\wdz@>\bigaw@\global\bigaw@\wdz@\fi
 \ifdim\wd@ne>\bigaw@\global\bigaw@\wd@ne\fi                                
 \ifCD@\enskip\fi                                                           
 \ifdim\wd\tw@>\z@
  \mathrel{\mathop{\hbox to\bigaw@{\rightarrowfill@\displaystyle}}%
    \limits^{#1}_{#2}}
 \else\mathrel{\mathop{\hbox to\bigaw@{\rightarrowfill@\displaystyle}}%
    \limits^{#1}}\fi                                                        
 \ifCD@\enskip\fi                                                          
 \ampersand@}                                                              
\atdef@<#1<#2<{\ampersand@\setboxz@h{$\m@th\ssize\;\;{#1}\;$}%
 \setbox@ne\hbox{$\m@th\ssize\;\;{#2}\;$}\setbox\tw@\hbox{$\m@th#2$}%
 \ifCD@\global\bigaw@\minCDaw@\else\global\bigaw@\minaw@\fi
 \ifdim\wdz@>\bigaw@\global\bigaw@\wdz@\fi
 \ifdim\wd@ne>\bigaw@\global\bigaw@\wd@ne\fi
 \ifCD@\enskip\fi
 \ifdim\wd\tw@>\z@
  \mathrel{\mathop{\hbox to\bigaw@{\leftarrowfill@\displaystyle}}%
       \limits^{#1}_{#2}}\else
  \mathrel{\mathop{\hbox to\bigaw@{\leftarrowfill@\displaystyle}}%
       \limits^{#1}}\fi
 \ifCD@\enskip\fi\ampersand@}
\begingroup
 \catcode`\~=\active \lccode`\~=`\@
 \lowercase{%
  \global\atdef@)#1)#2){~>#1>#2>}
  \global\atdef@(#1(#2({~<#1<#2<}}
\endgroup
\atdef@ A#1A#2A{\llap{$\m@th\vcenter{\hbox
 {$\ssize#1$}}$}\Big\uparrow\rlap{$\m@th\vcenter{\hbox{$\ssize#2$}}$}&&}
\atdef@ V#1V#2V{\llap{$\m@th\vcenter{\hbox
 {$\ssize#1$}}$}\Big\downarrow\rlap{$\m@th\vcenter{\hbox{$\ssize#2$}}$}&&}
\atdef@={&\enskip\mathrel
 {\vbox{\hrule width\minCDaw@\vskip3\ex@\hrule width
 \minCDaw@}}\enskip&}
\atdef@|{\Big\Vert&&}
\atdef@\vert{\Big\Vert&&}
\def\pretend#1\haswidth#2{\setboxz@h{$\m@th\scriptstyle{#2}$}\hbox
 to\wdz@{\hfill$\m@th\scriptstyle{#1}$\hfill}}
\message{poor man's bold,}
\def\pmb{\RIfM@\expandafter\mathpalette\expandafter\pmb@\else
 \expandafter\pmb@@\fi}
\def\pmb@@#1{\leavevmode\setboxz@h{#1}%
   \dimen@-\wdz@
   \kern-.5\ex@\copy\z@
   \kern\dimen@\kern.25\ex@\raise.4\ex@\copy\z@
   \kern\dimen@\kern.25\ex@\box\z@
}
\def\binrel@@#1{\ifdim\wd2<\z@\mathbin{#1}\else\ifdim\wd\tw@>\z@
 \mathrel{#1}\else{#1}\fi\fi}
\newdimen\pmbraise@
\def\pmb@#1#2{\setbox\thr@@\hbox{$\m@th#1{#2}$}%
 \setbox4\hbox{$\m@th#1\mkern.5mu$}\pmbraise@\wd4\relax
 \binrel@{#2}%
 \dimen@-\wd\thr@@
   \binrel@@{%
   \mkern-.8mu\copy\thr@@
   \kern\dimen@\mkern.4mu\raise\pmbraise@\copy\thr@@
   \kern\dimen@\mkern.4mu\box\thr@@
}}
\def\documentstyle#1{\W@{}\input #1.sty\relax}
\message{syntax check,}
\font\dummyft@=dummy
\fontdimen1 \dummyft@=\z@
\fontdimen2 \dummyft@=\z@
\fontdimen3 \dummyft@=\z@
\fontdimen4 \dummyft@=\z@
\fontdimen5 \dummyft@=\z@
\fontdimen6 \dummyft@=\z@
\fontdimen7 \dummyft@=\z@
\fontdimen8 \dummyft@=\z@
\fontdimen9 \dummyft@=\z@
\fontdimen10 \dummyft@=\z@
\fontdimen11 \dummyft@=\z@
\fontdimen12 \dummyft@=\z@
\fontdimen13 \dummyft@=\z@
\fontdimen14 \dummyft@=\z@
\fontdimen15 \dummyft@=\z@
\fontdimen16 \dummyft@=\z@
\fontdimen17 \dummyft@=\z@
\fontdimen18 \dummyft@=\z@
\fontdimen19 \dummyft@=\z@
\fontdimen20 \dummyft@=\z@
\fontdimen21 \dummyft@=\z@
\fontdimen22 \dummyft@=\z@
\def\fontlist@{\\{\tenrm}\\{\sevenrm}\\{\fiverm}\\{\teni}\\{\seveni}%
 \\{\fivei}\\{\tensy}\\{\sevensy}\\{\fivesy}\\{\tenex}\\{\tenbf}\\{\sevenbf}%
 \\{\fivebf}\\{\tensl}\\{\tenit}}
\def\font@#1=#2 {\rightappend@#1\to\fontlist@\font#1=#2 }
\def\dodummy@{{\def\\##1{\global\let##1\dummyft@}\fontlist@}}
\def\nopages@{\output{\setbox\z@\box\@cclv \deadcycles\z@}%
 \alloc@5\toks\toksdef\@cclvi\output}
\let\galleys\nopages@
\newif\ifsyntax@
\newcount\countxviii@
\def\syntax{\syntax@true\dodummy@\countxviii@\count18
 \loop\ifnum\countxviii@>\m@ne\textfont\countxviii@=\dummyft@
 \scriptfont\countxviii@=\dummyft@\scriptscriptfont\countxviii@=\dummyft@
 \advance\countxviii@\m@ne\repeat                                           
 \dummyft@\tracinglostchars\z@\nopages@\frenchspacing\hbadness\@M}
\def\first@#1#2\end{#1}
\def\printoptions{\W@{Do you want S(yntax check),
  G(alleys) or P(ages)?}%
 \message{Type S, G or P, followed by <return>: }%
 \begingroup 
 \endlinechar\m@ne 
 \read\m@ne to\ans@
 \edef\ans@{\uppercase{\def\noexpand\ans@{%
   \expandafter\first@\ans@ P\end}}}%
 \expandafter\endgroup\ans@
 \if\ans@ P
 \else \if\ans@ S\syntax
 \else \if\ans@ G\galleys
 \else\message{? Unknown option: \ans@; using the `pages' option.}%
 \fi\fi\fi}
\def\alloc@#1#2#3#4#5{\global\advance\count1#1by\@ne
 \ch@ck#1#4#2\allocationnumber=\count1#1
 \global#3#5=\allocationnumber
 \ifalloc@\wlog{\string#5=\string#2\the\allocationnumber}\fi}
\def\document{\def\alloclist@{}\def\fontlist@{}}
\let\enddocument\bye

\let\proclaim\undefined
\let\footnote\undefined
\let\=\undefined
\let\>\undefined

\catcode`\@=\active
\message{... finished}

\expandafter\ifx\csname mathdefs.tex\endcsname\relax
  \expandafter\gdef\csname mathdefs.tex\endcsname{}
\else \message{Hey!  Apparently you were trying to
  \string\input{mathdefs.tex} twice.   This does not make sense.} 
\errmessage{Please edit your file (probably \jobname.tex) and remove
any duplicate ``\string\input'' lines}\endinput\fi




\catcode`\X=12\catcode`\@=11

\def\n@wcount{\alloc@0\count\countdef\insc@unt}
\def\n@wwrite{\alloc@7\write\chardef\sixt@@n}
\def\n@wread{\alloc@6\read\chardef\sixt@@n}
\def\r@s@t{\relax}\def\v@idline{\par}\def\@mputate#1/{#1}
\def\l@c@l#1X{\firstpart.#1}\def\gl@b@l#1X{#1}\def\t@d@l#1X{{}}

\def\crossrefs#1{\ifx\all#1\let\tr@ce=\all\else\def\tr@ce{#1,}\fi
   \n@wwrite\cit@tionsout\openout\cit@tionsout=\jobname.cit 
   \write\cit@tionsout{\tr@ce}\expandafter\setfl@gs\tr@ce,}
\def\setfl@gs#1,{\def\@{#1}\ifx\@\empty\let\next=\relax
   \else\let\next=\setfl@gs\expandafter\xdef
   \csname#1tr@cetrue\endcsname{}\fi\next}
\def\m@ketag#1#2{\expandafter\n@wcount\csname#2tagno\endcsname
     \csname#2tagno\endcsname=0\let\tail=\all\xdef\all{\tail#2,}
   \ifx#1\l@c@l\let\tail=\r@s@t\xdef\r@s@t{\csname#2tagno\endcsname=0\tail}\fi
   \expandafter\gdef\csname#2cite\endcsname##1{\expandafter
     \ifx\csname#2tag##1\endcsname\relax?\else\csname#2tag##1\endcsname\fi
     \expandafter\ifx\csname#2tr@cetrue\endcsname\relax\else
     \write\cit@tionsout{#2tag ##1 cited on page \folio.}\fi}
   \expandafter\gdef\csname#2page\endcsname##1{\expandafter
     \ifx\csname#2page##1\endcsname\relax?\else\csname#2page##1\endcsname\fi
     \expandafter\ifx\csname#2tr@cetrue\endcsname\relax\else
     \write\cit@tionsout{#2tag ##1 cited on page \folio.}\fi}
   \expandafter\gdef\csname#2tag\endcsname##1{\expandafter
      \ifx\csname#2check##1\endcsname\relax
      \expandafter\xdef\csname#2check##1\endcsname{}%
      \else\immediate\write16{Warning: #2tag ##1 used more than once.}\fi
      \multit@g{#1}{#2}##1/X%
      \write\t@gsout{#2tag ##1 assigned number \csname#2tag##1\endcsname\space
      on page \number\count0.}%
   \csname#2tag##1\endcsname}}

\def\multit@g#1#2#3/#4X{\def\t@mp{#4}\ifx\t@mp\empty%
      \global\advance\csname#2tagno\endcsname by 1 
      \expandafter\xdef\csname#2tag#3\endcsname
      {#1\number\csname#2tagno\endcsnameX}%
   \else\expandafter\ifx\csname#2last#3\endcsname\relax
      \expandafter\n@wcount\csname#2last#3\endcsname
      \global\advance\csname#2tagno\endcsname by 1 
      \expandafter\xdef\csname#2tag#3\endcsname
      {#1\number\csname#2tagno\endcsnameX}
      \write\t@gsout{#2tag #3 assigned number \csname#2tag#3\endcsname\space
      on page \number\count0.}\fi
   \global\advance\csname#2last#3\endcsname by 1
   \def\t@mp{\expandafter\xdef\csname#2tag#3/}%
   \expandafter\t@mp\@mputate#4\endcsname
   {\csname#2tag#3\endcsname\lastpart{\csname#2last#3\endcsname}}\fi}
\def\t@gs#1{\def\all{}\m@ketag#1e\m@ketag#1s\m@ketag\t@d@l p
\let\realscite\scite
\let\realstag\stag
   \m@ketag\gl@b@l r \n@wread\t@gsin
   \openin\t@gsin=\jobname.tgs \re@der \closein\t@gsin
   \n@wwrite\t@gsout\openout\t@gsout=\jobname.tgs }
\outer\def\localtags{\t@gs\l@c@l}
\outer\def\globaltags{\t@gs\gl@b@l}
\outer\def\newlocaltag#1{\m@ketag\l@c@l{#1}}
\outer\def\newglobaltag#1{\m@ketag\gl@b@l{#1}}

\newif\ifpr@ 
\def\m@kecs #1tag #2 assigned number #3 on page #4.%
   {\expandafter\gdef\csname#1tag#2\endcsname{#3}
   \expandafter\gdef\csname#1page#2\endcsname{#4}
   \ifpr@\expandafter\xdef\csname#1check#2\endcsname{}\fi}
\def\re@der{\ifeof\t@gsin\let\next=\relax\else
   \read\t@gsin to\t@gline\ifx\t@gline\v@idline\else
   \expandafter\m@kecs \t@gline\fi\let \next=\re@der\fi\next}
\def\pretags#1{\pr@true\pret@gs#1,,}
\def\pret@gs#1,{\def\@{#1}\ifx\@\empty\let\n@xtfile=\relax
   \else\let\n@xtfile=\pret@gs \openin\t@gsin=#1.tgs \message{#1} \re@der 
   \closein\t@gsin\fi \n@xtfile}

\newcount\sectno\sectno=0\newcount\subsectno\subsectno=0
\newif\ifultr@local \def\ultralocal{\ultr@localtrue}
\def\firstpart{\number\sectno}
\def\lastpart#1{\ifcase#1 \or a\or b\or c\or d\or e\or f\or g\or h\or 
   i\or k\or l\or m\or n\or o\or p\or q\or r\or s\or t\or u\or v\or w\or 
   x\or y\or z \fi}

\def\resetall{\global\advance\sectno by 1\subsectno=0
   \gdef\firstpart{\number\sectno}\r@s@t}
\def\resetsub{\global\advance\subsectno by 1
   \gdef\firstpart{\number\sectno.\number\subsectno}\r@s@t}
\def\newsection#1\par{\resetall\vskip0pt plus.3\vsize\penalty-250
   \vskip0pt plus-.3\vsize\bigskip\bigskip
   \message{#1}\leftline{\bf#1}\nobreak\bigskip}
\def\subsection#1\par{\ifultr@local\resetsub\fi
   \vskip0pt plus.2\vsize\penalty-250\vskip0pt plus-.2\vsize
   \bigskip\smallskip\message{#1}\leftline{\bf#1}\nobreak\medskip}


\newdimen\marginshift

\newdimen\margindelta
\newdimen\marginmax
\newdimen\marginmin

\def\margininit{       
\marginmax=3 true cm                  
				      
\margindelta=0.1 true cm              
\marginmin=0.1true cm                 
\marginshift=\marginmin
}    

\def\t@gsjj#1,{\def\@{#1}\ifx\@\empty\let\next=\relax\else\let\next=\t@gsjj
   \def\@@{p}\ifx\@\@@\else
   \expandafter\gdef\csname#1cite\endcsname##1{\citejj{##1}}
   \expandafter\gdef\csname#1page\endcsname##1{?}
   \expandafter\gdef\csname#1tag\endcsname##1{\tagjj{##1}}\fi\fi\next}
\newif\ifshowstuffinmargin
\showstuffinmarginfalse
\def\jjtags{\ifx\shlhetal\relax 
  \else
\ifx\shlhetal\undefinedcontrolseq
\else
\showstuffinmargintrue
\ifx\all\relax\else\expandafter\t@gsjj\all,\fi\fi \fi
}

\def\tagjj#1{\realstag{#1}\mginpar{\zeigen{#1}}}
\def\citejj#1{\zeigen{#1}\mginpar{\rechnen{#1}}}

\def\rechnen#1{\expandafter\ifx\csname stag#1\endcsname\relax ??\else
                           \csname stag#1\endcsname\fi}

\newdimen\theight

\def\marginfont{\sevenrm}

\def\trymarginbox#1{\setbox0=\hbox{\marginfont\hskip\marginshift #1}%
		\global\marginshift\wd0 
		\global\advance\marginshift\margindelta}

\def \mginpar#1{%
\ifvmode\setbox0\hbox to \hsize{\hfill\rlap{\marginfont\quad#1}}%
\ht0 0cm
\dp0 0cm
\box0\vskip-\baselineskip
\else 
             \vadjust{\trymarginbox{#1}%
		\ifdim\marginshift>\marginmax \global\marginshift\marginmin
			\trymarginbox{#1}%
                \fi
             \theight=\ht0
             \advance\theight by \dp0    \advance\theight by \lineskip
             \kern -\theight \vbox to \theight{\rightline{\rlap{\box0}}%
\vss}}\fi}


\def\t@gsoff#1,{\def\@{#1}\ifx\@\empty\let\next=\relax\else\let\next=\t@gsoff
   \def\@@{p}\ifx\@\@@\else
   \expandafter\gdef\csname#1cite\endcsname##1{\zeigen{##1}}
   \expandafter\gdef\csname#1page\endcsname##1{?}
   \expandafter\gdef\csname#1tag\endcsname##1{\zeigen{##1}}\fi\fi\next}
\def\verbatimtags{\showstuffinmarginfalse
\ifx\all\relax\else\expandafter\t@gsoff\all,\fi}
\def\zeigen#1{\hbox{$\langle$}#1\hbox{$\rangle$}}
\def\margincite#1{\ifshowstuffinmargin\mginpar{\rechnen{#1}}\fi}

\def\(#1){\edef\dot@g{\ifmmode\ifinner(\hbox{\noexpand\etag{#1}})
   \else\noexpand\eqno(\hbox{\noexpand\etag{#1}})\fi
   \else(\noexpand\ecite{#1})\fi}\dot@g}

\newif\ifbr@ck
\def\eat#1{}
\def\[#1]{\br@cktrue[\br@cket#1'X]}
\def\br@cket#1'#2X{\def\temp{#2}\ifx\temp\empty\let\next\eat
   \else\let\next\br@cket\fi
   \ifbr@ck\br@ckfalse\br@ck@t#1,X\else\br@cktrue#1\fi\next#2X}
\def\br@ck@t#1,#2X{\def\temp{#2}\ifx\temp\empty\let\neext\eat
   \else\let\neext\br@ck@t\def\temp{,}\fi
   \def\teemp{#1}\ifx\teemp\empty\else\rcite{#1}\fi\temp\neext#2X}
\def\resetbr@cket{\gdef\[##1]{[\rtag{##1}]}}
\def\references{\resetbr@cket\newsection References\par}

\newtoks\symb@ls\newtoks\s@mb@ls\newtoks\p@gelist\n@wcount\ftn@mber
    \ftn@mber=1\newif\ifftn@mbers\ftn@mbersfalse\newif\ifbyp@ge\byp@gefalse
\def\defm@rk{\ifftn@mbers\n@mberm@rk\else\symb@lm@rk\fi}
\def\n@mberm@rk{\xdef\m@rk{{\the\ftn@mber}}%
    \global\advance\ftn@mber by 1 }
\def\rot@te#1{\let\temp=#1\global#1=\expandafter\r@t@te\the\temp,X}
\def\r@t@te#1,#2X{{#2#1}\xdef\m@rk{{#1}}}
\def\b@@st#1{{$^{#1}$}}\def\str@p#1{#1}
\def\symb@lm@rk{\ifbyp@ge\rot@te\p@gelist\ifnum\expandafter\str@p\m@rk=1 
    \s@mb@ls=\symb@ls\fi\write\f@nsout{\number\count0}\fi \rot@te\s@mb@ls}
\def\byp@ge{\byp@getrue\n@wwrite\f@nsin\openin\f@nsin=\jobname.fns 
    \n@wcount\currentp@ge\currentp@ge=0\p@gelist={0}
    \re@dfns\closein\f@nsin\rot@te\p@gelist
    \n@wread\f@nsout\openout\f@nsout=\jobname.fns }
\def\m@kelist#1X#2{{#1,#2}}
\def\re@dfns{\ifeof\f@nsin\let\next=\relax\else\read\f@nsin to \f@nline
    \ifx\f@nline\v@idline\else\let\t@mplist=\p@gelist
    \ifnum\currentp@ge=\f@nline
    \global\p@gelist=\expandafter\m@kelist\the\t@mplistX0
    \else\currentp@ge=\f@nline
    \global\p@gelist=\expandafter\m@kelist\the\t@mplistX1\fi\fi
    \let\next=\re@dfns\fi\next}
\def\symbols#1{\symb@ls={#1}\s@mb@ls=\symb@ls} 
\def\bigsymbol{\textstyle}
\symbols{\bigsymbol\ast,\dagger,\ddagger,\sharp,\flat,\natural,\star}
\def\ftnumbers{\ftn@mberstrue} \def\ftsymbols{\ftn@mbersfalse}
\def\paginal{\byp@ge} \def\resetftnumbers{\ftn@mber=1}
\def\ftnote#1{\defm@rk\expandafter\expandafter\expandafter\footnote
    \expandafter\b@@st\m@rk{#1}}

\long\def\jump#1\endjump{}
\def\ssum{\mathop{\lower .1em\hbox{$\textstyle\Sigma$}}\nolimits}

\def\qed{\nobreak\kern 1em \vrule height .5em width .5em depth 0em}
\def\newneq{\hbox{\rlap{\hbox to 1\wd9{\hss$=$\hss}}\raise .1em 
   \hbox to 1\wd9{\hss$\scriptscriptstyle/$\hss}}}
\def\subsetne{\setbox9 = \hbox{$\subset$}\mathrel{\hbox{\rlap
   {\lower .4em \newneq}\raise .13em \hbox{$\subset$}}}}
\def\supsetne{\setbox9 = \hbox{$\subset$}\mathrel{\hbox{\rlap
   {\lower .4em \newneq}\raise .13em \hbox{$\supset$}}}}

\def\vbar{\mathchoice{\vrule height6.3ptdepth-.5ptwidth.8pt\kern-.8pt}
   {\vrule height6.3ptdepth-.5ptwidth.8pt\kern-.8pt}
   {\vrule height4.1ptdepth-.35ptwidth.6pt\kern-.6pt}
   {\vrule height3.1ptdepth-.25ptwidth.5pt\kern-.5pt}}
\def\f@dge{\mathchoice{}{}{\mkern.5mu}{\mkern.8mu}}
\def\b@c#1#2{{\rm \mkern#2mu\vbar\mkern-#2mu#1}}
\def\b@b#1{{\rm I\mkern-3.5mu #1}}
\def\b@a#1#2{{\rm #1\mkern-#2mu\f@dge #1}}
\def\bb#1{{\count4=`#1 \advance\count4by-64 \ifcase\count4\or\b@a A{11.5}\or
   \b@b B\or\b@c C{5}\or\b@b D\or\b@b E\or\b@b F \or\b@c G{5}\or\b@b H\or
   \b@b I\or\b@c J{3}\or\b@b K\or\b@b L \or\b@b M\or\b@b N\or\b@c O{5} \or
   \b@b P\or\b@c Q{5}\or\b@b R\or\b@a S{8}\or\b@a T{10.5}\or\b@c U{5}\or
   \b@a V{12}\or\b@a W{16.5}\or\b@a X{11}\or\b@a Y{11.7}\or\b@a Z{7.5}\fi}}

\catcode`\X=11 \catcode`\@=12


\expandafter\ifx\csname citeadd.tex\endcsname\relax
\expandafter\gdef\csname citeadd.tex\endcsname{}
\else \message{Hey!  Apparently you were trying to
\string\input{citeadd.tex} twice.   This does not make sense.} 
\errmessage{Please edit your file (probably \jobname.tex) and remove
any duplicate ``\string\input'' lines}\endinput\fi

\def\sciteu{\sciteerror{undefined}}

\def\sciteerror#1#2{{\mathortextbf{\scite{#2}}}\complainaboutcitation{#1}{#2}}
\def\mathortextbf#1{\hbox{\bf #1}}
\def\complainaboutcitation#1#2{%
\vadjust{\line{\llap{---$\!\!>$ }\qquad scite$\{$#2$\}$ #1\hfil}}}

\sectno=-1   
\localtags
\NoBlackBoxes
\define\mr{\medskip\roster}
\define\sn{\smallskip\noindent}
\define\mn{\medskip\noindent}
\define\bn{\bigskip\noindent}
\define\ub{\underbar}
\define\wilog{\text{without loss of generality}}
\define\ermn{\endroster\medskip\noindent}
\define\dbca{\dsize\bigcap}
\define\dbcu{\dsize\bigcup}
\define \nl{\newline}
\magnification=\magstep 1
\documentstyle{amsppt}

{    
\catcode`@11

\ifx\alicetwothousandloaded@\relax
  \endinput\else\global\let\alicetwothousandloaded@\relax\fi

\gdef\subjclass{\let\savedef@\subjclass
 \def\subjclass##1\endsubjclass{\let\subjclass\savedef@
   \toks@{\def\usualspace{{\rm\enspace}}\eightpoint}%
   \toks@@{##1\unskip.}%
   \edef\thesubjclass@{\the\toks@
     \frills@{{\noexpand\rm2000 {\noexpand\it Mathematics Subject
       Classification}.\noexpand\enspace}}%
     \the\toks@@}}%
  \nofrillscheck\subjclass}
} 


\expandafter\ifx\csname alice2jlem.tex\endcsname\relax
  \expandafter\xdef\csname alice2jlem.tex\endcsname{\the\catcode`@}
\else \message{Hey!  Apparently you were trying to
\string\input{alice2jlem.tex}  twice.   This does not make sense.}
\errmessage{Please edit your file (probably \jobname.tex) and remove
any duplicate ``\string\input'' lines}\endinput\fi

\expandafter\ifx\csname bib4plain.tex\endcsname\relax
  \expandafter\gdef\csname bib4plain.tex\endcsname{}
\else \message{Hey!  Apparently you were trying to \string\input
  bib4plain.tex twice.   This does not make sense.}
\errmessage{Please edit your file (probably \jobname.tex) and remove
any duplicate ``\string\input'' lines}\endinput\fi

\def\renewcommand{\newcommand}	       
\edef\cite{\the\catcode`@}%
\catcode`@ = 11
\let\@oldatcatcode = \cite
\chardef\@letter = 11
\chardef\@other = 12
%
%
%
%
\def\@innerdef#1#2{\edef#1{\expandafter\noexpand\csname #2\endcsname}}%
%
%
\@innerdef\@innernewcount{newcount}%
\@innerdef\@innernewdimen{newdimen}%
\@innerdef\@innernewif{newif}%
\@innerdef\@innernewwrite{newwrite}%
%
%
%
\def\@gobble#1{}%
%
%
%
\ifx\inputlineno\@undefined
   \let\@linenumber = \empty 
\else
   \def\@linenumber{\the\inputlineno:\space}%
\fi
%
%
%
\def\@futurenonspacelet#1{\def\cs{#1}%
   \afterassignment\@stepone\let\@nexttoken=
}%
\begingroup 
\def\\{\global\let\@stoken= }%
\\ 
\endgroup
\def\@stepone{\expandafter\futurelet\cs\@steptwo}%
\def\@steptwo{\expandafter\ifx\cs\@stoken\let\@@next=\@stepthree
   \else\let\@@next=\@nexttoken\fi \@@next}%
\def\@stepthree{\afterassignment\@stepone\let\@@next= }%
%
%
%
\def\@getoptionalarg#1{%
   \let\@optionaltemp = #1%
   \let\@optionalnext = \relax
   \@futurenonspacelet\@optionalnext\@bracketcheck
}%
%
%
\def\@bracketcheck{%
   \ifx [\@optionalnext
      \expandafter\@@getoptionalarg
   \else
      \let\@optionalarg = \empty
      \expandafter\@optionaltemp
   \fi
}%
\def\@@getoptionalarg[#1]{%
   \def\@optionalarg{#1}%
   \@optionaltemp
}%
%
%
%
\def\@nnil{\@nil}%
\def\@fornoop#1\@@#2#3{}%
\def\@for#1:=#2\do#3{%
   \edef\@fortmp{#2}%
   \ifx\@fortmp\empty \else
      \expandafter\@forloop#2,\@nil,\@nil\@@#1{#3}%
   \fi
}%
\def\@forloop#1,#2,#3\@@#4#5{\def#4{#1}\ifx #4\@nnil \else
       #5\def#4{#2}\ifx #4\@nnil \else#5\@iforloop #3\@@#4{#5}\fi\fi
}%
\def\@iforloop#1,#2\@@#3#4{\def#3{#1}\ifx #3\@nnil
       \let\@nextwhile=\@fornoop \else
      #4\relax\let\@nextwhile=\@iforloop\fi\@nextwhile#2\@@#3{#4}%
}%
%
%
%
\@innernewif\if@fileexists
\def\@testfileexistence{\@getoptionalarg\@finishtestfileexistence}%
\def\@finishtestfileexistence#1{%
   \begingroup
      \def\extension{#1}%
      \immediate\openin0 =
         \ifx\@optionalarg\empty\jobname\else\@optionalarg\fi
         \ifx\extension\empty \else .#1\fi
         \space
      \ifeof 0
         \global\@fileexistsfalse
      \else
         \global\@fileexiststrue
      \fi
      \immediate\closein0
   \endgroup
}%
%
%
%
%
\def\bibliographystyle#1{%
   \@readauxfile
   \@writeaux{\string\bibstyle{#1}}%
}%
\let\bibstyle = \@gobble
%
%
\let\bblfilebasename = \jobname
\def\bibliography#1{%
   \@readauxfile
   \@writeaux{\string\bibdata{#1}}%
   \@testfileexistence[\bblfilebasename]{bbl}%
   \if@fileexists
      \nobreak
      \@readbblfile
   \fi
}%
\let\bibdata = \@gobble
%
%
\def\nocite#1{%
   \@readauxfile
   \@writeaux{\string\citation{#1}}%
}%
\@innernewif\if@notfirstcitation
%
%
\def\cite{\@getoptionalarg\@cite}%
%
%
\def\@cite#1{%
   \let\@citenotetext = \@optionalarg
   \printcitestart
   \nocite{#1}%
   \@notfirstcitationfalse
   \@for \@citation :=#1\do
   {%
      \expandafter\@onecitation\@citation\@@
   }%
   \ifx\empty\@citenotetext\else
      \printcitenote{\@citenotetext}%
   \fi
   \printcitefinish
}%
\def\@onecitation#1\@@{%
   \if@notfirstcitation
      \printbetweencitations
   \fi
   \expandafter \ifx \csname\@citelabel{#1}\endcsname \relax
      \if@citewarning
         \message{\@linenumber Undefined citation `#1'.}%
      \fi
      \expandafter\gdef\csname\@citelabel{#1}\endcsname{%
\strut
\vadjust{\vskip-\dp\strutbox
\vbox to 0pt{\vss\parindent0cm \leftskip=\hsize 
\advance\leftskip3mm
\advance\hsize 4cm\strut\openup-4pt 
\rightskip 0cm plus 1cm minus 0.5cm ?  #1 ?\strut}}
         {\tt
            \escapechar = -1
            \nobreak\hskip0pt
            \expandafter\string\csname#1\endcsname
            \nobreak\hskip0pt
         }%
      }%
   \fi
   \csname\@citelabel{#1}\endcsname
   \@notfirstcitationtrue
}%
%
%
\def\@citelabel#1{b@#1}%
%
%
\def\@citedef#1#2{\expandafter\gdef\csname\@citelabel{#1}\endcsname{#2}}%
%
%
%
\def\@readbblfile{%
   \ifx\@itemnum\@undefined
      \@innernewcount\@itemnum
   \fi
   \begingroup
      \def\begin##1##2{%
         \setbox0 = \hbox{\biblabelcontents{##2}}%
         \biblabelwidth = \wd0
      }%
      \def\end##1{}
      %
      %
      \@itemnum = 0
      \def\bibitem{\@getoptionalarg\@bibitem}%
      \def\@bibitem{%
         \ifx\@optionalarg\empty
            \expandafter\@numberedbibitem
         \else
            \expandafter\@alphabibitem
         \fi
      }%
      \def\@alphabibitem##1{%
         \expandafter \xdef\csname\@citelabel{##1}\endcsname {\@optionalarg}%
         \ifx\biblabelprecontents\@undefined
            \let\biblabelprecontents = \relax
         \fi
         \ifx\biblabelpostcontents\@undefined
            \let\biblabelpostcontents = \hss
         \fi
         \@finishbibitem{##1}%
      }%
      \def\@numberedbibitem##1{%
         \advance\@itemnum by 1
         \expandafter \xdef\csname\@citelabel{##1}\endcsname{\number\@itemnum}%
         \ifx\biblabelprecontents\@undefined
            \let\biblabelprecontents = \hss
         \fi
         \ifx\biblabelpostcontents\@undefined
            \let\biblabelpostcontents = \relax
         \fi
         \@finishbibitem{##1}%
      }%
      \def\@finishbibitem##1{%
         \biblabelprint{\csname\@citelabel{##1}\endcsname}%
         \@writeaux{\string\@citedef{##1}{\csname\@citelabel{##1}\endcsname}}%
         \ignorespaces
      }%
      %
      %
      \let\em = \bblem
      \let\newblock = \bblnewblock
      \let\sc = \bblsc
      \frenchspacing
      \clubpenalty = 4000 \widowpenalty = 4000
      \tolerance = 10000 \hfuzz = .5pt
      \everypar = {\hangindent = \biblabelwidth
                      \advance\hangindent by \biblabelextraspace}%
      \bblrm
      \parskip = 1.5ex plus .5ex minus .5ex
      \biblabelextraspace = .5em
      \bblhook
      \input \bblfilebasename.bbl
   \endgroup
}%
%
%
\@innernewdimen\biblabelwidth
\@innernewdimen\biblabelextraspace
%
%
%
\def\biblabelprint#1{%
   \noindent
   \hbox to \biblabelwidth{%
      \biblabelprecontents
      \biblabelcontents{#1}%
      \biblabelpostcontents
   }%
   \kern\biblabelextraspace
}%
%
%
%
\def\biblabelcontents#1{{\bblrm [#1]}}%
%
%
\def\bblrm{\rm}%
%
%
\def\bblem{\it}%
%
%
\def\bblsc{\ifx\@scfont\@undefined
              \font\@scfont = cmcsc10
           \fi
           \@scfont
}%
%
%
\def\bblnewblock{\hskip .11em plus .33em minus .07em }%
%
%
\let\bblhook = \empty
%
%
%
\def\printcitestart{[}
\def\printcitefinish{]}
\def\printbetweencitations{, }
\def\printcitenote#1{, #1}
%
%
%
\let\citation = \@gobble
%
%
%
\@innernewcount\@numparams
%
%
\def\newcommand#1{%
   \def\@commandname{#1}%
   \@getoptionalarg\@continuenewcommand
}%
%
%
\def\@continuenewcommand{%
   \@numparams = \ifx\@optionalarg\empty 0\else\@optionalarg \fi \relax
   \@newcommand
}%
%
%
\def\@newcommand#1{%
   \def\@startdef{\expandafter\edef\@commandname}%
   \ifnum\@numparams=0
      \let\@paramdef = \empty
   \else
      \ifnum\@numparams>9
         \errmessage{\the\@numparams\space is too many parameters}%
      \else
         \ifnum\@numparams<0
            \errmessage{\the\@numparams\space is too few parameters}%
         \else
            \edef\@paramdef{%
               \ifcase\@numparams
                  \empty  No arguments.
               \or ####1%
               \or ####1####2%
               \or ####1####2####3%
               \or ####1####2####3####4%
               \or ####1####2####3####4####5%
               \or ####1####2####3####4####5####6%
               \or ####1####2####3####4####5####6####7%
               \or ####1####2####3####4####5####6####7####8%
               \or ####1####2####3####4####5####6####7####8####9%
               \fi
            }%
         \fi
      \fi
   \fi
   \expandafter\@startdef\@paramdef{#1}%
}%
%
%
%
%
\def\@readauxfile{%
   \if@auxfiledone \else 
      \global\@auxfiledonetrue
      \@testfileexistence{aux}%
      \if@fileexists
         \begingroup
            \endlinechar = -1
            \catcode`@ = 11
            \input \jobname.aux
         \endgroup
      \else
         \message{\@undefinedmessage}%
         \global\@citewarningfalse
      \fi
      \immediate\openout\@auxfile = \jobname.aux
   \fi
}%
%
%
\newif\if@auxfiledone
\ifx\noauxfile\@undefined \else \@auxfiledonetrue\fi
%
%
%
%
\@innernewwrite\@auxfile
\def\@writeaux#1{\ifx\noauxfile\@undefined \write\@auxfile{#1}\fi}%
%
%
%
\ifx\@undefinedmessage\@undefined
   \def\@undefinedmessage{No .aux file; I won't give you warnings about
                          undefined citations.}%
\fi
%
%
\@innernewif\if@citewarning
\ifx\noauxfile\@undefined \@citewarningtrue\fi
%
%
%
\catcode`@ = \@oldatcatcode


\def\widestnumber#1#2{}

\def\rm{\fam0 \tenrm}

\def\fakesubhead#1\endsubhead{\bigskip\noindent{\bf#1}\par}



%
%
%

%

\font\textrsfs=rsfs10
\font\scriptrsfs=rsfs7
\font\scriptscriptrsfs=rsfs5

\newfam\rsfsfam
\textfont\rsfsfam=\textrsfs
\scriptfont\rsfsfam=\scriptrsfs
\scriptscriptfont\rsfsfam=\scriptscriptrsfs

\edef\oldcatcodeofat{\the\catcode`\@}
\catcode`\@11

\def\Cal@@#1{\noaccents@ \fam \rsfsfam #1}

\catcode`\@\oldcatcodeofat


\expandafter\ifx \csname margininit\endcsname \relax\else\margininit\fi

\pageheight{8.5truein}
\topmatter
\title{Long iterations for the continuum  \\
Sh707} \endtitle
\author {Saharon Shelah \thanks {\null\newline I would like to thank 
Alice Leonhardt for the beautiful typing. \null\newline
Partially supported by the United States-Israel Binational Science
Foundation \null\newline
Latest Revision - 01/Oct/9 \null\newline
First version written 99/10;\S9;done 98/10 and \S10 - 12/99} \endthanks} \endauthor 
  
\affil{Institute of Mathematics\\
 The Hebrew University\\
 Jerusalem, Israel
 \medskip
 Rutgers University\\
 Mathematics Department\\
 New Brunswick, NJ  USA} \endaffil

\abstract  We deal with an iteration theorem of forcing notion with a
kind of countable support of nice enough forcing notion which is
proper $\aleph_2$-c.c. forcing notions.  We then look at some special
cases $(\Bbb Q_D$'s preceded by random forcing). \endabstract
\endtopmatter
\document

\head {Anotated content} \endhead  \resetall 
\bn
\S0 $\quad$ Introduction
\bn
\S1 $\quad$ Trunk Controller
\sn
[We define ``trunk controller" which will serve as the ``apure", not a
name of a kind of CS iteration.  We define standard; what is being
based and fully based on $\langle X_\beta:\beta < \alpha \rangle$ (at
in Definition \scite{it.1}).  Also we define simple, semi-simple,
semi-simply based on (Definition \scite{it.1a}).  We then define an
${\Cal F}$-forcing $\Bbb Q$ (\scite{it.2}), being clear, basic,
straight ${\Cal F}$-clear and weakly clear (Definition \scite{it.2a}).
We define ${\Cal F}$-iteration (\scite{it.4}).  We then prove some
basic claims.  Lastly, we define $(\theta,\sigma)$-pure decidability.]
\bn
\S2 $\quad$ Being ${\Cal F}$-pseudo c.c.c. is preserved by ${\Cal F}$-iterations
\sn
[We define ${\Cal F}$-psc, condition guaranteeing $\aleph_1$ is not
collapsed, an explicit form of properness and variants
(\scite{ct.1},\scite{ct.2},\scite{ct.3},\scite{ct.2b}) give sufficient
conditions (\scite{ct.3a},\scite{ct.2c}).  We give sufficient
conditions for pure decidability in claim \scite{ct.3}.  We prove if
$\bar Q$ is ${\Cal F}$-psc iteration then Lim$_{\Cal F}(\bar{\Bbb Q})$
is a ${\Cal F}$-psc forcing (+ variants, in Lemma \scite{ct.4}).  Give
a definition of witnesses for c.c.c. by sets of pairs.]
\bn
\S3 $\quad$ Nicer pure properness and pure decidability
\sn
[We return to condition for pure properness (claim \scite{mr.1}).]
\bn
\S4 $\quad$ Averages by an ultrafilter and restricted non null trees
\sn
[We consider the relationships of a (non principal) ultrafilter $D$ on
$\omega$ and subtrees $T$ of ${}^{\omega >} 2$ which positive Lebesgue
measure, considering $T = \text{ Lim}_D\langle T_n:n < \omega
\rangle$.  We concentrate on the case of the rate of convergence of
$T$ and $T_n$ to their measure is bounded from above a function $g$
from a family ${\Cal G} \subseteq {}^{\omega >} \omega$ of reasonable
candidates.]
\bn
\S5 $\quad$ On $\Bbb Q_{\bar D}$ and iterations
\bn
\S6 $\quad$ On a relative of Borel conjecture with large ${\frak b}$
\bn
\S7 $\quad$ Continuing \cite{Sh:592}
\bn
\S8 $\quad$ On ``$\eta$ is $R$-big over $M$"
\sn
[We generalize the ``$\eta$ is ${\Cal G}$-continuous over $M$", to
``$\eta$ is ${\Cal L}$-big over $M$".  So ``$\eta \in \text{ lim } T$"
is replaced by $(\eta,\nu) \in R \subseteq  {}^\omega 2 \times
{}^\omega 2$.]
\newpage

\head {\S0 Introduction} \endhead  \resetall \sectno=0
\bigskip

This is a modest try to investigate iterations $\bar{\Bbb Q} = \langle
\Bbb P_\alpha,{\underset\tilde {}\to {\Bbb Q}_\alpha}:\alpha <
\alpha^* \rangle$ which increase the continuum arbitrarily.  The
support is countable, but defining $p \le q$, only for finitely many
$\alpha \in \text{ Dom}(p)$ we are allowed to fail to have pure
extension.  More explicitly for every $p \in \Bbb Q_\alpha$, has a
``trunk" tr$(p)$, the apure part, and we demand that $\langle
tr(p(\alpha)):\alpha \in \text{ Dom}(p)\rangle$ is an ``old" element,
i.e. a function from $\bold V$.  In this context we have a quite
explicit form of properness which guaranteed $\aleph_1$ is not
collapsed.  Assuming CH there are reasonable conditions guaranteeing
the $\aleph_2$-c.c.

We may be more liberal in the first step of the iteration.  We then
concentrate on more specific content.  We let $\Bbb Q_0$ be
Random$_A$, adding a sequence of random reals $\langle
{\underset\tilde {}\to \nu_\gamma}:\gamma \in A \rangle$, and each
$\Bbb Q_\alpha = Q_{1 + \beta}$ is $\Bbb Q_{\underset\tilde {}\to
{\bar D}_\alpha},
{\underset\tilde {}\to {\bar D}_\alpha} = \langle {\underset\tilde
{}\to D^\alpha_\eta}:\eta \in {}^{\omega >} \omega
\rangle,{\underset\tilde {}\to D^\alpha_\eta}$ a $\Bbb P_\alpha$-anme
of a non principal ultrafilter on $\omega$.  However, for the results
we have in mind, ${\underset\tilde {}\to D^\alpha_\eta}$ should
satisfy some special properties: in the direction of being a Ramsey
ultrafilter.  If $\Bbb Q_0 =$ Random$_\lambda$, we may try to demand
that for every $r \in \bold V^{\lim(\bar{\Bbb Q})}$, for more $\beta <
\lambda,{\underset\tilde {}\to \nu_\beta}$ is random over $\bold
V[\bold r]$.  We do not know to do it, \ub{but} if we can restrict
ourselves to measure 1 set of the form $\cup \lim(T^{<n>}),T$ a
subtree of ${}^{\omega >}2$ with the fastness of convergence of
$\langle|2^n \cap T|/2^{-n}:n < \omega \rangle$ to Leb(lim$(T)$) by $g
\in \bold V$, moreover this holds above any $\eta \in {}^{\omega >}2$.
This is a ``poor relative" of the ``Borel conjecture + ${\frak b}$
large". \nl
The method seems to me more versatile than the method of first forcing
whatever and then forcing with the random algebra.  Note that:

We lastly deal with relative of \cite{Sh:592}.
\newpage

\head {\S1 Trunk Controllers} \endhead  \resetall \sectno=1
\bigskip

The reader may in \scite{it.1}(2) use $\zeta \le 1$, use only the
fully based case, and ignore \scite{it.9} (associativity).
\definition{\stag{it.1} Definition}:  1) A trunk controller ${\Cal F}$ is a set or
a class with quasi-orders $\le = \le^{\Cal F}$ (which we denote also
by $\le_{us}$) and $\le_{\text{pr}} = 
\le^{\Cal F}_{\text{pr}}$
and $\le_{\text{apr}} = \le^{\Cal F}_{\text{apr}}$ such that:
\mr
\item "{$(a)$}"  $\le_{\text{pr}} \subseteq \le$ and $\le_{\text{apr}}
\subseteq \le$.
\ermn
2) We define the $(\bar \alpha,X,\zeta)$-standard trunk controller 
${\Cal F} = {\Cal F}_{{\bar \alpha},\zeta}[X]$ in $\bold V$ by induction on 
$\zeta \le \omega_1$, where $\bar \alpha = \langle
\alpha_\varepsilon:\varepsilon \le \zeta \rangle,\alpha_\varepsilon$
an ordinal and $X$ is a trunk controller, but we may write $\bar
\alpha'$ with $\bar \alpha' \restriction (\zeta + 1) = \bar \alpha$ instead
of $\bar \alpha$: the $(\bar \alpha,X,\zeta)$-standard trunk
controller ${\Cal F}$ in $\bold V$ is:
\mr
\item "{$(a)$}"  the set of elements is the set of functions $f$ from a
countable subset of $\alpha_\zeta$ into $X \cup \bigcup \{{\Cal F}_{{\bar
\alpha},\varepsilon}
[X]:\varepsilon < \zeta\}$, abusing notation we assume that 
$\langle X \rangle \char 94 \langle{\Cal F}_{{\bar \alpha},\varepsilon}[X]:
\varepsilon < \zeta \rangle$ is an
increasing sequence of structures and for $\varepsilon = 0$ we stipulate 
${\Cal F}_{\bar \alpha,\varepsilon -1}[X]=X$
\sn
\item "{$(b)$}" $f_1 \le_{\text{pr}} f_2$ \ub{iff} Dom$(f_1) \subseteq
\text{ Dom}(f_2)$ and $\beta \in \text{ Dom}(f_1) \Rightarrow 
\dsize \bigvee_{\varepsilon \in [-1,\zeta)} \, {\Cal F}_{{\bar \alpha},\varepsilon}
[X] \models f_1(\beta) \le_{\text{pr}} f_2(\beta)$
\sn
\item "{$(c)$}"  $f_1 \le f_2$ iff
{\roster
\itemitem{ $(i)$ }   Dom$(f_1) \subseteq \text{ Dom}(f_2)$
\sn
\itemitem{ $(ii)$ }  $\beta \in \text{ Dom}(f_1) \Rightarrow 
\dsize \bigvee_{\varepsilon \in [-1,\zeta)} \, {\Cal F}_{{\bar \alpha},\varepsilon}
[X] \models f_1(\beta) \le f_2(\beta)$
\sn
\itemitem{ $(iii)$ }  the set $\{\beta \in \text{ Dom}(f_1):
\dsize \bigvee_{\varepsilon \in [-1,\zeta]} \, {\Cal F}_{{\bar \alpha},\varepsilon}
[X] \models f_1(\beta) \le f_2(\beta) \wedge \neg 
[f_1(\beta) \le_{\text{pr}} f_2(\beta)]\}$ is finite
\sn
\endroster}
\item "{$(d)$}"   $f_1 \le_{\text{apr}} f_2$ \ub{iff} 
{\roster
\itemitem{ $(i)$ }  $f_1 \le f_2$
\sn
\itemitem{ $(ii)$ }  Dom$(f_1) = \text{ Dom}(f_2)$
\sn
\itemitem{ $(iii)$ }  for all but finitely many $\beta \in \text{ Dom}
(f_1)$ we have $f_1(\beta) = f_2(\beta)$ and for the rest
$\dsize \bigvee_{\varepsilon \in [-1,\zeta)} f_1(\beta) \le_{\text{apr}}
f_2(\beta)$.
\endroster}
\ermn
3) A trunk controller ${\Cal F}$ is $\aleph_1$-complete if $({\Cal F},
\le_{\text{pr}})$ is $\aleph_1$-complete. \nl
4) A trunk controller ${\Cal F}$ is \ub{based} on $\langle X_\beta:\beta < \alpha
\rangle$ \ub{if}:
\mr
\item "{$(a)$}"  each $X_\beta$ is a trunk controller
\sn
\item "{$(b)$}"  $f \in {\Cal F} \Rightarrow \text{ Dom}(f) 
\subseteq \alpha$ and $f \in {\Cal F} \and
\beta \in \text{ Dom}(f) \Rightarrow f(\beta) \in X_\beta$
\sn
\item "{$(c)$}"  if $f_1,f_2 \in {\Cal F}$ has domain $\subseteq \alpha$,
\ub{then} $f_1 \le_{\text{pr}} f_2$ iff $\text{ Dom}(f_1) \subseteq
\text{ Dom}(f_2) \and (\forall \beta \in \text{ Dom}(f_1))[X_\beta \models
f_1(\beta) \le_{\text{pr}} f_2(\beta)]$
\sn 
\item "{$(d)$}"  $f_1 \le f_2$ \ub{iff}
{\roster
\itemitem{ $(i)$ }   Dom$(f_1) \subseteq \text{ Dom}(f_2)$,
\sn
\itemitem{ $(ii)$ }  $\beta \in \text{ Dom}(f_1) \Rightarrow 
X_\beta \models f_1(\beta) \le f_2(\beta)$ and 
\sn
\itemitem{ $(iii)$ }  $\{\beta \in \text{ Dom}(f_1):X_\beta \models f_1(\beta) 
\nleq_{\text{pr}} f_2(\beta)\}$ is finite
\endroster}
\item "{$(e)$}"  if $f,g \in {\Cal F}$, Dom$(f) \subseteq \beta < \alpha,g 
\restriction \beta \le_x f$, \ub{then} $f \cup (g \restriction [\beta,
\alpha]) \in {\Cal F}$ is $\le_x$-lub of $\{f,g\}$ for
$x \in \{us,pr\}$, see \scite{it.1b}; also if $f_n \subseteq f_{n+1} \in {\Cal F}$ for 
$n < \omega$ and $\dbcu_n f_n \in {\Cal F}$ then $\dbcu_n f_n$ is a
$\le_{pr}$-lub of $\{f_n:n < \omega\}$ and is a $\le$-upper bound of
$\{f_n:n < \omega\}$ 
\sn
\item "{$(f)$}"  $f_1 \le_{\text{apr}} f_2$ \ub{iff}
{\roster
\itemitem{ $(i)$ }  Dom$(f_1) = \text{ Dom}(f_2)$,
\sn
\itemitem{ $(ii)$ }  $\beta \in \text{ Dom}(f_1) \Rightarrow X_\beta
\models f_1(\beta) \le f_2(\beta)$ 
\sn
\itemitem{ $(iii)$ }  the set $\{\beta \in \text{ Dom}(f_1):f_1(\beta) \ne
f_2(\beta)\}$ is finite and for those $\beta$'s, $X_\beta \models
f_1(\beta) \le_{\text{apr}} f_2(\beta)$
\endroster}
\sn
\item "{$(g)$}"  if $f \in {\Cal F},\beta < \alpha$ \ub{then} $f
\restriction \beta \in {\Cal F}$.
\ermn
5) We say a trunk controller ${\Cal F}$ is \ub{full} over $\langle X_\beta:\beta
< \alpha \rangle$ or \ub{fully based} on it if: 
\mr
\item "{$(\alpha)$}"  ${\Cal F}$ is based on $\langle X_\beta:\beta < \alpha
\rangle$
\sn
\item "{$(\beta)$}"  whenever $f$ is a function with 
domain a countable subset of $\alpha$ such that $\beta \in \text{ Dom}(f)
\Rightarrow f(\beta) \in X_\beta$ then $f \in {\Cal F}$.
\ermn
6) We say a trunk controller ${\Cal F}$ is finitely based over
$\langle X_\beta:\beta < \alpha \rangle$ if:
\mr
\item "{$(\alpha)$}"  ${\Cal F}$ is based on $\langle X_\alpha:\beta <
\alpha \rangle$
\sn
\item "{$(\beta)$}"  $0 \in X_\alpha$ minimal
\sn
\item "{$(\gamma)$}"  $f \in {\Cal F}$ iff $f$ is a function with
domain a countable subset of $\alpha$ and $\{\beta \in \text{
Dom}(f):\neg(0 \le_{pr} f(\beta)\}$ is finite.
\ermn
7) We say ${\Cal F}$ is the trivial trunk controller if: its set of
elements is ${\Cal H}(\aleph_0)$ and $\le = \le_p = \le_{apr}$ are the
equality on ${\Cal H}(\aleph_0)$. \nl
8) We say ${\Cal F}$ is transparent if $p_0 \le_{pr} p_1 \and p_0
\le_{pr} p_2 \Rightarrow (\exists p_3)(p_1 \le_{pr} p_3 \and p_2
\le_{pr} p_3)$. \nl
9) In part (4), (5), replacing $\langle X_\beta:\beta < \alpha \rangle$
by $\alpha$ means ``for some $\langle X_\beta:\beta < \alpha
\rangle$". \nl
10) In part (4), for $\beta < \alpha$ let ${\Cal F}^{[\beta]} = {\Cal F}[\beta]$ be
$X_\beta$, (clearly normally uniquely defined).
\enddefinition
\bigskip

\definition{\stag{it.1a} Definition}  
1) We say a trunk controller ${\Cal F}$ is \ub{simple} if: for any
sequence $\langle y_\varepsilon:\varepsilon < \omega_2 \rangle$ with
$y_\varepsilon \in {\Cal F}$ for $\varepsilon < \omega_2$ for some club $E$
of $\omega_2$ and pressing down $h:E \rightarrow \omega_2$ we have:
for any $\varepsilon < \zeta$ from $E$ of cofinality $\aleph_1,h(\varepsilon) = 
h(\zeta) \Rightarrow y_\varepsilon,y_\zeta$ have a common 
$\le^{\Cal F}_{\text{pr}}$-upper bound.
\nl
2) We say the trunk controller ${\Cal F}$ is semi simple if: for any
$y_\varepsilon \in {\Cal F}$ for $\varepsilon < \omega_2$ for some 
$\varepsilon < \zeta < \omega_2$, there is a common $\le_{\text{pr}}$-upper bound of
$y_\varepsilon,y_\zeta$. \nl
3) We say the trunk controller ${\Cal F}$ is semi simply based on
$\langle X_\beta:\beta < \alpha^* \rangle$ \ub{if} it is based on it,
$X_0$ is semi simple and every $X_{1 + \beta}$ is simple. 
We say ``simply based" if also $X_0$ is simple. \nl
4) If ${\Cal F}$ is a trunk controller based on $\langle
X_\gamma:\gamma < \alpha \rangle$ and $\beta \le \alpha$ then let
${\Cal F} \restriction \beta = \{f \in {\Cal F}:\text{Dom}(f)
\subseteq \beta\}$.
\enddefinition
\bigskip

\proclaim{\stag{it.1c} Claim}   Suppose that $\bar{\Cal F} = \langle
{\Cal F}_\beta:\beta < \alpha^* \rangle$ is a sequence of trunk controllers.
\nl
1) There is a unique trunk controller ${\Cal F}$ fully based on
$\bar{\Cal F}$. \nl
2) Assume CH.  If ${\Cal F}_0$ is a semi simple trunk controller and each 
${\Cal F}_{1+ \beta}$ is a simple trunk controller and ${\Cal F}$ is
semi simply based on $\bar{\Cal F}$ \ub{then} the ${\Cal F}$ from part 
(1) is semi simple. \nl
3) In part (2) if also ${\Cal F}_0$ is simple, \ub{then} ${\Cal F}$ is simple. \nl
4) For every trunk controller $X$ and $\zeta$ and $\bar \alpha =
\langle \alpha_\varepsilon:\varepsilon \le \zeta \rangle,{\Cal
F}_{{\bar \alpha},\zeta}[X]$, is a well defined trunk controller,
simple if $X$ is simple. \nl
5) If each ${\Cal F}_\beta$ is $\aleph_1$-complete, \ub{then} in part
(1) also ${\Cal F}$ is $\aleph_1$-complete.  Similarly in part (4), if
$X$ is $\aleph_1$-complete, \ub{then} ${\Cal F}_{\alpha,\zeta}[X]$ is
$\aleph_1$-complete.
\endproclaim
\bigskip

\proclaim{\stag{it1.d} Claim}  1) If ${\Cal F}$ is a trunk controller
based on $\langle X_\gamma:\gamma < \alpha \rangle$ and $\beta \le
\alpha$ \ub{then} ${\Cal F} \restriction \beta$ is a trunk controller
based on $\langle X_\gamma:\gamma < \beta \rangle$. \nl
Similarly for ``fully based". \nl
2) If ${\Cal F}$ is simple, \ub{then} ${\Cal F}$ is semi-simple.
\endproclaim  
\bigskip

\demo{\stag{it.1b} Convention}:   Let $\le^{\Cal F}_{\text{us}} = 
\le^{\Cal F}$ and we write $\le^{\Cal F}_x$ for $x$ varying on 
$\{\text{us,pr,apr}\}$.  Similarly in Definition \scite{it.2}.
\enddemo
\bigskip

\definition{\stag{it.2} Definition}   A ${\Cal F}$-forcing notion $\Bbb Q$
means:  a tuple $(Q,\le,\le_{\text{pr}},\le_{\text{apr}},\text{val})$ 
(we may put superscript $\Bbb Q$ to clarify) satisfying:
\mr
\item "{$(a)$}"  $Q$ is a nonempty set (- the set of conditions) (we may
write $p \in \Bbb Q$ instead of $p \in Q$ and say $\Bbb Q$-names, etc.
and $Q$ instead $(Q,\le)$)
\sn
\item "{$(b)$}"  $\le,\le_{\text{pr}},\le_{\text{apr}}$ are quasi-orders on 
$\Bbb Q$ (called the usual, the pure and the apure)
\sn
\item "{$(c)$}"  $\le_{\text{pr}} \subseteq \le$
and $\le_{\text{apr}} \subseteq \le$ 
\sn
\item "{$(d)$}"  val is a function from $Q$ to ${\Cal F}$, a trunk
controller
\sn
\item "{$(e)$}"  $\Bbb Q \models p \le_x q \Rightarrow {\Cal F} \models
\text{ val}^{\Bbb Q}(p) \le_x \text{ val}^{\Bbb Q}(q)$ for $x = 
\text{ us,pr,apr}$
\sn
\item "{$(f)$}"  if $p_0 \le p_2$ \ub{then} for some $p_1$ we have $p_0
\le_{\text{pr}} p_1 \le_{\text{apr}} p_2$ and we can (i.e. ${\Cal F}$
give the additional information how to) compute val$(p_1)$ from
$(\text{val}(p_0),\text{val}(p_2))$; we write val$(p_1) = \text{
inter}_{\Cal F}(\text{val}(p_0),\text{val}(p_2))$; so we are abusing
our notation as we should have expanded ${\Cal F}$ by the function inter.
\endroster
\enddefinition
\bigskip

\definition{\stag{it.2a} Definition}  1) $\Bbb Q$, an ${\Cal
F}$-forcing is very clear (as an ${\Cal F}$-forcing) or a very clear
${\Cal F}$-forcing if:
\mr
\item "{$(*)$}"  if $p_0,p_1 \in \Bbb Q$ and val$^{\Bbb Q}(p_0)$,val$^{\Bbb Q}(p_0)$ has
a common $\le^{\Cal F}_{\text{pr}}$-upper bound $y$ \ub{then} for some 
$q \in \Bbb Q$ we have $p_0 \le_{\text{pr}} q,p_1 \le_{\text{pr}} q$ and
val$^{\Bbb Q}(q) = y$.
\ermn
2) ${\Bbb Q}$ is basic when: if $p_0 \le p_2$ \ub{then} for some $p_1$ we
have $p_0 \le_{pr} p_1 \le_{apr} p_2$ 
and val$^{\Bbb Q}(p_1) =
\text{ val}^{\Bbb Q}(p_0)$.
\nl
3) Let $\Bbb Q$ be an ${\Cal F}$-forcing, it is \ub{straight}, or
${\Cal F}$-straight when: \nl

if $p_1 \le_{apr} q_1,p_1 \le_{pr} p_2$ and
$p_2,q_1$ are compatible, \ub{then} there is $q_2$ such that $q_1 \le
q_2,p_2 \le_{apr} q_2$ which is a $\le^{\Bbb Q}$-lub of $p_2,q_1$ and
val$^{\Bbb Q}(q_2)$ can be computed from $\langle \text{val}^{\Bbb
Q}(p_1),\text{val}^{\Bbb Q}(p_2),\text{val}^{\Bbb Q}(q_1)\rangle$, and
we stipulate that this computation is a function which is part of the trunk
controller.  We call it amal$_{\Cal F}(-,-)$  (the
point is that when we iterate over $\bold V$ this function will be
$\bold V$).  If $p_2,q_1$ are incompatible, we use $q_2 = q_1$.  [Used
in \scite{ct.3a}, \scite{ct.2c}].  \nl
4) An ${\Cal F}$-forcing $\Bbb Q$ is called clear or ${\Cal F}$-clear when: if
$p \le_{pr} p_1,p \le_{pr} p_2$ and $p_1,p_2$ are
$\le_{pr}$-compatible then they have a common $\le_{pr}$-upper bound
$q$ with val$^{\Bbb Q}(q)$ computable (see (3)) from $\langle
\text{val}^{\Bbb Q}(p),\text{val}^{\Bbb Q}(p_1),\text{val}^{\Bbb
Q}(p_2)\rangle$ and we denote it by pramal$(-,-,-)$.  
[Used in \scite{ct.3a}, \scite{ct.2c}; this is not
part (1)]. \nl
5) An ${\Cal F}$ forcing $\Bbb Q$ is weakly \footnote{this seems not
to imply clear} clear when: 
\mr
\item  If $p_0,p_1 \in Q$ and val$^{\Bbb Q}(p_0)$,
val$(p_1)$ are $\le_{\text{pr}}$-compatible in ${\Cal F}$, \ub{then}
$p_0,p_1$ are $\le_{pr}$-compatible. 
\ermn
6) We say $\Bbb Q$ is transparent \footnote{This simplifies quite a
number of definitions below.  Of course, instead for every $y_3$ it is
enough to have one such $y_3 = y_3(\text{val}^{\Bbb Q}(p_1):\ell <
3)$, this function being part of ${\Cal F}$}
(or ${\Cal F}$-transparent) when: if $p_0
\le_{pr} p_1,p_0 \le_{pr} p_2$,val$^{\Bbb Q}(p_1) \le_{pr}
y_3$,val$^{\Bbb Q}(p_2) \le_{pr} y_3 \in {\Cal F}$, \ub{then} there is
$p_3 \in \Bbb Q$ such that $p_1 \le_{pr} p_3,p_2 \le_{pr} p_3$ and
val$^{\Bbb Q}(p_3) = y_3$. 
\enddefinition
\bn
\ub{\stag{it.3a} Discussion}:  1) We may involve measures; this, of course,
also complicates the iteration.  Interest not clear at present. \nl
2) We can consider some variants: if $p \le_{\text{pr}} q_\ell$ for $\ell=1,2$, do we
just ask $q_1,q_2$ compatible?  Does it suffice to demand ``val$(q_1)$,
val$(q_2)$ are $\le^{\Cal F}_{\text{pr}}$-compatible"?  The example
satisfies this but the general theorems do not need it.
\bigskip

\definition{\stag{it.4} Definition}  Let ${\Cal F}$ be a trunk controller
based on $\langle X_\beta:\beta < \alpha^* \rangle$.  We define by induction
on the ordinal $\alpha \le \alpha^*$ what is an ${\Cal F}$-iteration
$\bar{\Bbb Q}$ of length $\alpha$ and what is Lim$_{\Cal F}(\bar{\Bbb Q})$
\mr
\item "{$(a)$}"  $\bar{\Bbb Q}$ is an ${\Cal F}$-iteration of length $\alpha$ if
{\roster
\itemitem{ $(\alpha)$ }  $\bar{\Bbb Q} = \langle \Bbb P_\beta,
{\underset\tilde {}\to {\Bbb Q}_\beta}:\beta < \alpha \rangle$
\sn
\itemitem{ $(\beta)$ }  if $\beta < \alpha$ then $\bar{\Bbb Q} \restriction
\beta$ is an ${\Cal F}$-iteration of length $\beta$
\sn
\itemitem{ $(\gamma)$ }  if $\alpha = \beta +1$ then $\Bbb P_\beta =
\text{ Lim}_{\Cal F}(\bar{\Bbb Q} \restriction \beta)$
\sn
\itemitem{ $(\delta)$ }  if $\alpha = \beta +1$ then 
${\underset\tilde {}\to {\Bbb Q}_\beta}$ is a ${\Bbb P}_\beta$-name of an
$X_\beta$-forcing notion
\endroster}
\item "{$(b)$}"  for $\bar{\Bbb Q} = \langle \Bbb P_\beta,
{\underset\tilde {}\to {\Bbb Q}_\beta}:\beta < \alpha \rangle$ an
${\Cal F}$-iteration of length $\alpha$ we define the ${\Cal
F}$-forcing notion $\Bbb P_\alpha =
\text{ Lim}_{\Cal F}(\bar{\Bbb Q})$ as follows (see \scite{it.5a}):
{\roster
\itemitem{ $(\alpha)$ }  the set of elements of $\Bbb P_\alpha$ is the set
of $p$ such that for some $f \in {\Cal F}$ we have
\sn
\itemitem{ ${{}}$ } $\qquad (i) \quad p$ a function
\sn
\itemitem{ ${{}}$ }  $\qquad (ii) \quad$ Dom$(p) = \text{ Dom}(f)$, so
is a countable subset of $\alpha$
\sn
\itemitem{ ${{}}$ }  $\qquad (iii) \quad$ for $\beta \in \text{ Dom}(p)$ we
have $p(\beta)$ is a $\Bbb P_\beta$-name of \nl

$\qquad \qquad \qquad$ a member of ${\underset\tilde {}\to {\Bbb Q}_\beta}$
\sn
\itemitem{ ${{}}$ }  $\qquad (iv) \quad \Vdash_{P_\beta} ``\text{val}
^{\underset\tilde {}\to {\Bbb Q}_\beta}(p(\beta)) = f(\beta)"$ for
$\beta \in \text{ Dom}(p)$
\sn
\itemitem{ ${{}}$ }  $\qquad (v)$ \footnote{actually follows}  $\quad 
\beta < \alpha \Rightarrow p \restriction \beta \in \Bbb P_\beta$. \nl
Clearly $f$ is unique and we call it $f^p$ or $f[p]$
\endroster}
\item "{$(\beta)$}"  $\le^{{\Bbb P}_\alpha}_{\text{pr}}$ is defined by: \nl
$p \le^{{\Bbb P}_\alpha}_{\text{pr}} q$ \ub{iff} $(p,q \in P_\alpha$ and)
Dom$(p) \subseteq \text{ Dom}(q)$ and \nl

$\qquad \qquad \quad \beta \in 
\text{ Dom}(p) \Rightarrow q \restriction \beta \Vdash_{{\Bbb P}_\beta} ``p(\beta) 
\le^{\underset\tilde {}\to {\Bbb Q}_\beta}_{\text{pr}} q(\beta)"$
\sn
\item "{$(\gamma)$}"  $\le^{{\Bbb P}_\alpha}$ is defined by: \nl
$p \le^{{\Bbb P}_\alpha} q$ \ub{iff} $(p,q \in P_\alpha$ and)
{\roster
\itemitem{ $(i)$ }  Dom$(p) \subseteq \text{ Dom}(q)$ and 
\sn
\itemitem{ $(ii)$ }  $\beta \in \text{ Dom}(p) \Rightarrow q 
\restriction \beta \Vdash_{{\Bbb P}_\beta} ``p(\beta) 
\le^{\underset\tilde {}\to {\Bbb Q}_\beta} q(\beta)"$ and
\sn
\itemitem{ $(iii)$ }   for some finite
$w \subseteq \text{ Dom}(\beta)$ we have $\beta \in \text{ Dom}(p) \backslash
w \Rightarrow q \restriction \beta \Vdash_{{\Bbb P}_\beta} ``p(\beta) 
\le^{\underset\tilde {}\to {\Bbb Q}_\beta}_{\text{pr}} q(\beta)"$
\endroster}
\item "{$(\delta)$}"  $\le^{{\Bbb P}_\alpha}_{\text{apr}}$ is 
defined by \nl
$p \le^{{\Bbb P}_\alpha}_{\text{apr}} q$ \ub{iff} $(p,q \in {\Bbb P}_\alpha$ 
and)
{\roster
\itemitem{ $(i)$ }   Dom$(p) = \text{ Dom}(q)$ and 
\sn
\itemitem{ $(ii)$ }  $p \le q$ and
\sn
\itemitem{ $(iii)$ }   $\beta \in \text{ Dom}(p) \Rightarrow q 
\restriction \beta \Vdash_{{\Bbb P}_\beta} ``p(\beta) \le_{\text{apr}} q(\beta)$ in 
${\underset\tilde {}\to {\Bbb Q}_\beta}"$
\sn
\itemitem{ $(iv)$ }  for all but finitely \footnote{actually follows
by clause (ii), when (as usual) $\Vdash_{{\Bbb P}_\beta}$ ``if $p'
\le^{{\Bbb Q}_\beta}_{pr} q'$ and $p' \le^{{\Bbb Q}_\beta}_{apr} q'$
then $p' = q'$"} many $\beta \in \text{ Dom}(p)$
we have $p(\beta) = q(\beta)$.
\endroster}
\endroster
\enddefinition
\bn
\ub{\stag{it.5} Convention}:  If ${\Cal F}$ and $\bar{\Bbb Q} = \langle
{\Bbb P}_\beta,{\underset\tilde {}\to {\Bbb Q}_\beta}:\beta < \alpha \rangle$
are as in \scite{it.4} then $\Bbb P_\alpha = \text{Lim}_{\Cal F}
(\bar{\Bbb Q})$.
\bigskip

\proclaim{\stag{it.5a} Claim}  If $\bar{\Bbb Q}$ is an ${\Cal
F}$-iteration and $\beta \le \ell g(\bar{\Bbb Q})$, \ub{then} $\Bbb
P_\beta$ is a $({\Cal F} \restriction \beta)$-forcing.
\endproclaim
\bigskip

\demo{Proof}  Straight, the least obvious is checking clause (f) of
Definition \scite{it.2}. 
\bn
\ub{Clause (f)}:  Assume $p_0 \le p_2$ and we shall define $p_1$.  Let
Dom$(p_1) = \text{ Dom}(p_2)$ and let the finite $w \subseteq \text{
Dom}(f_1)$ be as in $(\gamma)(iii)$ of \scite{it.4}(b), and we choose
$p_1(\alpha)$ for $\alpha \in \text{ Dom}(p_1)$ a follows.  If $\alpha
\in \text{ Dom}(p_2) \backslash \text{ Dom}(p_0)$ we let $p_1(\alpha)
= p_2(\alpha)$, and if $\alpha \in \text{ Dom}(p_0)$ and $\alpha
\notin w$ we let $p_1(\alpha) = p_2(\alpha)$.  If $\alpha \in w$ by
Definition \scite{it.2}(f), we know that $p_2 \restriction \alpha
\Vdash_{{\Bbb P}_\alpha} ``p_0(\alpha) \le p_2(\alpha)"$ hence $(\exists p)[p_0(\alpha))
\le^{{\Bbb Q}_\alpha}_{pr} p \le^{\underset\tilde {}\to {\Bbb
Q}_\alpha}_{apr} p_2(\alpha) \and \text{ val}^{\underset\tilde {}\to
{\Bbb Q}_\beta}(p) = \text{ interp}_{{\Cal F}[\alpha]}
(\text{val}_{{\Cal F}[\alpha]}(p_0(\alpha)),\text{val}^{{\Cal F}[\alpha]}(p_2(\alpha)))]$
and choose $p(\alpha)$ as such $p$.  Now check.
\enddemo
\bigskip

\proclaim{\stag{it.6} Claim}  Assume ${\Cal F}$ is a trunk controller based
on $\bar X = \langle X_\beta:\beta < \alpha^* \rangle$; moreover is fully based on
$\bar X$ (see Definition \scite{it.1}(5)) and $\alpha \le \alpha^*$ and 
$\bar{\Bbb Q}$ is an ${\Cal F}$-iteration of length $\alpha$ and $\gamma \le 
\beta \le \alpha$. \nl 
1) If $p \in \Bbb P_\beta$ \ub{then} $p \restriction \gamma \in \Bbb
P_\gamma$ and $\Bbb P_\beta \models ``p \restriction 
\gamma \le_{\text{pr}} p"$.
\nl
2) $\Bbb P_\gamma \subseteq \Bbb P_\beta$, i.e. $p \in \Bbb P_\gamma
\Rightarrow p \in \Bbb P_\beta$ and $\le^{{\Bbb P}_\gamma}_x =
\le^{{\Bbb P}_\beta}_x \restriction \Bbb P_\gamma$ (see convention
\scite{it.1b}). \nl
3) If $p \in \Bbb P_\beta,x \in \{\text{us,pr,apr}\},
p \restriction \gamma \le^{{\Bbb P}_\gamma}_x q \in
\Bbb P_\gamma$ and $r = q \cup (p \restriction (\gamma,\beta))$, \ub{then} $r$
is $\le^{{\Bbb P}_\beta}_x$-lub of $\{p,q\}$ when $x \in \{us,pr\}$
and $p \le_{apr} r \and q \le_{pr} r$ when $x = apr$. \nl
4) $\Bbb P_\gamma \lessdot \Bbb P_\beta$. 
\endproclaim
\bigskip

\demo{Proof}  Straight, by induction on $\alpha$.
\enddemo
\bigskip

\proclaim{\stag{it.7} Claim}  Assume $\bar{\Bbb Q}$ is an ${\Cal
F}$-iteration of length $\alpha$ and $\Bbb P = \text{ Lim}(\bar{\Bbb
Q})$. \nl
1) The property ``very clear", \scite{it.2a}(1) is preserved, i.e. if
each ${\underset\tilde {}\to {\Bbb Q}_\beta} \, (\beta < \alpha)$ is
very clear, then so is $\Bbb P$. \nl
2) The property ``straight", \scite{it.2a}(3) is preserved. \nl
3) The property ``clear", \scite{it.2a}(4) is preserved.
\endproclaim
\bigskip

\demo{Proof}  Straight.
\enddemo
\bigskip

\proclaim{\stag{it.7a} Claim}  1) For an ${\Cal F}$-forcing $\Bbb Q$;
very clear implies clear and implies weakly clear. \nl
2) Assume $\Bbb Q$ is a ${\Cal F}$-forcing, $\Bbb Q$ is weakly clear
(\scite{it.2a}(5)), and ${\Cal F}$ is semi-simple, \ub{then} $\Bbb Q$
and even $(\Bbb Q,\le_{pr})$ satisfies the $\aleph_2$-c.c. \nl
3) Assume $\Bbb Q$ is an ${\Cal F}$-forcing, $\Bbb Q$ is weakly clear
and ${\Cal F}$ is simple, \ub{then} $\Bbb Q$ and even $(\Bbb
Q,\le_{pr})$ satisfies the regressive $\aleph_2$-c.c.
\endproclaim
\bigskip

\demo{Proof}  Straight.
\enddemo
\bigskip

\proclaim{\stag{it.8} Claim}  Assume ${\Cal F}$ is a trunk controller based
on $\alpha^*$. \nl
0) The empty sequence is an ${\Cal F}$-iteration. \nl 
1) If $\bar{\Bbb Q}$ is an ${\Cal F}$-iteration of length $\alpha,\alpha + 1 
\le \alpha^*$ and $\underset\tilde {}\to {\Bbb Q}$ is a ${\Bbb P}_\alpha$-name
of an $X_\alpha$-forcing notion, \ub{then} there is a 
${\Cal F}$-iteration $\bar{\Bbb Q}'$ of
length $\alpha + 1$ such that $\bar{\Bbb Q}' \restriction \alpha = 
\bar{\Bbb Q}$ and ${\underset\tilde {}\to {\Bbb Q}'_\alpha} = 
\underset\tilde {}\to {\Bbb Q}$ that is $\bar{\Bbb Q} \char 94 
\langle \text{Lim}_{\Cal F}(\bar{\Bbb Q}),\underset\tilde {}\to {\Bbb Q}
\rangle$ is an ${\Cal F}$-iteration. \nl
2) If $\bar{\Bbb Q} = \langle {\Bbb P}_\beta,
{\underset\tilde {}\to {\Bbb Q}_\beta}:\beta <
\alpha \rangle$ and $\alpha$ is a limit ordinal and $\bar{\Bbb Q} 
\restriction \beta$ is an ${\Cal F}$-iteration for every $\beta < \alpha$ \ub{then}
$\bar{\Bbb Q}$ is an ${\Cal F}$-iteration.  \nl
3) For any function $\bold F$ and ordinal $\alpha \le \alpha^*$ there is a
unique ${\Cal F}$-iteration $\bar{\Bbb Q}$ such that:
\mr
\item "{$(\alpha)$}"  $\ell g(\bar{\Bbb Q}) \le \alpha$
\sn
\item "{$(\beta)$}"  $\beta < \ell g(\bar{\Bbb Q}) \Rightarrow
{\underset\tilde {}\to {\Bbb Q}_\beta} = \bold F(\bar{\Bbb Q} \restriction \beta)$
\sn
\item "{$(\gamma)$}"  if $\ell g(\bar{\Bbb Q}) < \alpha$ \ub{then}
$\bold F(\bar{\Bbb Q})$ is not a $(\text{Lim}_{\Cal F}(\bar{\Bbb Q}))$-name
of an ${\Cal F}^{[\beta]}_\beta$-forcing.
\endroster
\endproclaim
\bigskip

\demo{Proof}  Straight.
\enddemo
\bn
Not really necessary, but natural and aesthetic, is
\proclaim{\stag{it.9} Claim}  Associativity holds, that is assume
\mr
\item "{$(a)$}"  $\bar{\Bbb Q} = \langle \Bbb P_\beta,
{\underset\tilde {}\to {\Bbb Q}_\beta}:\beta < \alpha^* \rangle$ is an
${\Cal F}$-iteration so $\Bbb P_{\alpha^*} = \text{ Lim}_{\Cal F}(\bar{\Bbb Q})$
\sn
\item "{$(b)$}"  $\langle \alpha_\varepsilon:\varepsilon \le \varepsilon^*
\rangle$ is increasing continuous, $\alpha_0 = 0,\alpha_{\varepsilon^*} =
\alpha^*$
\sn
\item "{$(c)$}"  for $\gamma \le \beta \le \alpha^*$ we define $\Bbb P_\beta/
\Bbb P_\gamma$, an ${\Cal F}$-forcing, naturally: it is a $\Bbb P_\gamma$-name
so for $G_\gamma \subseteq \Bbb P_\gamma$ generic over $\bold V$ its interpretation is:
\block
set of elements is $\{p \in \Bbb P_\beta:\text{Dom}(p) \subseteq [\gamma,
\beta)\}$ \nl
val:  inherited from $\Bbb P_\beta$ \nl
$\le_x: p \le_x q$ iff for some $r \in G_\gamma$ we have $\Bbb P_\beta
\models r \cup p \le_x r \cup q$
\endblock
\sn
\item "{$(d)$}"  let ${\Cal F}' = \{f:\text{for some } g \in {\Cal F},f$ is
a function with domain $\{\varepsilon < \varepsilon^*:\text{Dom}(g) \cap
[\alpha_\varepsilon,\alpha_{\varepsilon +1}) \ne \emptyset\}$ and
$\varepsilon \in \text{ Dom}(f) \Rightarrow f(\varepsilon) = g \restriction
[\alpha_\varepsilon,\alpha_{\varepsilon +1})\}$, the orders of ${\Cal F}'$
are natural.
\ermn
\ub{Then} we can find an ${\Cal F}'$-iteration $\bar{\Bbb Q}' =
\langle \Bbb P'_\varepsilon,{\underset\tilde {}\to {\Bbb Q}'_\varepsilon}:
\varepsilon < \varepsilon^* \rangle$ and 
$\langle F_\varepsilon:\varepsilon \le \varepsilon^* \rangle$ such that
\mr
\item "{$(\alpha)$}"  $F_\varepsilon$ is an isomorphism from
$\Bbb P_{\alpha_\varepsilon}$ onto $\Bbb P'_\varepsilon$
\sn
\item "{$(\beta)$}"  when $\varepsilon < \varepsilon^*,F_\varepsilon$ maps
the $\Bbb P_{\alpha_\varepsilon}$-name $\Bbb P_{\alpha_{\varepsilon +1}}/
\Bbb P_{\alpha_\varepsilon}$ to the $\Bbb P'_\varepsilon$-name
${\underset\tilde {}\to {\Bbb Q}'_\varepsilon}$.
\endroster
\endproclaim
\bigskip

\remark{Remark}  For standard ${\Cal F}$ we can use ${\Cal F}' = {\Cal F}$.
\endremark
\bigskip

\demo{Proof}  Straight.
\enddemo
\bn
\ub{\stag{it.10} Discussion}:  May we be interested in nonstarndard ${\Cal F}$?
Maybe, if CH holds, 
$\eta_\alpha \in {}^{\omega_1}2$ for $\alpha < \alpha^*$ are
pairwise distinct, and we are interested in ${\Cal F}$ such that $\{f(\alpha):
f \in {\Cal F}\}$ depend on $\eta_\alpha$ continuously.  So for proving
$\aleph_2$-c.c. we have a stronger handle.
\bigskip

\definition{\stag{it.11} Definition}  1) $\Bbb Q$ has $(\theta,\sigma)$-pure 
decidability if:

if $p \in \Bbb Q$ and $p \Vdash_{\Bbb Q} 
``\underset\tilde {}\to \tau \in \theta"$, \ub{then} for some 
$A \subseteq \theta,|A| < \sigma$ and $q$ we have $p \le_{\text{pr}} 
q \in \Bbb Q$ and $q \Vdash ``\underset\tilde {}\to \tau = A"$. \nl
2) We write ``$\theta$-pure decidability" for ``$(\theta,\theta)$-pure
decidability".
\enddefinition
\newpage

\head {\S2 \\
Being ${\Cal F}$-Pseudo c.c.c. is preserved by ${\Cal F}$-iterations}
\endhead  \resetall 
\bn
The reader may concentrate on the transparent case, \scite{ct.2b},
\scite{ct.2d}, \scite{ct.3}(3) (but read the proofs of
\scite{ct.3}(1),(2)) and \scite{ct.4} and continue using the Knaster
explicit version. \nl
In the main cases, $\bold H$ disappears but $\Game_p$ is needed for the iteration.
\bigskip

\definition{\stag{ct.1} Definition}  1) Let ${\Cal F}$ be a trunk controller,
$\Bbb Q$ be an ${\Cal F}$-forcing notion.  We say that $\Bbb Q$ is 
${\Cal F}$-psc (${\Cal F}$-pseudo c.c.c. forcing in full) as witnessed by $\bold H$ if:

for every $p \in \Bbb Q$ in the following game $\Game_p = \Game_{p,\Bbb Q,
\bold H} = \Game_p[\Bbb Q,\bold H]$ between the player interpolator
and extendor which lasts $\omega_1$ moves,
the interpolator has a winning strategy.

In the $\zeta$-th move:
\mr
\item "{$\boxtimes$}"   the interpolator chooses a condition $p'_\zeta$ such
that $p \le_{\text{pr}} p'_\zeta$
and val$^{\Bbb Q}(p'_\zeta) = \bold H(\langle \text{val}^{\Bbb Q}
(p_\xi)$, val$^{\Bbb Q}(q_\xi)):\xi < \zeta \rangle)$ and then the extendor 
chooses $q_\zeta \in \Bbb Q$ such that $p'_\zeta \le
q_\zeta$ and lastly the interpolator chooses a condition $p_\zeta$ such that
$p'_\zeta \le_{\text{pr}} p_\zeta \le_{\text{apr}} q_\zeta$
and val$^{\Bbb Q}(p_\zeta) = \bold H(\langle(\text{val}^{\Bbb Q}(p_\xi)$, 
val$^{\Bbb Q}(q_\xi)):\xi < \zeta \rangle \char 94
\langle \text{val}^{\Bbb Q}(q_\zeta) \rangle)$.  \nl
[For future notation let $p^*_\zeta$ be $p'_\zeta$ and let $q_{-1} = p$].
\ermn
A play is won by the interpolator if:
\mr
\item "{$(\alpha)$}" for any stationary $A \subseteq \omega_1$, for
some $B \subseteq A$ we have
{\roster
\itemitem{ $(*)$ }  $B$ is stationary and $\bold H(\langle (\text{val}^{\Bbb Q}(p_\xi)$,
val$^{\Bbb Q}(q_\xi)):\varepsilon < \omega_1 \rangle \char 94 \langle
B \rangle) =1$
\endroster} 
\sn
\item "{$(\beta)$}"   if $B \subseteq \omega_1$ satisfy $(*)$ then: 
for every $\varepsilon < \zeta$ from $B$ we have:
$q_\varepsilon,q_\zeta$ are compatible in $\Bbb Q$ if
val$^{\Bbb Q}(q_\varepsilon)$, val$^{\Bbb Q}(q_\zeta)$ are compatible in
${\Cal F}$
\sn
\item "{$(\gamma)$}"  for $E = \omega_1$ or just $E$ a club of
$\omega_1$ computed from $\langle \text{val}^{\Bbb
Q}(p_\varepsilon),\text{val}^{\Bbb Q}(q_\varepsilon)):\varepsilon <
\omega_1 \rangle$ we have: if $\varepsilon < \zeta$ are from $E$,
$p_\varepsilon \le_{\text{pr}} q_\varepsilon$ and $p_\zeta \le_{\text{pr}}
q_\zeta$, \ub{then} $q_\varepsilon,q_\zeta$ has a common
$\le_{\text{pr}}$-upper bound $q$, with val$^{\Bbb Q}(q) = \bold H
(\varepsilon,\zeta,\langle (\text{val}^{\Bbb Q}(p_\xi),\text{val}^{\Bbb Q}
(q_\varepsilon):\xi \le \zeta \rangle)$. \nl
[Why not just $\varepsilon < \zeta$ from $B$?  For the iteration theorem \scite{ct.4}.]
\ermn
2) We define ``${\Bbb Q}$ is $({\Cal F},{\Cal P})$-psc" 
as above but at the end defining when a play is won by the
interpolator, we make the changes:
\mr
\item "{$(\alpha)'$}"  in every limit stage he has a legal move
\ub{or} the sequence $\langle(\text{val}^{\Bbb
Q}(p_\xi),\text{val}^{\Bbb Q}(q_\xi)):\xi < \zeta \rangle$ is not in
${\Cal P}$
\sn
\item "{$(\beta)'$}"  if $\langle(\text{val}^{\Bbb Q}(p_\zeta),
\text{val}^{\Bbb Q}(q_\zeta)):\zeta < \omega_1 
\rangle \in {\Cal P}$ \ub{then} for every stationary set
$A \in {\Cal P}$ of $\omega_1$, there is a stationary subset $B \in {\Cal P}$
as there. 
\ermn
3)  We can replace $\langle(\text{val}^{\Bbb Q}
(p_\varepsilon)$,val$^{\Bbb Q}(q_\varepsilon)):\varepsilon < \zeta \rangle$
by $\langle \text{val}^{\Bbb Q}(q_\varepsilon):-1 \le \varepsilon < \zeta
\rangle$. \nl
4) If we omit $\bold H$ and say \ub{bare}, this means that: we just 
omit the relevant demands on the interpolator in $\boxtimes$ 
and in $(*)$ of $(\alpha)$ of part (1), just requires
that $(\beta)$ and $(\gamma)$ holds: if we omit clause $(\gamma)$ of
(1) we say ``weakly" ${\Cal F}$-psc. \nl
5) We say $\bar{\Bbb Q}$ is ${\Cal F}$-straight if each
${\underset\tilde {}\to {\Bbb Q}_\beta}$ is ${\Cal
F}^{[\beta]}$-straight; similarly ${\Cal F}$-clear. \nl
6) We say $\bar{\Bbb Q}$ is semi-straight (or semi ${\Cal F}$-straight)
if each $\Bbb Q_{1 + \beta}$ is ${\Cal F}^{[\beta]}$-straight
similarly semi ${\Cal F}$-clear.
\enddefinition
\bigskip

\remark{Remark}  Note that \scite{mr.1}(1) speak actually on any semi
straight ${\Cal F}$-psc forcing $\Bbb P$.
\endremark 
\bigskip

\definition{\stag{ct.2} Definition}  1) We say 
$\bar{\Bbb Q}$ is an ${\Cal F}$-psc iteration as witnessed by
$\bar{\bold H}$ \ub{if}:
\mr
\item "{$(a)$}"  ${\Cal F}$ is a trunk controller, fully based on some
$\alpha' \ge \ell g(\bar{\Bbb Q})$
\sn
\item "{$(b)$}"  $\bar{\Bbb Q}$ is an ${\Cal F}$-iteration and
\sn
\item "{$(c)$}"   for every $\beta
< \ell g(\bar Q)$ we have $\Vdash_{{\Bbb P}_\beta} ``{\underset\tilde
{}\to {\Bbb Q}_\beta}$ is an $({\Cal F}^{[\beta]},\bold V)$-psc as
witnessed by $\bold H_\beta"$ and $\bar{\bold H} = \langle \bold
H_\beta:\beta < \ell g(\bar Q) \rangle$; \nl
note that $\bold H_\beta \in \bold V$ and is an object, not a $\Bbb
P_\beta$-name
\sn
\item "{$(d)$}"  let us stress that the witnesses for ``$\Bbb Q_{1 + \beta}$ is
straight/clear" are given in $\bold V$ (not as name).
\ermn
2) We say $\bar{\Bbb Q}$ is an essentially ${\Cal F}$-psc 
iteration as witnessed by $\bold H$ \ub{if}
\mr 
\item "{$(a)$}" it is an ${\Cal F}$-iteration 
\sn 
\item "{$(b)$}" $\bar{\bold H} = \langle \bar{\bold H}_\beta:\beta <
\alpha^* \rangle$ or $\bar{\bold H} = \langle \bold H_\beta:0 < \beta
< \alpha^* \rangle$ and $\alpha^* \ge \ell g(\bar{\Bbb Q})$ 
\sn 
\item "{$(c)$}" for every nonzero $\beta < \ell g(\bar{\Bbb Q})$ we have 
$$
\Vdash_{{\Bbb P}_\beta} ``{\underset\tilde {}\to {\Bbb Q}_\beta}
\text{ is an } ({\Cal F}^{[\beta]},\bold V)\text{-psc forcing notion as
witnessed by } \bold H_\beta" 
$$
\item "{$(d)$}"  ${\underset\tilde {}\to {\Bbb Q}_0}$ satisfies the c.c.c.
and $\le^{{\Bbb Q}_0}_{\text{pr}},\le^{{\Bbb Q}_0}_{\text{apr}}$ are
equality and val$^{{\Bbb Q}_0}$ is constantly 0. 
\ermn
3) We say $\bar{\Bbb Q}$ is a semi simple ${\Cal F}$-psc iteration if:
\mr
\item "{$(a)$}"  it is an ${\Cal F}$-iteration
\sn
\item "{$(b)$}"  for every non zero $\beta < \ell g(\bar{\Bbb Q})$ we have
$\Vdash_{{\Bbb P}_\beta} ``{\underset\tilde {}\to {\Bbb Q}_\beta}$ is a simple
$({\Cal F}^{[\beta]},\bold V)$-psc forcing notion".
\endroster
\enddefinition
\bigskip

\remark{\stag{it.2z} Remark}  No harm demanding
\mr
\item "{$(c)$}"  ${\underset\tilde {}\to {\Bbb Q}_0}$ satisfies the c.c.c.
and $\le^{{\Bbb Q}_0}_{\text{pr}}$ is equality, $\le^{{\Bbb Q}_0}_{\text{apr}}$ is 
$\le^{{\Bbb Q}_0}$ and val$^{{\Bbb Q}_0}$ is constantly 0.
\endroster
\endremark
\bigskip

\definition{\stag{ct.2a} Definition}   Let ${\Cal F}$ be at trunk
controller and $\Bbb Q$ an ${\Cal F}$-forcing. \nl
1) We say that $\Bbb Q$ is strongly ${\Cal F}$-psc forcing notion as
witnessed by $\bold H$ if for every $p \in \Bbb Q$ in the game
$\Game'_p = \Game'_{p,{\Bbb Q},\bold H} = \Game'_p[\Bbb Q,\bold H]$
the interpolator has a winning strategy, where the game is defined as
in \scite{ct.1} except that in addition we demand

$$
\varepsilon < \zeta \Rightarrow \Bbb Q \models p'_\varepsilon \le_{pr}
p'_\zeta.
$$
\mn
2) Adding the adjective ``semi" (in \scite{ct.1}(1) hence in (1) here) means
just that in clause $(\beta)$ we just ask for some $\varepsilon
< \zeta$ in $B$, the conditions $q_\varepsilon,q_\zeta$ are compatible
in $\Bbb Q$; so we call the games semi-$\Game_p$, semi-$\Game'_p$.
In Definition \scite{ct.2} and \scite{ct.2a}(2) adding ``semi" means,
$\Bbb Q_0$ satisfies the semi version, each $\Bbb Q_{1 + \alpha}$ the
regular one. \nl
3) Saying ``strongly" in Definition \scite{ct.2} means that the
demands are for each ${\underset\tilde {}\to {\Bbb Q}_\beta}$.
\enddefinition
\bigskip

\proclaim{\stag{ct.3a} Claim}  1) If in Definition \scite{ct.1},
$\Bbb Q$ is bare ${\Cal F}$-psc forcing and is ${\Cal F}$-clear or at
least straight (see Definition \scite{it.2a}) \ub{then} $\Bbb Q$ is
${\Cal F}$-psc (as witnessed by some $\bold H$, so $\bold H$ is
redundant), similarly with ${\Cal P}$. \nl
2) Assume
\mr
\item "{$(a)$}"  $\Bbb Q$ is a $\sigma$-centered forcing notion, i.e.
$\Bbb Q = \dbcu_{n < \omega} R_n$ each $R_n$ directed and for
simplicity may assume $n \ne m \Rightarrow R_m \cap R_m = \emptyset$
\sn
\item "{$(b)$}"  ${\Cal F}$ is standard
\sn
\item "{$(c)$}"  $\Bbb Q$ is defined by: 
$(Q^{\Bbb Q},\le^{\Bbb Q})$ is $Q$, \nl
$\le^{\Bbb Q}_{\text{pr}}$ is equality \nl
$\le^{\Bbb Q}_{\text{apr}}$ is $\le^{\Bbb Q}$ \nl
val$^{\Bbb Q}(q) = \text{ Min}\{n:q \in R_n\}$.
\ermn
\ub{Then} $\Bbb Q$ is a ${\Cal F}$-psc forcing as witnessed by a
trivial $\bold H$, and is purely proper. \nl
3) Assume
\mr
\item "{$(a)$}"  $\Bbb Q$ is a forcing notion, $\alpha^*$ is an ordinal,
$\bar h = \langle h_q:q \in \Bbb Q \rangle$ is such that each $h_q$ is a 
finite (partial) function from $\alpha^*$ to $\omega$ and: if
$h = h_{q_1} \cup h_{q_2}$ is a function then
$q_1,q_2$ has a common upper bound $q$ with $h_q = h_{q_1} \cup h_{q_2}$
\sn
\item "{$(b)$}"  ${\Cal F}$ is a trunk controller whose set of elements
include $\{h_q:q \in Q\}$ such that $h_{q_1} \subseteq h_{q_2} \Rightarrow
{\Cal F} \models h_{q_1} \le h_{q_2}$
\sn
\item "{$(c)$}"  $\Bbb Q$ is defined as in clause (c) above except that val$^{\Bbb Q}(q) =
h_q$.
\ermn
\ub{Then} $\Bbb Q$ is ${\Cal F}$-psc.
\endproclaim
\bigskip

\demo{Proof}  Straightforward.  In part (1), the ``$\Bbb Q$ is
a clear ${\Cal F}$-forcing" is used for clause $(\gamma)$ in
\scite{ct.1}(1); for using ``straight" note that in clause $(\gamma)$
also $p_\varepsilon \le_{pr} q_\varepsilon$ hold by the demands in
\scite{ct.1}$(\gamma)$, as well as $p_\varepsilon \le_{apr}
q_\varepsilon$ by $\boxtimes$ of \scite{ct.1}(1).  \hfill$\square_{\scite{ct.3a}}$\margincite{ct.3a}
\enddemo
\bigskip

\definition{\stag{ct.2b} Definition}  Let ${\Cal F}$ be a trunk
controller. \nl
1) ${\Cal F}$ satisfies the psc if:
\mr
\item "{$(*)$}"  if ${\Cal F} \models ``x \le_{pr} y_\varepsilon
\le_{apr} z_\varepsilon"$ for $\varepsilon < \omega_1$ and $A
\subseteq \omega_1$ is stationary then for some stationary $B
\subseteq A$ we have: \nl
if $\varepsilon < \zeta$ are from $B$ then $z_\varepsilon,z_\zeta$ has a
common upper bound $z$ such that:
{\roster
\itemitem{ $(\alpha)$ }  val$(z) = \text{ glue}_{{\Cal
F},\varepsilon,\zeta}(z_\varepsilon,z_\zeta,y_\varepsilon,y_\zeta)$
\sn
\itemitem{ $(\beta)$ }  if $y_\varepsilon \le_{pr}
z_\varepsilon,y_\zeta \le_{pr} z_\zeta$ then $z_\zeta \le_{pr} z$.
\endroster}
\ermn
2) ${\Cal F}$ satisfies the semi-psc if:
\mr
\item "{$(*)$}"  if ${\Cal F} \models ``x \le_{pr} y_\varepsilon
\le_{apr} z_\varepsilon"$ for $\varepsilon < \omega_1$, \ub{then} for
some $\varepsilon < \zeta < \omega_1,z_\varepsilon,z_\zeta$ has a
common upper bound $z$ such that $(\alpha) + (\beta)$ above holds.
\ermn
3) In parts (1), (2) we add the adjective ``continuous" if in $(*)$ we
add $\varepsilon < \zeta < \omega_1 \Rightarrow {\Cal F} \models
``y_\varepsilon \le_{pr} y_\zeta"$. 
We add Knaster if we replace stationary by unbounded (this is
alternative to semi$^7$, but not for iteration!)
\nl
4) We add the adjective ``finished" if we omit clause $(\gamma)$ in
\scite{ct.1}(1) and its variants. \nl
5) We say that an ${\Cal F}$-forcing $\Bbb Q$ is [Knaster] explicitly
[semi]
${\Cal F}$-psc forcing \ub{if}
\mr
\item "{$(a)$}"  $\Bbb Q$ is an ${\Cal F}$-forcing
\sn
\item "{$(b)$}"  ${\Cal F}$ satisfies the [Knaster][semi] psc
\sn
\item "{$(c)$}"  if $z = \text{ glue}_{{\Cal
F},\varepsilon,\zeta}(z',z'',y',y'')$ and $q',q'',p',p'' \in \Bbb Q$
and $p' \le_{pr} p'',p' \le_{apr} q',p'' \le_{apr} q'',y' = \text{
val}^{\Bbb Q}(p'),y'' = \text{ val}^{\Bbb Q}(p''),z' = \text{
val}^{\Bbb Q}(q'),z'' = \text{ val}^{\Bbb Q}(q'')$, \ub{then} there is
$q \in \Bbb Q$ such that:
{\roster
\itemitem{ $(\alpha)$ }  val$^{\Bbb Q}(q) = \text{ glue}_{\Cal
F}(z,z',z'',y',y'')$
\sn
\itemitem{ $(\beta)$ }  if $p' \le_{pr} q',p'' \le_{pr} q''$ then $q''
\le_{pr} q$.
\endroster}
\ermn
6) We say $\bar{\Bbb Q}$ is a Knaster explicit [semi] ${\Cal F}$-psc
\mr
\item "{$(a)$}"  ${\Cal F}$ is a trunk controller, fully based on some
$\alpha' \ge \ell g(\bar{\Bbb Q})$
\sn
\item "{$(b)$}"  $\bar{\Bbb Q}$ is an ${\Cal F}$-iteration
\sn
\item "{$(c)$}"  ${\Cal F}^{[0]}$ satisfies the Knaster [semi]
psc; $\aleph_1$-complete
\sn
\item "{$(d)$}"  ${\Cal F}^{[1 + \beta]}$ satisfies the [Knaster] psc
when $1 + \beta < \ell g(\bar Q);\aleph_1$-complete
\sn
\item "{$(e)$}"  $\Bbb Q_0$ is explicitly [Knaster][semi] ${\Cal
F}^{[0]}$-psc,$\le_{pr}-\aleph_1$-complete
\sn
\item "{$(f)$}"  $\Bbb Q_{1 + \beta}$ is (forced to be) an 
explicitly [Knaster] ${\Cal F}$-psc forcing
$\le_{pr}-\aleph_1$-complete.
\ermn
7) We (in (5), (6)) add continuous if so are the ${\Cal F}$'s.
\enddefinition
\bigskip

\proclaim{\stag{ct.2c} Claim}  1) Assume
\mr
\item "{$(a)$}"  ${\Cal F}$ is a [semi]-psc trunk controller
\sn
\item "{$(b)$}"  $\Bbb Q$ is ${\Cal F}$-forcing notion.
\ermn
\ub{Then} $\Bbb Q$ is a [semi] ${\Cal F}$-psc forcing notion.
\nl
2) Assume
\mr
\item "{$(a)$}"   ${\Cal F}$ is a continuous [semi]-psc trunk controller
\sn
\item "{$(b)$}"  $\Bbb Q$ is a straight ${\Cal F}$-forcing notion
\sn
\item "{$(c)$}"  $(\Bbb Q,\le_{pr})$ is $\aleph_1$-complete.
\ermn
\ub{Then} $\Bbb Q$ is a [semi] strong ${\Cal F}$-psc forcing notion.
\endproclaim
\bigskip

\proclaim{\stag{ct.2d} Claim}  Assume $\bar{\Bbb Q}$ is an explicitly
[Knaster/semi] ${\Cal F}$-psc iteration. \nl
1) If $\beta \le \ell g(\bar{\Bbb Q})$, then $\bar{\Bbb Q}
\restriction \beta$ is an explicitly [Knaster/semi] ${\Cal F}$-psc
iteration. \nl
2) If the trunk controller ${\Cal F}$ is fully based on $\langle {\Cal
F}^\beta:\beta < \alpha^* \rangle$ each $F^{[ 1+ \beta]}$ satisfies
[Knaster]-psc and ${\Cal F}^{[0]}$ satisfies [Knaster/semi]-psc,
\ub{then} ${\Cal F}$ satisfies [Knaster/semi]-psc. \nl
3) If $\Bbb Q$ is explicitly [semi] ${\Cal F}$-psc forcing, \ub{then}
$\Bbb Q$ is strongly [semi] ${\Cal F}$-psc forcing.
\endproclaim
\bigskip

\demo{Proof}  Should be clear.
\enddemo
\bigskip

\proclaim{\stag{ct.3} Claim}  1) If $\Bbb Q$ is ${\Cal F}$-psc and $p \in
\Bbb Q$ and ${\underset\tilde {}\to \tau_m}$ is a $\Bbb Q$-name of an
ordinal for $m < \omega$,
\ub{then} for some $q$ and $\langle \alpha_n:n < \omega \rangle$ we have:
\mr
\item "{$(a)$}"  $p \le_{\text{pr}} q$
\sn
\item "{$(b)$}"  $q \Vdash ``{\underset\tilde {}\to \tau_m} \in \{\alpha_n:
n < \omega\}$ for $m < \omega"$.
\ermn
2) If $\Bbb Q$ is strong ${\Cal F}$-psc, \ub{then} $\Bbb Q$ is
purely proper. \nl
3) If $\Bbb Q$ is explicitly [Knaster/semi] ${\Cal F}$-psc, \ub{then}
$\Bbb Q$ is purely proper.
\endproclaim
\bigskip

\remark{\stag{ct.3b} Remark}  1) If $\Bbb Q$ is straight (see
\scite{it.2a}(3)) we can add in \scite{ct.3}(1):
\mr
\item "{$(c)$}"  ${\Cal I}_q = \{r:q \le_{\text{apr}} r$ and $r$ forces a
value to $\underset\tilde {}\to \tau\}$ is predense over $p$; see \S3.
\ermn
2) Note that if every stationary ${\Cal S} \subseteq 
[2^{|{\Bbb Q}|}]^{\aleph_0}$ reflects in some $A \subseteq 2^{|{\Bbb
Q}|}$ of cardinality $\aleph_1$ then
also in \scite{ct.3}(1) we can get purely properness.  So the difference
is very small and still strange.
\endremark
\bigskip

\demo{Proof}  1) Assume not and let $\bold H$ be a witness for
``$\Bbb Q$ is psc".  So simulate a play of the game $\Game_p = 
\Game_{p,{\Bbb Q},\bold H}$, where
the interpolator plays using a fixed winning strategy whereas the extendor
chooses $q_\zeta$ such that:
\mr
\item "{$(\alpha)$}"  $p'_\zeta \le q_\zeta$ (see notation in
\scite{ct.1}(1)) \nl
(i.e. a legal move)
\sn
\item "{$(\beta)$}"  for some $m_\zeta < \omega$ the condition 
$q_\zeta$ forces a value to ${\underset\tilde {}\to \tau_{m_\zeta}}$,
call it $j_\zeta$
\sn
\item "{$(\gamma)$}"  $j_\zeta \notin \{j_\varepsilon:\varepsilon <
\zeta\}$.
\ermn
If the extendor can choose $q_\zeta$ for every $\zeta < \omega_1$, in the
end $\varepsilon < \zeta \and m_\varepsilon = m_\zeta 
\Rightarrow q_\varepsilon,q_\zeta$ are incompatible
(as $j_\varepsilon \ne j_\zeta$) but the interpolator has to win the play (as
he has used his winning strategy); contradiction by $(\beta) +
(\gamma)$ of \scite{ct.1}(1).

So necessarily for some $\zeta < \omega_1$ there is no $q$ as
required.  Let $p^*$ be $p'_\zeta$ so
$p \le_{\text{pr}} p^*$ by the definition of the game.  By our assumption
toward contradiction, for some $m$ we 
have $p^* \nVdash ``{\underset\tilde {}\to \tau_m} \in
\{j_\varepsilon:\varepsilon < \zeta\}"$ hence for some $q$ we have $p^* \le q$
and $q \Vdash ``{\underset\tilde {}\to \tau_m} \notin \{j_\varepsilon:
\varepsilon < \zeta\}"$ so \wilog \, for some $j$ we have
$q \Vdash ``{\underset\tilde {}\to \tau_m} = j"$.  But then the extendor could
have chosen $q_\zeta = q,j_\zeta=j$ and so clauses
$(\alpha),(\beta),(\gamma)$ holds, a contradiction. \nl

Note that we could have replaced clause $(\gamma)$ by
\mr
\item "{$(\gamma)^-$}"  $q_\zeta$ is incompatible with $q_\varepsilon$
whenever $\varepsilon < \zeta$ and $m_\xi = m_\varepsilon$. 
\ermn
2) Let $N \prec ({\Cal H}(\chi),\in,<^*_\chi)$ be countable, $\Bbb Q \in N$ and
$p \in N \cap {\Bbb Q}$.  Let $\langle {\underset\tilde {}\to \tau_n}:
n < \omega \rangle$ list the $\Bbb Q$-names of ordinals which belongs to $N$.
We define a strategy {\bf St}$_e$ for the extendor in the game 
$\Game_p$; so in stage $\zeta$ he has to choose $q_\zeta$ such 
that $p'_\zeta \le
q_\zeta$. If for some $n$, for no countable set $X$ of ordinals, $p^*_\zeta
\Vdash ``{\underset\tilde {}\to \tau_n} \in X"$ let $n(\zeta)$ be the minimal
such $n$, and as by the proof of the first part there are $(q,j)$ such that:
$p'_\zeta \le q$ and $q \Vdash ``{\underset\tilde {}\to \tau_{n(\zeta)}}
=j"$ and $j \notin \{i:\text{for some } \varepsilon < \zeta$ we have 
$q_\varepsilon \Vdash ``{\underset\tilde {}\to \tau_{n(\zeta)}} = i"\}$, choose
a $<^*_\chi$-minimal such pair, call it $(q_\zeta,j_\zeta)$, so the extendor
will choose $q_\zeta$.  If there is no such $n$, let $n(\zeta) = \omega$ and
$q_\zeta = p'_\zeta$.  Now in $\Game_p$ the interpolator has a winning
strategy {\bf St}$_i$, \wilog \, {\bf St}$_i \in N$.  Let $\langle p_\zeta,
q_\zeta:\zeta < \omega_1 \rangle$ be a play where the interpolator uses the
strategy {\bf St}$_i$ and the extendor uses the strategy {\bf St}$_e$, clearly
it exists and the interpolator wins.  
Clearly $\varepsilon < \zeta < \omega_1 \Rightarrow n(\varepsilon) \le
n(\zeta)$ (read the choices above).
Now we can prove by induction on $n$
that for some $\zeta < \omega_1,n(\zeta) > n$ and let $\zeta_n$ be the
minimal such $n$, so $p_{\zeta_n} \Vdash ``{\underset\tilde {}\to \tau_n}
\in X_n"$ for some countable $X_n$ set of ordinals. 
If we fail for $n$, then $\dbcu_{m < n} \zeta_m < \zeta < \varepsilon <
\omega_1 \Rightarrow (q_\zeta \Vdash \underset\tilde {}\to \tau = 
j_\zeta) \and (q_\xi \Vdash \underset\tilde {}\to \tau = j_\varepsilon) 
\Rightarrow (q_\zeta,q_\varepsilon$ incompatible), 
(see choice of $j_\varepsilon$), but
this contradicts the use of {\bf St}$_i$.  Now $\langle p_\zeta,q_\zeta:
\zeta < \zeta_n \rangle$ can be defined from $p$, {\bf St}$_i,
\langle {\underset\tilde {}\to \tau_\ell}:\ell \le n \rangle$ and
$\bold H$ (read the
definition of {\bf St}$_e$) hence $\langle p_\zeta,q_\zeta:\zeta < \zeta_n
\rangle \in N$, so $\zeta_n \in N$, and similarly $p_{\zeta_n} \in N$.
So as $p_{\zeta_n} \Vdash ``{\underset\tilde {}\to \tau_n} \in X_n"$, the set
$\{\xi:p_{\zeta_n} \nVdash ``{\underset\tilde {}\to \tau_n} \ne \xi"\}$ is
countable and it belongs to $N$.  So $p_{\zeta_n} \Vdash
``{\underset\tilde {}\to \tau_n} \in N \cap \text{ Ord}"$.  Now if $\zeta <
\omega_1$ is $\ge \dbcu_{n < \omega} \zeta_n$ then $p_\zeta$ is as
required. \nl
3) Same proofs, only easier (in some cases follows by \scite{ct.2d}(3)).
  \hfill$\square_{\scite{ct.3}}$\margincite{ct.3}.
\enddemo
\bigskip

\remark{\stag{ct.3c} Remark}  Of course, if $(Q^{\Bbb Q},\le^{\Bbb
Q}_{pr})$ is $\aleph_1$-complete, part (2) of \scite{ct.3} follows
from part (1) easily.
\endremark
\bigskip

\proclaim{\stag{ct.4} Lemma}  Assume that $\bar{\Bbb Q}$ is a [strong], [semi]
${\Cal F}$-psc-iteration.
\ub{Then} for every $\beta \le \ell g
(\bar{\Bbb Q})$ we have ${\Bbb P}_\beta$ is a [strong], [semi] ${\Cal
F}$-psc forcing notion.
\endproclaim
\bigskip

\demo{Proof}   Let $\bar{\bold H} = \langle \bold H_\beta:\beta < \ell
g(\bar{\Bbb Q})$ be a witness for ``$\bar{\Bbb Q}$ is a [strong], [semi]
${\Cal F}$-psc iteration".
Let $\beta \le \ell g(\bar{\Bbb Q})$ and $p \in \Bbb P_\beta$, and define
$\bold H^\beta$ naturally composing the $\langle \bold H_\gamma:\gamma <
\beta \rangle$ and we shall describe a winning strategy for the 
interpolator in the game $\Game'_p = \Game'_{p,{\Bbb P}_\beta,{\bold H}^\beta}$;
he just guarantees that:
\mr
\item "{$(a)$}"  Dom$(p_\zeta) = \text{ Dom}(q_\zeta)$
\sn
\item "{$(b)$}"  if $p_\zeta(\gamma) \ne q_\zeta(\gamma)$ then $\gamma \in
\dbcu_{\varepsilon < \zeta} \text{ Dom}(p_\varepsilon)$
\sn
\item "{$(c)$}"  if $\gamma \in \dbcu_{\zeta < \omega_1} \text{ Dom}(p_\zeta)$
and $\xi(\gamma) = \text{ Min}\{\zeta:\zeta \in \text{ Dom}(p_\zeta)\}$
\ub{then} $\Vdash_{{\Bbb P}_\gamma} ``\langle p_{\xi(\gamma)+1+ \zeta}
(\gamma),q_{\xi(\gamma)+1+\zeta}(\gamma):\zeta < \omega_1 \rangle$ 
is a play of $\Game'_{p_{\xi(\gamma)}(\gamma)},
[{\underset\tilde {}\to {\Bbb Q}_\gamma},\bold H^\gamma]$ or the semi version in which 
the interpolator uses a fixed winning strategy 
${\underset\tilde {}\to {\bold St}^\gamma_{p_\xi(\gamma)}}$"
\sn
\item "{$(d)$}"  Dom$(p'_\zeta) =
\cup\{\text{Dom}(p_\varepsilon):\varepsilon < \zeta\}$.
\ermn
Let $\langle p'_\zeta,p_\zeta,q_\zeta:\zeta < \omega_1 \rangle$ be a play of the game
$\Game_p[{\Bbb P}_\beta]$ in which the interpolator uses the strategy
described above.  To prove that the interpolator wins the play, let $A
\subseteq \omega_1$ be stationary and $A \in \bold V$, of course.  
For $\zeta \in \dbcu_{\zeta < \omega_1} \text{ Dom}(p_\zeta)$ let
$E_\zeta$ be as in clause $(\gamma)$ of \scite{ct.1}(1) (or
the relevant variant) and $E = \{\delta:\delta \text{ limit and }
\zeta \in \dbcu_{\varepsilon < \delta} \text{ Dom}(p_\varepsilon)
\Rightarrow \delta \in E_\zeta\}$.
Define $w_\zeta = \{\gamma \in
\text{ Dom}(p_\zeta):q_\zeta \restriction \gamma
\nVdash_{{\Bbb P}_\gamma} ``p_\zeta(\gamma) 
\le_{\text{pr}} q_\zeta(\gamma)"\}$, so by the definition of 
${\Cal F}$-iteration, $w_\zeta$ is finite, and by the strategy of the 
interpolator we have $w_\zeta \subseteq \dbcu_{\varepsilon < \zeta}
\text{ Dom}(p_\varepsilon)$.  So by Fodor lemma for some stationary $A_0
\subseteq A \cap E$ we have $\zeta \in A_0 \Rightarrow w_\zeta = w^*$.

Letting $w^* = \{\gamma_\ell:\ell < k\}$ such that 
$\gamma_\ell < \gamma_{\ell+1}$, we choose by induction on 
$\ell \le k$ a stationary set $A_\ell \subseteq \omega_1$ (from $\bold V$, of
course) such that $A_{\ell +1} \subseteq A_\ell$ and 
$\bold H_{\gamma_\ell}(\langle 
\text{val}^{\underset\tilde {}\to Q_{\gamma_\ell}} (q_\zeta(\gamma_\ell)):-1
\le \zeta < \omega_1 \rangle,A_{\ell +1})$ is 1.  
For $\ell = 0,A_0$ has already been chosen, 
for $\ell +1$, in $\bold V^{P_{\gamma_\ell}}$, we know that
$\langle(p_{\xi(\gamma_\ell)+1+ \zeta}(\gamma_\ell),
q_{\xi(\gamma_\ell)+1+ \zeta}(\gamma_\ell)):\zeta < \omega_1 \rangle$ is a
play of the game $\Game_{p_\xi(\gamma)}[\Bbb Q_{\gamma_\ell},
\bold H_{\gamma_\ell}]$ in which the interpolator
uses a winning strategy and
$\langle(\text{val}^{{\Bbb Q}_{\gamma_\ell}}(p_{\xi(\gamma_\ell) + 1 +
\zeta}(\gamma_\ell)$, val$^{{\Bbb
Q}_{\gamma_\ell}}(q_{\xi(\gamma_\ell)+1+ \zeta})):\zeta < \omega_1
\rangle$ is a sequence of pairs from ${\Cal F}^{[\gamma_\ell]}$, and the
sequence is from $\bold V$.  So there is a stationary $A_{\ell +1} \subseteq
A_\ell$ as required above.
Lastly, let $B =: A_{k^*}$, we shall prove that $B$ is as required,
that is, it satisfies clause $(\beta)$ of \scite{ct.1}(1), 
we concentrate on $(\beta)$; the other clause, $(\gamma)$ is similar and not
really needed, suffice to get the finished variant.
Let $\varepsilon < \zeta$ be from $B$ and we shall find $r$ as required.
Stipulate $\gamma_k = \alpha^*$ and we now let $Y = \text{ Dom}(q_\zeta)
\cup \{\alpha^*\}$ and we choose by induction on $\gamma \in Y$ a condition
$r_\gamma \in \Bbb P_\gamma$ such that (we ignore the ``semi" version,
which requires just special treatment for $\gamma = 0$)
\mr
\widestnumber\item{$(iii)$}
\item "{$(i)$}"  $q_\varepsilon \restriction \gamma \le r_\gamma$ and
\sn
\item "{$(ii)$}"  $q_\zeta \restriction \gamma \le r_\gamma$ and
\sn
\item "{$(iii)$}"  $\gamma = \beta + 1 \and \beta \in \text{ Dom}(r_\gamma) \backslash w^* 
\Rightarrow r_\gamma \Vdash_{{\Bbb P}_\gamma} ``q_\varepsilon(\beta)
\le_{\text{pr}} r_\gamma(\beta)$ and $q_\zeta(\beta)
\le_{\text{pr}} r_\gamma(\beta)"$
\sn
\item "{$(iv)$}"  if $\beta \in Y \cap \gamma$ and $[\beta,\gamma]
\cap w^* = \emptyset$ then $r_\beta = r_\gamma
\restriction \beta$
\endroster
\enddemo
\bn
\ub{Case 1}:  $\gamma = \text{ Min}(Y)$.

Let $r_\gamma = \emptyset$, the empty function.
\bn
\ub{Case 2}:  $\gamma \ne \text{ Min}(Y)$ but $\gamma \cap Y$ has no last
element.

Let $r_\gamma = \cup\{r_\beta:\beta \in \gamma \cap Y\}$, now check; 
it is a well defined function by clause (iv).
We use here ``${\Cal F}$ is fully based on $\alpha$", see Definition
\scite{it.1}(5) to show $r_\gamma \in \Bbb P_\beta$.  Lastly
note that for checking $q_\zeta \restriction \gamma \le r_\gamma$ and
$q_\varepsilon \restriction \gamma \le r_\gamma$ we are using
clause (iii).
\bn
\ub{Case 3}:  $Y \cap \gamma$ has last element $\beta,\beta \notin w^*$ and
$\beta \notin \text{ Dom}(q_\varepsilon)$.

We let $r_\gamma(\beta) = q_\zeta(\gamma)$.
\bn
\ub{Case 4}:  $Y \cap \gamma$ has last element $\beta,\beta \notin w^*$ but
$\beta \in \text{ Dom}(q_\varepsilon)$.

Now, $r_\beta$ necessarily forces that
\mr
\item "{$(*)$}"  $p_\varepsilon(\beta) \le_{\text{pr}} p_\zeta(\beta)
\le_{\text{pr}} q_\zeta(\beta)$ \nl
and $p_\varepsilon(\beta) \le_{\text{pr}} q_\varepsilon(\beta)$.
\ermn
So by clause $(\gamma)$ of Definition \scite{ct.1}(1), $r_\beta$ forces that
in ${\underset\tilde {}\to {\Bbb Q}_\beta}$ there is 
$r$ such that $q_\zeta(\beta) \le_{\text{pr}} r,q_\varepsilon(\beta) \le_{\text{pr}} r$, and
val$^{\underset\tilde {}\to {\Bbb Q}_\beta}(r)$ is as it should be by
clause $(\delta)$ of \scite{ct.1}(1) (so is an object, not just a
$\Bbb P_\beta$-name).  Lastly let $r_\gamma
(\beta)$ be a $\Bbb P_\beta$-name of such a condition.
\bn
\ub{Case 5}:  $Y \cap \gamma$ has last element $\beta,\beta = \gamma_\ell
\in w^*$.

By the choice of $A_{\ell +1}$, we know $\Vdash_{{\Bbb P}_\beta}
``p_\zeta(\beta),q_\zeta(\beta)$ are $\le^{\underset\tilde {}\to {\Bbb
Q}_\beta}$ compatible".  So, for some $r'_\beta,r_\beta \le^{{\Bbb
P}_\beta} r'$, for some $x \in {\Cal F}^{[\beta]} r'$ forces that some
$\le^{\underset\tilde {}\to {\Bbb Q}_\beta}$-common upper bound
$\underset\tilde {}\to r"$ of $p_\zeta(\beta),q_\zeta(\beta)$,
satisfies
val$^{\underset\tilde {}\to {\Bbb Q}_\beta}({\underset\tilde {}\to r''})
= x$, and let $r_\gamma = r'_\beta \cup \{(\beta,{\underset\tilde
{}\to r''})\}$.

So we are done.  \hfill$\square_{\scite{ct.4}}$\margincite{ct.4}
\bigskip

\proclaim{\stag{ct.4a} Claim}  Assume $\bar{\Bbb Q}$ is a [Knastner]
explicit [semi] ${\Cal F}$-psc iteration.  \ub{Then} Lim$_{\Cal
F}(\bar{\Bbb Q})$ satisfies the [Knastner/semi] explicit ${\Cal F}$-psc.
\endproclaim
\bigskip

\demo{Proof}  Similar to \scite{ct.4}.
\enddemo
\bn
\ub{\stag{ct.5} Discussion}:   We may consider replacing stationary 
$A,B \subseteq \omega_1$ by a subset of $[\omega_1]^2$, so we use:
\definition{\stag{ct.6} Definition}  We call $({\frak H},<_{\frak H})$ is c.c.c. witness
if:
\mr
\item "{$(a)$}"  $\le_{\frak H}$ is a partial order of ${\frak H}$
\sn
\item "{$(b)$}"  ${\frak H} \subseteq {\Cal P}([\omega_1]^2) \backslash
\{\emptyset\}$
\sn
\item "{$(c)$}"  for any $X \in {\frak H}$, club $E$ of $\omega_1$ 
and a pressing
down function $h$ on $E$ for some $Y \le_{\frak H} X$ we have $h \restriction
\cup \bigl\{ \{\alpha,\beta\}:\{\alpha,\beta\} \in {\frak Y} \bigr\}$ is
constant.
\endroster
\enddefinition
\newpage

\head {\S3 Nicer pure properness and pure decidability} \endhead  \resetall \sectno=3
\bigskip

\proclaim{\stag{mr.1} Claim}   Let ${\Cal F}$ be a trunk controller
over $\alpha^{**}$ and $\bar{\Bbb Q}$ be a semi ${\Cal F}$-psc iteration
and is semi-straight (see \scite{ct.1}(6)) and $\alpha^* = \ell g(\bar{\Bbb Q})$. \nl
1) If $p \Vdash_{{\Bbb P}_{\alpha^*}} ``\underset\tilde {}\to \tau:\omega
\rightarrow \text{ Ord}"$, \ub{then} we can find $q$ and 
${\Cal I}_n(n < \omega)$ such that: $p \le_{\text{pr}} q \in 
\Bbb P_{\alpha^*}$ and for each $n < \omega$ we have:
\medskip
\noindent
$(*)_{q,{\Cal I}_n,\underset\tilde {}\to \tau(n)}$
\mr
\item "{$(a)$}"  $q \in \Bbb P_{\alpha^*}$
\sn
\item "{$(b)$}"  ${\Cal I}_n \subseteq \{r:q \le_{\text{apr}} r
\text{ and } r \text{ forces a value to } \underset\tilde {}\to \tau(n)\}$
\sn
\item "{$(c)$}"  ${\Cal I}_n$ is countable
\sn
\item "{$(d)$}"  ${\Cal I}_n$ is predense over $q$.
\ermn
2) If above $\bar{\Bbb Q}$ is a strong semi ${\Cal F}$-psc iteration and
$\{p,\bar{\Bbb Q}\} \in N \prec ({\Cal H}(\chi),\in)$ and $N$ is
countable and $\langle \underset\tilde {}\to \tau(n):n < \omega
\rangle$ list the $\Bbb P_{\alpha^*}$-names of ordinals from $N$,
\ub{then} we can add
\mr
\item "{$(e)$}"  $q$ is $(N,\Bbb P_{\alpha^*})$-generic. 
\ermn
3) Assume $\Bbb Q$ is an ${\Cal F}$-iteration and:
\mr
\item "{$(a)$}"  $\beta^* \le \alpha^*$ and $\beta \in [\beta^*,\alpha^*)$ 
implies $\Vdash_{{\Bbb P}_\beta} ``{\underset\tilde {}\to {\Bbb Q}_\beta}$ has
pure $(2,2)$-decidability", see Definition \scite{it.11}
\sn
\item "{$(b)$}" ${\underset\tilde {}\to {\bold t}}$ is a
$\Bbb P_{\alpha^*}$-name, $\Vdash_{{\Bbb P}_{\alpha^*}} ``{\underset\tilde {}\to {\bold t}} \in \{0,1\}"$
\sn
\item "{$(c)$}"  $(*)_{p,{\Cal I},{\underset\tilde {}\to {\bold t}}}$ holds
\sn
\item "{$(d)$}"  each ${\underset\tilde {}\to {\Bbb Q}_\beta}$ (for $\beta \in 
[\beta^*,\alpha^*))$ satisfies $p' \le_{\text{pr}} p'' \Rightarrow
\text{ val}^{\underset\tilde {}\to {\Bbb Q}_\beta}(p') = 
\text{ val}^{\underset\tilde {}\to {\Bbb Q}_\beta}(p'')$.
\ermn
\ub{Then} there are ${\underset\tilde {}\to {\bold t}'},q'$ such that:

$$
p \le_{\text{pr}} q'
$$

$$
q' \restriction \beta^* = p \restriction \beta^*
$$

$$
\bold t' \text{ is a } \Bbb P_{\beta^*}\text{-name}
$$

$$
q' \Vdash_{{\Bbb P}_{\alpha^*}} ``{\underset\tilde {}\to {\bold t}'} =
{\underset\tilde {}\to {\bold t}''}.
$$
\mn
4) Similarly for pure $(\aleph_0,2)$-decidability.
\endproclaim
\bigskip

\demo{Proof}  1) Assume not and let $\bar{\bold H}$ be a witness for
``$\bar{\Bbb Q}$ is ${\Cal F}$-psc" hence by \scite{ct.4} some $\bold
H$ witness $P_{\alpha^*}$ is ${\Cal F}$-pcf.  So simulate 
a play of the game $\Game_p = \Game_{p,{\Bbb P}_{\alpha^*},\bold H}$, where
the interpolator plays using a fixed winning strategy whereas the extendor
chooses $q_\zeta$ and $n_\zeta \le \omega$ such that:
\mr
\item "{$(\alpha)$}"  $p'_\zeta \le q_\zeta$ (see notation in
\scite{ct.1}(1)) \nl
(i.e. a legal move)
\sn
\item "{$(\beta)$}"  $n_\zeta \le \omega$ is minimal $n$ such that
$\{q_\varepsilon:\varepsilon < \zeta,n_\varepsilon = n$, and
$q_\varepsilon$ forces a value of $\underset\tilde {}\to \tau(n)\}$ is
not predense over $p^*_\zeta$
\sn
\item "{$(\gamma)$}"  if $n_\zeta < \omega,q_\zeta$ 
forces a value to $\underset\tilde {}\to \tau(n_\zeta)$,
call it $j_\zeta$
\sn
\item "{$(\delta)$}"  if $n_\zeta < \omega$ then $q_\zeta$ is
incompatible with $q_\varepsilon$ if $\varepsilon < \zeta \and
n_\varepsilon = n_\zeta$.
\ermn
Now
\mr
\item "{$\boxtimes$}"  for some $\zeta,n_\zeta = \omega$. \nl
Why?  Otherwise the extendor can choose $q_\zeta$ for every $\zeta <
\omega_1$, and $n_\zeta = n^*$ for every $\zeta \in [\zeta^*,\omega_1)$
for some $\zeta^*$; in the end 
$\zeta^* < \varepsilon < \zeta \Rightarrow q_\varepsilon,q_\zeta$ are incompatible
(as $j_\varepsilon \ne j_\zeta$) but the interpolator has to win the play (as
he has used his winning strategy), contradicting clauses $(\beta) +
(\gamma)$ of \scite{ct.1}(1).
\ermn
So necessarily for some $\zeta < \omega_1,n_\zeta = \omega$.  
Let $p^*$ be: $p$ if $\zeta=0,p_{\zeta -1}$ if $\zeta$ is a 
successor ordinal and $p'_\zeta$ if $\zeta$ is a limit ordinal, so
$p \le_{\text{pr}} p^*$ by the definition of the game.  
For each $\varepsilon < \zeta$ let $q'_\varepsilon$ be a $\le$-lub of
$q_\varepsilon,p^*$, exists as $\Bbb Q$ is straight, so $p^* \le_{pr}
q'_\varepsilon$.  Let ${\Cal I}_n = \{q'_\varepsilon:\varepsilon <
\zeta^*$ and $n_\varepsilon=n\}$ and we shall show that $p^*,{\Cal
I}_n$ are as required ($p^*$ standing for $q$).  Now clauses (a), (b),
(c) are obvious, toward clause (d) assume $n < \omega,p^* \le q$ and
$q$ is incompatible with all members of ${\Cal I}_n$ and let $\zeta_n
= \text{ Min}\{\zeta:n_\zeta > n\}$, so $q$ could not have been a good
candidate for $q_{\zeta_n}$ hence is compatible with some
$q_\varepsilon,n_\varepsilon = n$.  So by the choice of
$q'_\varepsilon$ clearly $q'_\varepsilon \le q$ and $q'_\varepsilon
\in {\Cal I}_n$, contradiction. \nl
2) Essentially the same proof.
\nl
3) \ub{Case 1}:  Dom$(p) \subseteq \beta^*$.

Trivial.
\mn
\ub{Case 2}:  Dom$(p) \backslash \beta^* \ne \emptyset$.

For each $ \beta \in [\beta^*,\alpha^*)$ we define $\Bbb P_\beta$-names
${\underset\tilde {}\to {\bold t}^0_\beta},
{\underset\tilde {}\to {\bold t}^1_\beta},
{\underset\tilde {}\to {\bold t}^2_\beta},
{\underset\tilde {}\to q_\beta}$ as follows.  
Let Dom$^+(p) = \text{ Dom}(p) \cup \{\cup\{\gamma +1:\gamma \in
\text{ Dom}(p)\}\}$.  Let $G_\beta \subseteq
\Bbb P_\beta$ be generic over $\bold V$.
\mn
\ub{Possibility A}:  There are $q \in \Bbb P_{\alpha^*}/G_\beta$ 
such that $p \le_{\text{pr}} q$, Dom$(q) = \text{ Dom}(p)$ and
$q \Vdash_{{\Bbb P}_{\alpha^*}/G_\beta} ``{\underset\tilde {}\to {\bold t}} = \ell"$.

Let ${\underset\tilde {}\to {\bold t}^1_\beta}[G_\beta]$ be 1
and $\bold t^0_\beta[G_\beta] = \ell$ and let
${\underset\tilde {}\to q_\beta}[G_\beta]$ be $(p \restriction \beta) \cup
(q \restriction [\beta,\alpha^*]$.
\mn
\ub{Possibility B}:  Not possibility A.

Let ${\underset\tilde {}\to {\bold t}^0_\beta}[G_{\beta^*}] = 0,
{\underset\tilde {}\to {\bold t}^1_\beta}[G_\beta]=0$ and
${\underset\tilde {}\to q_\beta}[G_\beta] = p$.

In the first possibility, note that we have demanded Dom$(q) = \text{ Dom}(p)$; in the
second posibility this holds automatically.

Let ${\underset\tilde {}\to {\bold t}^2_\beta}[G_\beta] \in \{0,1\}$ be 1 
iff
${\underset\tilde {}\to {\bold t}^1_\beta}[G_\beta] = 0$ and for no
$\gamma \in \beta \cap \text{ Dom}(p) \backslash \beta^*$ do we have
${\underset\tilde {}\to {\bold t}^1_\gamma}[G_\beta \cap \Bbb
P_\gamma]=1$, clearly also ${\underset\tilde {}\to {\bold t}^2_\beta}$
is a $\Bbb P_\beta$-name of a number $\in \{0,1\}$.
\sn
Now note 
\mr
\item "{$\boxtimes$}"  $\Vdash_{P_{\alpha^*}}$ ``there is one and only one
$\beta \in \text{ Dom}(p) \backslash \beta^*$ such that $\bold t^2_\beta
[G_{\alpha^*} \cap \Bbb P_\beta] = 1$  \nl

$\qquad \quad$ call it $\underset\tilde {}\to \beta$"
\ermn
(the ``at most one" follows by the definition of ${\underset\tilde
{}\to {\bold t}^2_\beta}$, and the ``at least one" is by the 
finiteness clause in the definition of
$\le^{{\Bbb P}_{\alpha^*}}$).
So we can define $q \in \Bbb P_{\alpha^*}:\text{Dom}(p) = \text{ Dom}(q)$,
for $\gamma \in \text{ Dom}(p)$: if $\gamma \ge \beta^*$ and $\gamma \ge
\underset\tilde {}\to \beta$ then $q(\gamma) =
{\underset\tilde {}\to q_{\underset\tilde {}\to \beta}}(\gamma)$,
otherwise $g(\gamma) = p(\gamma)$.

Now for each $\beta \in \text{ Dom}(p) \backslash \beta^*$ we define a
$\Bbb P_\beta$-name $q^{[\beta]} \in {\underset\tilde {}\to Q_\beta}$,
such that:
\mr
\item "{$(\alpha)$}"  $q(\beta) \le_{\text{pr}} q^{[\beta]}$ and
\sn
\item "{$(\beta)$}"   val$^{\Bbb Q}(q^{[\beta]}) =
\text{ val}^{\Bbb Q}(p(\beta))$ [note this does not work well with
associativity of iterations] and 
\sn
\item "{$(\gamma)$}"  if $\ell \le 2$, and there are 
$q',\underset\tilde {}\to i$ such that $\Vdash_{{\Bbb P}_\beta}
``q \le_{pr} q'$ and $q' \Vdash_{\underset\tilde {}\to {\Bbb Q}_\beta}
{\underset\tilde {}\to {\bold t}^\ell_{\beta +1}} = \underset\tilde
{}\to i"$ where $\underset\tilde {}\to i$ is a $\Bbb P_\beta$-name,
\ub{then}
$(q^{[\beta]},{\underset\tilde {}\to i_\ell})$ satisfy this for some
$\Bbb P_\beta$-name ${\underset\tilde {}\to i_\beta}$ of a number $< 2$.
\ermn
Let $q^* \in \Bbb P _\beta$, Dom$(q^*) = \text{ Dom}(q),q^* \restriction
\beta^* = q \restriction \beta^*$ and $\beta \in \text{ Dom}(q) \backslash
\beta^* \Rightarrow q^*(\beta) = q^{[\beta]}$. 
We can find $G_{\alpha^*} \subseteq
\Bbb P_{\alpha^*}$ generic over $\bold V$ satisfying $q^* \in
G_{\alpha^*}$ 
and ${\underset\tilde {}\to {\bold t}}[G^\ell_\alpha] = \ell$ for some
$\ell \in \{0,1\}$.  Now:
\mr
\item "{$(*)_1$}"  for some $\beta \in [\beta^*,\alpha^*] 
\cap \text{ Dom}^+(p)$ we have 
${\underset\tilde {}\to {\bold t}^1_\beta}[G^\ell_{\alpha^*} \cap
\Bbb P_\beta] = 1$ \nl
[why?  as $p \le_{\text{pr}} q^* \in G_{\alpha^*}$, by assumption (c) necessarily 
${\Cal I} \cap G^\ell_{\alpha^*} \ne \emptyset$ so let $r \in {\Cal I} \cap
G^\ell_{\alpha^*}$, and so $\beta = \alpha^*$ is as required.]
\ermn
So let $\beta \in [\beta^*,\alpha^*]$ be minimal such that
${\underset\tilde {}\to {\bold t}^1_\beta}[G_{\alpha^*} \cap
\Bbb P_\beta]=1$
\mr
\item "{$(*)_2$}"  $\beta$ cannot be a limit ordinal $> \beta^*$
\nl
[why? by the finiteness clause in the definition of order in 
$\Bbb P_{\alpha^*}$]
\sn
\item "{$(*)_3$}"  $\beta = \gamma + 1 > \beta^*$ is impossible.
\nl
[Why?  If $\gamma \notin \text{ Dom}(p)$ this is trivial by the choice
of $\beta$, so assume $\gamma \in \text{ Dom}(p)$.  Now in $\bold
V[G_{\alpha^*} \cap \Bbb P_\gamma]$ the forcing notion
${\underset\tilde {}\to {\Bbb Q}_\gamma}[G_{\alpha^*} \cap \Bbb
P_\gamma]$ has pure $(2,2)$-decidability hence clearly
$q^{[\gamma]}[G_{\alpha^*} \cap \Bbb P_\gamma] \Vdash_{{\Bbb
Q}_\gamma}
``{\underset\tilde {}\to {\bold t}^\ell_{\gamma +1}} =
{\underset\tilde {}\to i^\ell_\gamma}[G_{\alpha^*} \cap \Bbb
P_\gamma]"$.  Now ${\underset\tilde {}\to {\bold t}^1_{\gamma
+1}}[G_{\alpha^*} \cap \Bbb P_{\gamma +1}] = 1$ by the choice of
$\beta = \gamma +1$, hence ${\underset\tilde {}\to
i^1_\gamma}[G_{\alpha^*} \cap \Bbb P_\gamma] = 1$.  Define $q':q'
\restriction \beta = q^* \restriction \beta,p' \restriction
(\text{Dom}(p)) \backslash \beta) = {\underset\tilde {}\to q_\beta}
\restriction (\text{Dom}(p) \backslash \beta)$ it proves $\gamma$
could serve instead $\beta$, contradiction.] 
\ermn
So $\beta_\ell = \beta^*$ and we are easily done by the definition of
$\underset\tilde {}\to p$. \nl
4) Similar to the proof of part (3).  \hfill$\square_{\scite{mr.1}}$\margincite{mr.1}
\enddemo
\bigskip

\definition{\stag{mr.2} Definition}  1) Let $S \subseteq
[\lambda]^{\aleph_0}$ be stationary.  We say that $\Bbb Q$, is purely
$(S,{\Cal F})$-proper \ub{if} $N \prec ({\Cal H}(\chi),\in)$ is
countable, $N \cap \lambda \in S,\Bbb Q,{\Cal F} \in N,p \in \Bbb Q \cap
N$ \ub{then} there is $q$ such that:
\mr
\item "{$(a)$}"  $p \le_{pr} q \in \Bbb Q$
\sn
\item "{$(b)$}"  $q$ is $(N,\Bbb Q)$-generic
\sn
\item "{$(c)$}"  val$^{\Bbb Q}(q) = \text{ prop}_{\Cal F}(p,N \cap S)$,
as usual prop$_{\Cal F}$ is considered part of ${\Cal F}$.
\endroster
\enddefinition
\bigskip

\proclaim{\stag{mr.3} Claim}  In \scite{mr.1}(2) we can conclude $\Bbb Q$
is purely $(([\lambda]^{\aleph_0}),{\Cal F})$-proper.
\endproclaim
\bigskip

\demo{Proof}  Read the claims and definitions.
\enddemo
\bn
We can also generalize the preservation theorems.
\proclaim{\stag{mr.6} Claim}  1) Assume
\mr
\item "{$(a)$}"  $\bar{\Bbb Q},{\Cal F}$ are as in \scite{mr.2}
\sn
\item "{$(b)$}"  $(D,R,<)$ is a fine covering \footnote{can use also
the weaker version there}  model in the sense of 
\cite[VI,Definition 1.x]{Sh:f}
\sn
\item "{$(c)$}"   $E \Vdash_{{\Bbb P}_\beta} ``\Bbb Q_\beta$ is purely
$(D_1,R,<)$-preserving
\sn
\item "{$(d)$}"  each ${\underset\tilde {}\to {\Bbb Q}_\beta}$ is
purely preserving.
\ermn
2) Similary for $\bar{\Bbb Q}$ a Knaster explictly semi ${\Cal F}$-psc
iteration.
\endproclaim
\bigskip

\demo{Proof}  See \cite[VI,1.13A,p.270]{Sh:f}.
\enddemo
\newpage

\head {\S4 \\
Averages by an ultrafilter and restricted non null trees} \endhead  \resetall 
\bigskip

\proclaim{\stag{pr.1} Claim}  Assume
\mr
\item "{$(a)$}"  $\bold V_1 = \bold V[\bold r]$ where $\bold r$ is a random
real over $\bold V$
\sn
\item "{$(b)$}"  in $\bold V,D$ is a non principal ultrafilter on $\omega$.
\ermn
\ub{Then} we can find in $\bold V_1$ a non principal ultrafilter $D_1$ on
$\omega$ extending $D$ such that
\mr
\item "{$(*)_{\bold r,D,D_1}$}"  Assume in 
$\bold V,({\Cal B},\langle {\Cal B}_n:
n < \omega \rangle)$ is a candidate (or a $(D,\bold V)$-candidate which 
means that (all this in $\bold V$),
${\Cal B},{\Cal B}_n$ are Borel subsets of ${}^\omega 2$ and
$(\forall k < \omega)(\exists B \in D)(\forall n \in B)[\text{Leb}
({\Cal B} \Delta {\Cal B}_n) < 1/(k+1)]$.
\nl
\ub{Then} the following conditions are equivalent in $\bold V_1$:
{\roster
\itemitem{ $(i)$ }  $\bold r \in \dbcu_{n \in B} {\Cal B}_n$ 
for every $B \in
D$ (recall, ${\Cal B}_n$ being a Borel set, is actually a definition of a
set and so $\dbcu_{n \in B} {\Cal B}_n$ is a definition in $\bold V$ of a
Borel set, so it defines a Borel set in $\bold V[\bold r]$)
\sn
\itemitem{ $(ii)$ }  the set $\{n:\bold r \in {\Cal B}_n\}$ belongs to $D_1$.
\endroster}
\endroster
\endproclaim
\bigskip

\demo{Proof}  Clearly $(i) \Rightarrow (ii)$.  So it suffices to find
in $\bold V_1$ an ultrafilter $D_1$ over $\omega$ containing $D' = D \cup
\{\{n:\bold r \in \bold{\Cal B}_n\}:({\Cal B},\langle {\Cal B}_n:n < \omega
\rangle)$ is a candidate$\}$.  This is equivalent to ``any
intersection of finitely many members of $D'$ is not empty.  As $D$ is
closed under finite intersection. 
Clearly it suffices to prove:
\mr
\item "{$\boxtimes$}"  Assume in $\bold V$ that $m^* < \omega$ and for each
$m < m^*,({\Cal B}_m,\langle {\Cal B}_{m,n}:n < \omega \rangle$ 
is a candidate and $B \in D$. \nl
\ub{Then} for some $n \in B$ in $\bold V_1$ we have $m < m^* \Rightarrow
[\bold r \in {\Cal B}_m \equiv \bold r \in {\Cal B}_{m,n}]$.
\endroster
\enddemo
\bigskip

\demo{Proof of $\boxtimes$}  It is enough, given a positive real $\varepsilon
> 0$ to find a Borel set ${\Cal B} = {\Cal B}_\varepsilon \in \bold V$ 
of Lebesgue measure $< \varepsilon$ such that
\mr
\item "{$(*)$}"  $r \in {}^\omega 2 \backslash 
{\Cal B} \Rightarrow (\exists n \in B)(\forall m < m^*)
[r \in {\Cal B}_m \equiv r \in {\Cal B}_{m,n}]$
\ermn
(as then we can find in $\bold V$ a sequence $\langle {\Cal
B}_{1/(k+1)}:k < \omega \rangle$, each ${\Cal B}_{1/(k+1)}$ as above;
so $\bold r$ being random over $\bold V$ does not belong to $\dbca_k
{\Cal B}_{1/(k+1)}$ hence for some $k,\bold r \notin {\Cal
B}_{(1/(k+1)}$ and so there is $n$ as required by $(*)$ because $(*)$
holds also in $\bold V[\bold r]$ by absoluteness).

Given $\varepsilon$, for $m < m^*$, as $({\Cal B}_m,\langle {\Cal
B}_{m,n}:n < \omega \rangle)$ is a candidate  
we can find $B'_m \in D$ such that $(\forall n \in B'_m)
([\text{Leb}({\Cal B}_m \Delta {\Cal B}_{m,n}) < \varepsilon/(m^* +1)]$.
Let $B'' = \dbca_{m < m^*} B'_m \cap B$, so clearly $B'' \in D$; now
for any $n \in B''$, the 
Borel set ${\Cal B} = \dbcu_{m < m^*}({\Cal B} \Delta {\Cal B}_{m,n})$ will do:
clearly Leb$({\Cal B}) < \varepsilon$ and easily $n$ is as required.
\hfill$\square_{\scite{pr.1}}$\margincite{pr.1}
\enddemo
\bn
We can do the following more generally, but the following was our required
example.

\definition{\stag{pr.2} Definition}  If $g:\omega \rightarrow \omega$ satisfies
$n < \omega \Rightarrow g(n) > n$ and is increasing 
we define $\bold T_g$ as the family of subtrees 
$T$ of ${}^{\omega >}2$ such that for
every $n < \omega$ and $\eta \in T \cap {}^n 2$ we have
\footnote{recall
$T^{[\eta]} = \{\nu \in T:\nu \trianglelefteq \eta \text{ or } \eta
\triangleleft \nu\}$}

$$
(1-1/n)|T^{[\eta]} \cap {}^{g(n)}2|/2^{g(n)} \le \text{ Leb}(\text{lim}
(T^{[\eta]})) \le |T^{[\eta]} \cap {}^{g(n)}2|/2^{g(n)} \tag "{$(*)_n$}"
$$
\mn
(the second inequality holds automatically), equivalently, for
$m \ge g(n)$

$$
(1-1/n)|T^{[\eta]} \cap {}^{g(n)}2|/2^{g(n)} \le |T^{[\eta]} \cap {}^m2|/2^m \le 
|T^{[\eta]} \cap {}^{g(n)}2|/2^{g(n)}. \tag"{$(*)_{n,m}$}"
$$
\enddefinition
\bigskip

\definition{\stag{pr.3} Definition}:  1) For subtrees $T_n$ of 
${}^{\omega >}2$ for
$n < \omega$ and filter $D$ on $\omega$ we say $T = \text{ lim}_D \langle
T_n:n < \omega \rangle$ if $T = \{\eta \in {}^{\omega >}2:\{n < \omega:
\eta \in T_n\} \in D\} = \{\eta \in {}^{\omega >}2:\{n < \omega:\eta \notin
T_n\} \ne \emptyset \text{ mod } D\}$.  
Similarly for $T_n \subseteq {\Cal H}(\aleph_0)$.  We
may omit $D$ if $D$ is the family co-bounded subset of $\omega$. 
Note that lim$_D\langle T_n:n < \omega \rangle$ if it exists is uniquely defined and is
absolute and if $D$ is an ultrafilter it is always well defined.
\nl
2) Let ${\Cal G}^{\bold V}$ be $\{g \in \bold V:g$ is an increasing
function from $\omega$ to $\omega$ and $g(n) > n\}$,
let ${\Cal G}$ vary on subsets of ${\Cal G}^{\bold V}$.  
Let $\bold T_{\Cal G} = \cup\{\bold T_g:g \in {\Cal G}\}$. 
Let ${\Cal G}^{\bold V}_w = \{g \in {\Cal G}^{\bold V}:\text{Dom}(g) =w\}$.
\nl
3) For Borel subsets ${\Cal B},{\Cal B}_n \, (n < \omega)$ of ${}^\omega 2$ and filter $D$
on $\omega$ let ${\Cal B} = \text{ ms-lim}_D \langle {\Cal B}_n:n < \omega
\rangle$ means that for every $\varepsilon > 0$ for the set $\{n < \omega:
\text{Leb}({\Cal B} \Delta {\Cal B}_n) < \varepsilon\}$ belongs to $D$.
Note that here ${\Cal B}$ is not uniquely determined and this has been
used in \scite{pr.1}, $(*)_{\bold r,D,D_1}$. 
\enddefinition
\bigskip

\proclaim{\stag{pr.4} Claim}  1) For subtrees $T_n$ of ${}^{\omega >}2$ for
$n < \omega$ and ultrafilter $D$ on $\omega$, lim$_D \langle T_n:n < \omega
\rangle$ is well defined and is a subtree of ${}^{\omega >}2$. \nl
2) If in addition $n < \omega \Rightarrow T_n \in \bold T_g$ (as in
Definition \scite{pr.2}), \ub{then}
\mr
\item "{$(a)$}"  lim$_D \langle T_n:n < \omega \rangle$ belongs to
$\bold T_g$
\sn
\item "{$(b)$}" if $T = \text{ lim}_D \langle T_n:n < \omega \rangle$ then
lim$(T) = \text{ ms-lim}_D \langle \text{lim}(T_n):n < \omega \rangle$.
\endroster
\endproclaim
\bigskip

\demo{Proof}  Easy.
\enddemo
\bigskip

\definition{\stag{pr.5} Definition}  1) We say $\rho \in {}^\omega 2$ is
$(N,\bold T_{\Cal G},D)$-continuous or ${\Cal G}$-continuous over $N$ for
$D$ \ub{if}:
\mr
\item "{$(a)$}"  $N \subseteq \bold V$ a transitive class, a model of ZFC,
or $\prec({\Cal H}(\chi),\in)$ for some $\chi$; or more generally, a set or
a class, which is a model of enough set theory (say ZFC$^-$) and ${\Cal H}
(\aleph_0) \subseteq N,\omega \in N$, with reasonable absoluteness and
$D \in N$ a filter on $\omega$ containing the co-finite sets (so
$D^{\bold V}$ is the filter generated in $\bold V$ by $D \cap N$) 
\sn
\item "{$(b)$}"  ${\Cal G} \in N$ (and of course ${\Cal G} \subseteq
{\Cal G}^{\bold V}$) and if $m(*) < \omega$ and for each
$m < m(*)$ we have $g_m \in {\Cal G}$, and $\langle T^m_n:n < \omega 
\rangle \in N,T^m_n \in N \cap \bold T_{g_m},T^m \in N \cap
\bold T_g$ and $T^m = \text{ lim}_D\langle T^m_n:n < \omega \rangle$, \ub{then}
$\rho \in \dbca_{m< m^*} \text{ lim}(T^m) \Rightarrow \{n:\text{if } m
< m^*$ then $\rho \in \text{ lim}(T^m_n)\} \ne \emptyset$ mod $D^{\bold V}$.
\ermn
2) We define the ideal Null$_{{\Cal G},D}$ as the $\sigma$-ideal generated
by the sets of the form 
$\{\eta \in {}^\omega 2:\rho \in \text{ lim}(T^m)$ for $m < m(*)$
but $\{n:\text{if } m < m(*)$ then $\rho \in \text{ lim}(T^m_n)\} =
\emptyset$ mod $D\}$ 
with $T^m,T^m_n \in \bold T_{g_m},\text{lim}(T^m) = 
\text{ lim}_D \langle \text{lim}(T^m_n):n < \omega \rangle$, for some
$m < m(*)$ where $m(*) < \omega,g_m \in {\Cal G}$. \nl
3) If $D$ is the filter of co-finite subsets of $\omega$, we may omit
it.
\enddefinition
\bn
\ub{\stag{pr.6} Observation}: 1) Assume ${\Cal G} \in \bold V$ is $\ne
\emptyset$ and $\bold V_1$ extends $\bold V$.  If $({}^\omega 2)^{\bold
V}$ is not in the ideal (Null$_{\Cal G})^{\bold V_1}$, \ub{then} there is no $\rho \in
({}^\omega 2)^{\bold V_1}$ which is a Cohen real over $\bold V$. \nl
2) If $D$ is an ultrafilter on $\omega$ (in $\bold V$), \ub{then} in Definition
\scite{pr.5}(1),(2) the case $m(*)=1$ suffice.
\bigskip

\demo{Proof}  1) Choose $g \in {\Cal G}$ and choose $\langle m_i:i <
\omega \rangle$ by $m_0 = 0,m_{i+1} = 2g(m_i) > m_i$.

Assume toward contradiction
$\rho \in{}^\omega 2$ is Cohen over $\bold V$.  On ${}^\alpha
2,\alpha \le \omega$ let $+$ be addition mod 2 coordinatewise.

In $\bold V$ we can find a sequence $\langle T_i:i < \omega \rangle$
of subtrees of ${}^{\omega >}2$ such that:  Leb(lim$(T_i)) \ge
1-1/2^i,{}^{n \ge}2 \subseteq T_n$ and $i \le j < \omega \and \eta
\in T_i \cap {}^{n_j} 2 \rightarrow (\exists!\nu)(\eta \triangleleft
\nu \in {}^{(m_{j+1})}2 \and \nu \notin T_i)$.  So easily $T_i \in
\bold T_g$ and $T =: \text{ lim}_D\langle T_n:n < \omega \rangle$ is
${}^{\omega >}2$.  Now if $\nu \in ({}^\omega 2)^{\bold V}$ then $\rho +
\nu$ is Cohen over $\bold V$ hence $\notin$ lim$(T_n)$ hence $\nu
\notin \rho + \text{ lim}(T_n)$.  So letting $T'_n = \{\nu + \rho
\restriction k:\nu \in T_n \cap {}^k 2,k < \omega\}$, still lim$\langle
T'_n:n < \omega \rangle$ is $T = {}^{\omega >} 2$ and for every $\nu
\in ({}^\omega 2)^{\bold V}$ we have $n < \omega \Rightarrow \nu
\notin \text{ lim}(T'_n)$ but $\nu \in \text{ lim}(T)$.  So $T,\langle
T'_n:n < \omega \rangle$ exemplify that $({}^\omega 2)^{\bold V}$ is
in Null$_{\Cal G}$, contradiction. \nl
2) Easy to check.    \hfill$\square_{\scite{pr.6}}$\margincite{pr.6} 
\enddemo
\bigskip

\demo{\stag{pr.7} Conclusion}  Assume
\mr
\item "{$(a)$}"  $\bold V_1 \supseteq \bold V$
\sn
\item "{$(b)$}"  ${\Cal G} \in \bold V,{\Cal G} \subseteq {\Cal G}^{\bold V}$
\sn
\item "{$(c)$}"  in $\bold V,D$ is a non principal ultrafilter on
$\omega$.
\sn
\item "{$(d)$}"  $\bold r \in ({}^\omega 2)^{\bold V_1}$ is 
${\Cal G}$-continuous over $\bold V$ (recall \scite{pr.5}(3)).
\ermn
\ub{Then} we can find $D_1$ such that
\mr
\item "{$(\alpha)$}" $D_1$ is a non principal ultrafilter on $\omega$
extending $D$
\sn
\item "{$(\beta)$}"  if $g \in {\Cal G}^{\bold V}$ and $T,\langle T_n:n <
\omega \rangle \in \bold V,\{T,T_n\} \subseteq \bold T_g$ and \nl
$T = \text{ lim}_D \langle T_n:n < \omega \rangle$ then
$$
\{n:(\bold r \in \text{ lim}(T)) \equiv (\bold r \in \text{ lim}(T_n))\} \in
D_1.
$$
\endroster
\enddemo
\bigskip

\demo{Proof}  To find such $D_1$, it is enough to prove
\mr
\item "{$\boxtimes$}"  assume 
{\roster
\itemitem{ $(*)_1$ }  $m^* < \omega$ and for $m < m^*,g_m \in
{\Cal G},T^m,T^m_n \in \bold T_{g_m} (\subseteq \bold V))$ for $n <
\omega$ and $\langle T^m_n:n < \omega \rangle \in \bold V$ and $T^m =
\text{ lim}_D\langle T^n_m:m < \omega \rangle$ and lastly $B \in D$.
\endroster}
\ub{Then} for some $n \in B$ we have $(\forall m < m^*)(\bold r \in
\text{ lim}(T^m) \equiv \bold r \in \text{ lim}(T^m_n))$.
\ermn
For $m < m^*$ we can find $k_m < \omega$ such that $\bold r \notin \text{
lim}(T^m) \Rightarrow \bold r \restriction k_m \notin T^m$ hence $k =
\text{ max}\{k_m:m < m^*\} < \omega$ and let $u = \{k < m^*:\bold r \in
\text{ lim}(T^m)\}$.  For each $m < m^*,m \notin u$ we know that
$\bold r \restriction k \notin T^m$ hence $A_m =: \{n:\bold r \restriction k_m
\notin T^m_n\} \in D$, and clearly $n \in A_m \Rightarrow \bold r
\restriction k_m \notin T^m_n \Rightarrow \bold r \restriction k \in
T^m_n$.  Let $B_1 = B \cap \cap \{A_m:m < m^*,m \notin u\}$, clearly $B_1
\in \bold V$ and $B_1 \in D$ hence $B_1$ is infinite.  So we can choose
by induction on $i <
\omega$, a number $n_i \in B_1$ such that $n_i > n_j$ for $j < i$ and
$m < m^* \Rightarrow
T^m_{n_i} \cap {}^{i \ge}2 = T^m \cap {}^{i \ge}2$ moreover we do this
in $\bold V$ (possible as $\bold r \restriction k_m \in \bold V$ for
$m < m^*$).  By the definition
of ``$\bold r$ is ${\Cal G}$-continuous over $\bold V$", the Borel set

$$
\align
{\Cal B} = \{&\eta \in {}^\omega 2:\eta \in \cap\{\text{lim}(T^m):m
\in u\} \text{ but} \\
  &\{i < \omega:\eta \in \cap\{\text{lim}(T^m_{n_i}):m \in u\}\}
\text{ is finite}\}\}
\endalign
$$
\mn
satisfies: $\bold r \notin {\Cal B}^{{\bold V}_1}$.  But $\bold r \in
\cap\{\text{lim}(T^m):m \in u\}$ by the choice of $u$ hence by the
definition of ${\Cal B}$ we have $\bold r \in
\cap\{\text{lim}(T^m_{n_i}):m \in u\}$, for infinitely many $i$'s.
Hence for some $i$
\mr
\item "{$(*)$}"  $m \in u \Rightarrow \bold r \in \text{ lim}(T^m_{n_i})$.
\ermn
Now if $m < m^*,m \notin u$ then $\bold r \notin \text{ lim}(T^m)$ by the
choice of $u$ and $\bold r \notin \text{ lim}(T^m_{n_i})$ as $n_i \in B_1
\subseteq A_m$, see above.  So $n_i$ is as required.
\hfill$\square_{\scite{pr.7}}$\margincite{pr.7}    
\enddemo
\newpage

\head {\S5 On iterating $\Bbb Q_{\bar D}$} \endhead  \resetall \sectno=5
\bigskip

\definition{\stag{br.1} Definition}  1)  Let {\bf IF} be the family of $\bar D =
\langle D_\eta:\eta \in {}^{\omega >} \omega \rangle$ with each $D_\eta$ a
filter on $\omega$ containing all the co-finite subsets of $\omega$. \nl
2) {\bf IUF} is the family of $\bar D = \langle D_\eta:\eta \in {}^{\omega >}
\omega \rangle$ with each $D_\eta$ a nonprincipal ultrafilter on $\omega$.
\enddefinition
\bn
On $\Bbb Q_{\bar D}$ see \cite{JdSh:321}.
\definition{\stag{br.2} Definition}  1) For $\bar D \in$ {\bf IF} we define 
${\Bbb Q}_{\bar D}$ as follows:
\mr
\item "{$(\alpha)$}"  the set of elements is $Q_{\bar D} = 
\{T:T \subseteq {}^{\omega >}
\omega$ is closed under initial segments, is nonempty and for some member
tr$(T) \in T$, the trunk, we have: $\nu \in T \and \ell g(\nu) \le
\ell g(\eta) \Rightarrow \nu \trianglelefteq \eta$ and $\eta
\trianglelefteq \nu \in T \Rightarrow \{n:\nu \char 94 \langle n \rangle
\in T\} \ne \emptyset \text{ mod } D_\nu\}$
\sn
\item "{$(\beta)$}"  $\le = \le^{{\Bbb Q}_{\bar D}}$ is the inverse of
inclusion
\sn
\item "{$(\gamma)$}"  $\le_{\text{pr}} = \le^{Q_{\bar D}}_{\text{pr}}$ is
defined by $T_1 \le_{\text{pr}} T_2 \equiv (T_2 \subseteq T_1 \and
\text{ tr}(T_1) = \text{ tr}(T_2))$
\sn
\item "{$(\delta)$}"   $\le_{apr} = \{(p,q):p \le q$ and $q =
p^{[\eta]}$ for some $\eta \in p\}$ on $p^{[\eta]}$ see below
\sn
\item "{$(\varepsilon)$}"  val$(T) = \text{ tr}(T) \in {\Cal H}(\aleph_0)$.
\ermn
2) Let $\underset\tilde {}\to \eta = 
\underset\tilde {}\to \eta({\Bbb Q}_{\bar D}) = {\underset\tilde {}\to
\eta_{\bar D}}$
be $\cup\{\text{tr}(p):p \in {\underset\tilde {}\to G_{{\Bbb Q}_{\bar D}}}\}$,
this is a ${\Bbb Q}_{\bar D}$-name of a member of ${}^\omega
\omega$. \nl
For $p \in Q_{\bar D},\eta \in p$ let $p^{[\eta]} = \{\nu \in p:\nu
\trianglelefteq \eta \vee \eta \trianglelefteq \nu\}$; so we have $p \le
p^{[\eta]} \in Q_{\bar D}$, tr$(p^{[\eta]}) \in \{\eta,tr(p)\}$.
\enddefinition
\bn
\ub{\stag{br.3} Fact}:   For $\bar D \in$ {\bf IUF}, we have:
\mr
\item "{$(a)$}"  $\Bbb Q_{\bar D}$ is a straight clear simple ${\Cal
F}$-forcing, $\sigma$-centered, purely proper if ${}^{\omega >}\omega
\subseteq {\Cal F}$
\sn
\item "{$(b)$}"   ${\Bbb Q}_{\bar D}$ is an
${\Cal F}$-psc if ${}^{\omega >} \omega \subseteq {\Cal F}$
\sn
\item "{$(c)$}"  from 
$\underset\tilde {}\to \eta[{\underset\tilde {}\to G_{{\Bbb Q}_{\bar D}}}]$ we
can reconstruct ${\underset\tilde {}\to G_{{\Bbb Q}_{\bar D}}}$ so it is a generic
real
\sn
\item "{$(d)$}"  $p \le_{apr} p^{[\eta]} \in \Bbb Q_{\bar D}$ for
$\eta \in p \in \Bbb Q_{\bar D}$.
\endroster
\bigskip

\demo{Proof}  Immediate, being $\sigma$-centered using \scite{ct.3a}(2) with
$\{R_n:n < \omega\} = \{\{p \in Q_{\bar D}:\text{tr}(p) = \eta\}:\eta \in
{}^{\omega >} \omega\}$.  \hfill$\square_{\scite{br.3}}$\margincite{br.3}
\enddemo
\bn
For completeness we prove the basic properties of $\Bbb Q_{\bar D}$.
\proclaim{\stag{br.4} Claim}  For $\bar D \in$ {\bf IUF} letting $\Bbb
Q = \Bbb Q_{\bar D}$ we have \nl
1) ${\Bbb Q}$ has pure 2-decidability. \nl
2) If $p \in \Bbb Q$ and ${\Cal I} \subseteq \Bbb Q$ is dense above
$p$, \ub{then} for some $q$ we have $p \le_{pr} q$ and $Y_0 = 
\{\eta \in p:\text{there is } q \text{ such that }
p^{[\eta]} \le_{\text{pr}} q \in {\Cal I}\}$ contains a front of
$q$ (that is $\eta \in \text{ lim}(p) \Rightarrow (\exists n)[\eta
\restriction n \in Y])$. \nl
3) If $p \in \Bbb Q$ and $Y \subseteq p$ satisfies $\eta \in Y
\and \eta \triangleleft \nu \in p \Rightarrow \nu \in Y$, \ub{then}
there is $q$ such that $p \le_{pr} q$ and: either $q \cap Y =
\emptyset$ or there is a function $h:q \backslash Y \backslash
\{\nu:\nu \triangleleft tr(q)\} \rightarrow
\omega_1$ such that for $\eta \triangleleft \nu$ in Dom$(h),h(\eta) >
h(\nu)$ and $(\forall \eta \in \text{ Dom}(h))(\forall k < \ell g(\eta))[\eta
\restriction k \in \text{ Dom}(h)]$. \nl
4) Let $p \in \Bbb Q,{\Cal I} \subseteq \Bbb Q$.  \ub{Then} ${\Cal I}$ is
predense above $p$ (in $\Bbb Q$) \ub{iff} there is $\langle(p_\eta,h_\eta):tr(p) \trianglelefteq
\eta \in p \rangle,\langle q_\eta:\eta \in Y \rangle$ such that:
\mr
\item "{$(a)$}"  $p^{[\eta]} \le_{pr} p_\eta \in \Bbb Q$
\sn
\item "{$(b)$}"  if tr$(p) \trianglelefteq \eta \in p$ then $\eta \in Y$ or
$h_\eta$ is a function with domain a subset of $\{\nu:\eta
\trianglelefteq \nu \in p_\eta\}$ closed under initial segments, range
of $h_\eta$ is $\subseteq \omega_1,h_\eta$ decreasing (i.e. $\rho
\triangleleft \nu \Rightarrow h(\rho) > h(\nu)$ when $\rho,\nu \in
\text{ Dom}(h_\eta))$ and $\nu \in \text{ Dom}(h_\eta) \wedge \nu \char
94 \langle \ell \rangle \in p_\eta \backslash \text{ Dom}(h_\eta)
\Rightarrow \nu \in Y$
\sn
\item "{$(c)$}"  $q_\eta \in {\Cal I}$ and $tr(q_\eta) \trianglelefteq \eta
\in q_\eta$.
\endroster 
\endproclaim 
\bigskip

\demo{Proof}  Straight, but we give details. \nl
1) Let $p \Vdash_{\Bbb Q} ``\underset\tilde {}\to \tau \in \{0,1\}"$.  Let $Y_0
=: \{\eta \in p:tr(p) \trianglelefteq \eta$ and there is $r \in \Bbb
Q$ 
forcing a value to $\underset\tilde {}\to \tau$ such that
$p^{[\eta]} \le_{pr} q\}$ and let $Y =: \{\eta \in p:\text{ for some }
\nu \in Y_0$ we have $tr(p) \trianglelefteq \nu \trianglelefteq
\eta\}$.  We apply part (3), (trivially $Y$ is as assumed there) so let
$p \le_{pr} q \in \Bbb Q$ be as there.  If $q \cap Y =
\emptyset$ let $r$ be such that $q \le r$ and $r$ forces a value to
$\underset\tilde {}\to \tau$; hence $tr(r) \in q \cap Y$,
contradiction.    So there is $h$ as there.  
Stipulate $h(\nu) = -1$ if $\nu \in Y$.  We prove by induction on $\alpha < \omega_1$ (and
$\alpha \ge -1$) that:
\mr
\item "{$(*)_\alpha$}"  if $tr(q) \trianglelefteq \eta \in \text{
Dom}(h)$ and $h(\eta) = \alpha$ \ub{then} there is $r = r_\eta$ such that $tr(q)
\trianglelefteq tr(r_\eta) \trianglelefteq \eta$ and $r_\eta$ forces a
value to $\underset\tilde {}\to \tau$.
\ermn
Now if $\alpha = -1$ then (by the choice of (3)), $\eta \in Y$ hence
(by the definition of $Y$) for some $\nu$ we have $tr(q)
\trianglelefteq \nu \trianglelefteq \eta$ and $\nu \in Y_0$.  Hence
(by the definition of $Y_0$) there is $r$ such that $q^{[\nu]}
\le_{pr} r \in \Bbb Q$ and $r$ forces a value to $\underset\tilde {}\to
\tau$ as required.  If $\alpha \ge 0$ for each $\ell < \omega$ such that
$\eta \char 94 \langle \ell \rangle \in q$ let $r_{\eta \char 94 <\ell>}
\Vdash_{\Bbb Q} ``\underset\tilde {}\to \tau = i_\ell"$.  If for some
such $\ell,tr(r_{\eta \char 94 <\ell>}) \trianglelefteq \eta$ we are
done, otherwise for some $i < 2$ we have $A =: \{\ell:i_\ell = i$ and
so $\eta \char 94 <\ell> \in q \text{ and } tr(r_{\eta \char 94
<\ell>}) = \eta \char 94 <\ell>\} \in D_\eta$ and let $r_\eta =
\cup\{r_{\eta \char 94 <\ell>}:\ell \in A\}$.

Having carried the induction, for $\alpha = h(tr(q)),r_{tr(q)}$ is as
required: forced a value to $\underset\tilde {}\to \tau$ and 
$tr(q) \trianglelefteq tr(r_{tr(q)}) \trianglelefteq tr(q)$ we have
$tr(r_{tr(q)}) = tr(q)$ hence $q
\le_{pr} r_{tr(q)}$ but $p \le_{pr} q$ hence $p \le_{pr}
r_{tr(q)}$. \nl
2) Let $Y = \{\eta: \text{ for some } \nu$ we have $tr(p)
\trianglelefteq \nu \trianglelefteq \eta \in p$ and $\nu \in Y_0\}$.
Apply part (3) to $p$ and $Y$ (clearly $Y$ is as required there) and
get $q$ as there.  If $q \cap Y = \emptyset$ find $r$ such that $q \le
r \in {\Cal I}$, (exists by the density of ${\Cal I}$) so by our
definitions $tr(r) \in Y_0 \subseteq Y$ and $tr(r) \in r \subseteq q$
so $q \cap Y \ne \emptyset$, contradiction.  So assume $q \cap Y \ne
\emptyset$ hence necessarily there is $h$ as there
and for every $\eta \in \text{ lim}(q)$, as $\langle h(\eta
\restriction \ell):\ell \in [\ell g(tr(q)),\omega) \rangle$ cannot be
strictly decreasing sequence of ordinals, necessarily for some $\ell
\ge \ell g(tr(q)),\eta \restriction \ell \notin \text{ Dom}(h)$ hence
$\eta \restriction \ell \in Y$ hence for some $m \in [\ell
g(tr(q)),\ell]$ we have $\eta \restriction m \in Y_0$ contradicting
our present assumption. \nl
3) Let $Z = \{\eta:tr(p) \trianglelefteq \eta \in p$ and for
$p^{[\eta]} \in \Bbb Q$ there is $h$ as required in the claim$\}$.

Clearly
\mr
\item "{$(*)_1$}"  $Y \subseteq Z \subseteq \{\eta:tr(p)
\trianglelefteq \eta \in p\}$ (use $h_\eta$ with Dom$(h_\eta) =
\{\eta\},h_\eta(\eta)=0$
\sn
\item "{$(*)_2$}"  if $tr(p) \trianglelefteq \eta \in p$ and $A =
\{\ell:\eta \char 94 \langle \ell \rangle \in Z\} \in D_\eta$ then
$\eta \in Z$. \nl
[Why?  If $h_\ell$ witness $\eta \char 94 \langle \ell \rangle \in Z$
for $\ell \in A$, let $\alpha^* = \cup\{h_\ell(\eta \char 94 \langle
\ell \rangle) +1:\ell \in A\}$ and define $h:\text{Dom}(h) = \{\nu:\nu
\trianglelefteq \eta\} \cup \cup \{\text{Dom}(h_\ell) \backslash
\{\nu:\nu \trianglelefteq \eta\}:\ell \in A\},h \restriction
(\text{Dom}(h_\ell) \backslash \{\nu:\nu \trianglelefteq \eta\})$ is
$h_\ell,h(\nu) = \alpha^* + \ell g(\eta) - \ell g(\nu)$ if $\nu
\trianglelefteq \ell g(\eta)$.]
\ermn
If $tr(p) \in Z$ we get the second possibility in the conclusion.  If
$tr(p) \notin Z$, let $q = p \backslash \{\eta \in p$: there is no
$\nu \trianglelefteq \eta$ which belongs to $Z\}$, so $\{\eta:\eta
\trianglelefteq tr(p)\} \subseteq q$ (see $Z$'s definition + present
assumption) and $q$ is closed under initial segments (read its
definition) and by $(*)_2$ we can prove by induction on $\ell \ge \ell
g(tr(p))$ that $\eta \in q \cap {}^n \omega$ implies $\{n:\eta \char
94 <n> \in q\} \in D_\eta$.  So clearly $p \le_{pr} q \in
\Bbb Q_{\bar D},q \cap Y = \emptyset$ as required.  \nl
4) Should be clear.   \hfill$\square_{\scite{br.4}}$\margincite{br.4}
\enddemo
\bn
The following is natural to note if we are interested in the Borel conjecture.
(Of course, this claim does not touch the problem of preserving this by the
later forcings in the interation we intend to use.)  Compare with \scite{bs.7}.
\proclaim{\stag{bs.5} Claim}  Assume 
\mr
\item "{$(a)$}"  $\bar D= \langle D_\eta:\eta \in {}^\omega 2 \rangle$
where \nl
$D_\eta$ is a non-principal ultrafilter on $\omega$
\sn
\item "{$(b)$}"  $N \prec ({\Cal H}(\chi),\in)$ is countable, $\bar D \in N$
\sn
\item "{$(c)$}"  $\rho_m \in {}^\omega \omega \backslash N$ for $m < \omega$
\sn
\item "{$(d)$}"  $p \in \Bbb Q_{\bar D} \cap N$.
\ermn
\ub{Then} we can find $q$ such that
\mr
\item "{$(\alpha)$}"  $p \le_{\text{pr}} q \in \Bbb Q_{\bar D}$
\sn
\item "{$(\beta)$}"  $q \Vdash$ ``if $f \in \dsize \prod_{n < \omega}
(\underset\tilde {}\to \eta(\Bbb Q_{\bar D})
[{\underset\tilde {}\to G_{{\Bbb Q}_{\bar D}}}](n) + 1)$ and
$f \in N[{\underset\tilde {}\to G_{{\Bbb Q}_{\bar D}}}]$ and $m <
\omega$ \nl

\hskip20pt \ub{then} $(\forall^* n)(\neg f \restriction n \triangleleft \rho_m)$".
\endroster
\endproclaim
\bigskip

\demo{Proof}  Let $\langle {\underset\tilde {}\to f_\ell}:\ell < \omega
\rangle$ list the $\underset\tilde {}\to f \in N$ such that 
$\Vdash_{{\Bbb Q}_{\bar D}} ``\underset\tilde {}\to f \in
\dsize \prod_{n < \omega} ((\underset\tilde {}\to \eta(\Bbb Q_{\bar D}))(n)+1)"$.
\nl
Now for every $tr(p) \trianglelefteq 
\eta \in p$ and $\ell < \omega$ there is a function
$f_{\ell,\eta} \in {}^\omega(\omega +1)$ such that: for every $k < \omega$
there is $q_{\ell,\eta,k} \in \Bbb Q_{\bar D} \cap N$ such that $q^{[\eta]}
\le_{\text{pr}} q_{\ell,\eta,k}$ and $q_{\ell,\eta,k}
\Vdash_{{\Bbb Q}_{\bar D}}$ ``for $\ell,n < k$ we have: if 
$f_{\ell,\eta}(n) < \omega$ then
${\underset\tilde {}\to f_\ell}(n) = f_{\ell,\eta}(n)$ and
if $f_{\ell,\eta}(n) = \omega$ then ${\underset\tilde {}\to f_\ell}(m) > k"$. \nl
Now $q_{\ell,\eta,k} \in N$ and \wilog \, $\langle q_{\ell,\eta,k}:\eta \in
p$ and $k < \omega \rangle,\langle f_{\ell,\eta}:\eta \in p \rangle$ belongs
to $N$ (but we cannot have $\langle q_{\ell,\eta,k}:\eta \in p,k < \omega$
and $\ell < \omega \rangle \in N$).  Now for each $\ell,\eta,k$ as
$f_{\ell,\eta} \in N$ and $\rho_m \notin N$ clearly the set $\{n:f_{\ell,\eta}(n) \ne
\rho_m(n)\}$ is infinite so let $k_0(\ell,\nu,m) = \text{ Min}\{k:
\text{if } m' < m \text{ then } f_{\ell,\nu} \restriction k \ne
\rho_{m'} \restriction \ell\}$ and $k(\ell,\nu,m) = \text{ Min}\{k:k \ge
k_0(\ell,\nu,m)$ and if $k' < k,\rho_{\ell,\nu}(k') = \omega$ then
$k > \rho_m(k')\}$. \nl
Now define $q$ as $\{\eta \in p:\text{ if } \nu \triangleleft \eta,\ell
\le \ell g(\nu) \text{ and } m \le \ell g(\eta) \text{ then } \eta \in
q_{\ell,\nu,k(\ell,\nu,m)}\}$.

The checking is straightforward.  \hfill$\square_{\scite{bs.5}}$\margincite{bs.5}
\enddemo
\bigskip

\proclaim{\stag{bs.6} Claim}  1) Assume
\mr
\item "{$(a)$}"  $\bold V \subseteq \bold V_1$
\sn
\item "{$(b)$}"  ${\Cal G} \in \bold V$ and ${\Cal G} \subseteq {\Cal
G}^{\bold V}$
\sn
\item "{$(c)$}"  $\bar D \in$ {\bf IUF}$^{\bold V}$.
\ermn
If $\bar D_1 \in$ {\bf IUF}$^{{\bold V}_1}$ and
$\eta \in {}^{\omega >} \omega \Rightarrow D_\eta \subseteq D_{1,\eta}$
\ub{then}
\mr
\item "{$(\alpha)$}"  $\Bbb Q^{\bold V}_{\bar D} \subseteq \Bbb Q_{\bar D_1}$ (so
$\Bbb Q^{\bold V}_{\bar D} = \Bbb Q_{\bar D_1} \cap \bold V$ and
$\le_{{\Bbb Q}_{\bar D}}[\bold V] = \le_{{\Bbb Q}_{\bar D_1}}
\restriction
\Bbb Q^{\bold V}_{\bar D}$ and similarly for incompatibility)
\sn
\item "{$(\beta)$}"  if in $\bold V$ we have ``${\Cal I} \subseteq
\Bbb Q_{\bar D}$ is predense over $p \in \Bbb Q_{\bar D}"$ \ub{then} in $\bold V_1$
we have ``${\Cal I} \subseteq \Bbb Q^{\bold V}_{{\bar D}_1}$ is
predense over $p$"
\sn
\item "{$(\gamma)$}"  if $G_1 \subseteq \Bbb Q_{{\bar D}_1}$ is
generic over $\bold V_1$, \ub{then} $G =: G_1 \cap \Bbb Q^{\bold
V}_{\bar D}[\bold V]$ is a generic
subset of $\Bbb Q^{\bold V}_{\bar D}$ over
$\bold V$.
\ermn
2) Assume in addition
\mr
\item "{$(d)$}"  $\bold r \in ({}^\omega 2)^{{\bold V}_1}$ is
${\Cal G}$-continuous over $\bold V$.
\ermn
If $D_{1,\eta}$ extending $D_\eta$ is chosen as in
\scite{pr.5} for each $\eta \in {}^{\omega >} \omega$, 
\ub{then} $\Vdash_{{\Bbb Q}_{{\bar D}_1}} ``\bold r$ is
${\Cal G}$-continuous over $\bold V[\underset\tilde {}\to \eta
(\Bbb Q_{{\bar D}_1})]"$.
\endproclaim
\bigskip

\demo{Proof}  1) Clause $(\alpha)$ is obvious, clause $(\beta)$ holds
by \scite{br.4}(4), and clause $(\gamma)$ follows (this is as in \cite{Sh:700}). \nl
2) So assume that $p \in \Bbb Q_{{\bar D}_1},m^* < \omega$ and for
each $m < m^*,g_m \in {\Cal G}$ and ${\underset\tilde {}\to T^m},
\langle {\underset\tilde {}\to T^m_n}:n < \omega \rangle \in \bold V$
are $\Bbb Q^{\bold V}_{\bar D}$-names hence $\Bbb Q_{{\bar D}_1}$-names
such that:
\mr
\item "{$(*)_0$}"  $\Vdash_{{\Bbb Q}_{{\bar D}_1}} ``{\underset\tilde {}\to T^m},
{\underset\tilde {}\to T^m_n} \in {\underset\tilde {}\to {\bold T}_{g_m}}$ and
${\underset\tilde {}\to T^m} = \text{ lim}\langle {\underset\tilde {}\to T_n}:
n < \omega \rangle"$.
\ermn
By the definition (\scite{pr.5}) we have to prove, for a given $n(*) <
\omega$ for some $n(*) > n(**)$ and $q$ above $p$ (in $\Bbb Q_{\bar
D}$), $q$ forces that: $m < m^* \Rightarrow \bold r \in
\lim({\underset\tilde {}\to T^n}) \equiv \bold r \in
\lim({\underset\tilde {}\to T^m_{n(*)}})$. \nl
By the definition and what we need to prove, \wilog
\mr
\item "{$(*)_2$}"  $p \Vdash ``{\underset\tilde {}\to T^m} \cap {}^{n
\ge} 2 = {\underset\tilde {}\to T^m_n} \cap {}^{n \ge}2"$ for $n <
\omega,m < m^*$. 
\ermn
We shall for some $n < \omega$ find $p' \ge p$ in $\Bbb Q_{{\bar D}_1}$ such that $p'
\Vdash ``(\bold r \in \text{ lim}({\underset\tilde {}\to T^m})) 
\equiv (\bold r \in \text{ lim}
({\underset\tilde {}\to T^m_n}))$ for every $m < m^*"$, 
this suffices (see \scite{pr.6}(2)); work in $\bold V$.  Let
$q_0 = ({}^{\omega >}\omega)$, so $q_0 \in \Bbb Q_{\bar D}$, now
we find $\langle T^m_\eta,T^m_{n,\eta}:\eta \in q_0,n < \omega \rangle$ of
course in $\bold V$ such that:
\mr
\item "{$(*)_3(i)$}"  $T^m_\eta,T^m_{n,\eta} \subseteq {}^{\omega >}
2$, for $n < \omega,m < m^*$
\sn
\item "{$(ii)$}"  for every $\eta \in q_0$ and $k < \omega$ we can find
$q^m_{\eta,k},q^m_{n,\eta,k} \in \Bbb Q_{\bar D}$ such that: \nl
$q^{[\eta]}_0 \le_{\text{pr}} q^m_{\eta,k}$, \nl
$q^{[\eta]}_0 \le_{\text{pr}} q^m_{n,\eta,k}$ \nl
$q^m_{\eta,k} \Vdash_{{\Bbb Q}_{\bar D}} ``{\underset\tilde {}\to T^m} \cap
{}^{k \ge} 2 = T^m_\eta \cap {}^{k \ge}2"$
\nl
$q^m_{n,\eta,k} \Vdash_{{\Bbb Q}_{\bar D}} ``{\underset\tilde {}\to T^m_n} \cap
{}^{k \ge} 2 = T^m_{n,\eta} \cap {}^{k \ge}2"$.
\ermn
Now clearly
\mr
\widestnumber\item{$(*)_4(iii)$}
\item "{$(*)_4$(i)}"  $T^m_\eta,T^m_{n,\eta} \in \bold T_{g_m}$
\sn
\item "{${}(ii)$}"  $T^m_\eta = \text{ lim}_{D_\eta}\langle T^m_{\eta
\char 94 <k>}:k < \omega \rangle$
\sn
\item "{${}(iii)$}"  $T^m_{n,\eta} = \text{ lim}_{D_\eta}\langle
T^m_{n,\eta \char 94 <k> \rangle}:k < \omega \rangle$.
\ermn
Next note that
\mr
\item "{$(*)_5(a)$}"   $T^m_\eta,T^m_{n,k}$ belong to $\bold T_g$
\sn
\item "{${{}}(b)$}"  $T^m_\eta = \text{ lim}\langle T^m_{n,\eta}:
n < \omega \rangle$. \nl
[Why does clause $(a)$ hold?  Let $T^m_\eta \cap {}^{g(\ell)}2 = t$ then
$q^m_{\eta,k}$ forces that ${\underset\tilde {}\to T^m} \cap
{}^{g(\ell)}2 = t$ but it also forces that ${\underset\tilde {}\to
T^m}$ satisfies the condition $(*)_\ell$ from Definition \scite{pr.2},
hence in fact it forces that $t$ satisfies it (as it gives all the
relevant information) hence $T^m_\eta$ satisfies it.  Similarly for
$T^m_{n,\eta}$.  Concerning clause (b) there is $q$ satisfying
$p^{[\eta]} \le_{pr}
q \in \Bbb Q_{\bar D}$ forcing ${\underset\tilde {}\to T^m_\eta} \cap
{}^\ell 2 = t_m,T^m_n \cap {}^\ell 2 = 
t_{m,n}$ so by $(*)_2$, if $n \ge \ell$ they
are equal.  As any two (even finitely many)  pure extensions of
$p^{[\eta]}$ are compatible, we have $T^m_\eta \cap {}^\ell 2 =
t_m,T^m_{n,\eta} \cap {}^\ell 2 = t_{m,n} = t_n$.  This is clearly
enough.]
\ermn
Hence by assumption (d) we have for $u \subseteq m^*$ and $\eta \in p$
\mr
\item "{$(*)^{u,n}_6$}"  $\bold r \in \dbca_{m \in u} \text{ lim}(T^m_\eta)$
implies that $(\exists^\infty n)(\bold r \in \dbca_{m \in u} \text{
lim } T^m_{n,\eta})$ and moreover \nl
$(\forall A \in ([\omega]^{\aleph_0})^{\bold V})(\exists^\infty n \in A)[\bold r \in
\dbca_{m \in u} \text{ lim}(T^m_{\eta,n})]$.
\ermn
But if $\bold r \notin \text{ lim}(T^m_\eta)$ then for some $k^* <
\omega,\bold r \restriction k^* \in T^m_\eta$ hence $n^* < n < \omega
\Rightarrow \bold r \restriction k^* \in T^m_{n,\eta}$ so by
$(*)_5(b)$ we have
\mr
\item "{$(*)^{m,\eta}_7$}"  if $\bold r \notin \text{ lim}(T^m_\eta)$ then
$(\exists^{< \aleph_0} n)[\bold r \in \text{ lim}(T^m_{n,\eta})]$.
\ermn
By the assumption on $D_{1,\eta}$ we have
\mr
\item "{$(*)_8(i)$}"  $\bold r \in \text{ lim}(T^m_\eta)$ iff $\bold r
\in \text{ lim}_{D_\eta} \langle T^m_{\eta \char 94 \langle k \rangle}:k <
\omega \rangle$
\sn
\item "{$(ii)$}"  $\bold r \in \text{ lim}(T^m_{n,\eta})$ iff $\bold r
\in \text{ lim}_{D_\eta} \langle T^m_{n,\eta \char 94 \langle k \rangle}:k <
\omega \rangle$.
\ermn
By $(*)_6 + (*)_7$ applied to $\eta = \text{ tr}(p)$, we can find
$n(*) > n(**)$, see $(*)_1$, such that
$(\forall m < m^*)[\bold r \in \text{ lim}(T^m_{\text{tr}(p)}) \equiv
\bold r \in \text{ lim}(T^m_{n(*),\text{tr}(p)})]$.  Next let

$$
q =: \{\nu \in p:\text{if } \ell g(tr(p)) \le 
\ell \le \ell g(\nu) \text{ and } m < m^* \text{ then }
(\bold r \in T^m_{\nu \restriction \ell}) \equiv (\text{ and } \bold r \in
T^m_{n(*),\nu \restriction \ell})\}.
$$
\mn
Now $p \le_{\text{pr}} q \in \Bbb Q_{\bar D}$ by $(*)_8$.  Lastly, let
$q^* =: \{\nu \in q:\text{if tr}(p) \trianglelefteq \nu$, \ub{then}
$\nu \in q_{\nu \restriction \ell}\}$.

Does $q^* \Vdash_{{\Bbb Q}_{{\bar D}_1}} ``\bold r \in \text{ lim}
({\underset\tilde {}\to T^m}) \equiv \bold r \in 
\text{ lim}({\underset\tilde {}\to T^m_{n(*)}})"$? if not, 
then for some $q^{**}$ we have $q^* \le q^{**}$ and
$q^{**} \Vdash_{Q_{{\bar D}_1}} ``\bold r \in \text{ lim}
({\underset\tilde {}\to T^m}) \equiv r \notin \text{ lim}
({\underset\tilde {}\to T^m_{n(*)}})"$ and moreover, for some $k$ we have
$q^{**} \Vdash_{{\Bbb Q}_{{\bar D}_1}} ``\bold r \restriction k \in
{\underset\tilde {}\to T^m} \equiv \bold r \restriction k \notin \text{ lim}
({\underset\tilde {}\to T^m_{n(*)}})"$.  But $q^{**},q_{\text{tr}(q^{**}),k},
q_{n(*),\text{tr}(q^{**}),k}$ are compatible having the same trunk, so let
$q'$ be a common upper bound with tr$(q') = \text{ tr}(q^{**})$ and we get a
contradiction.  \hfill$\square_{\scite{bs.6}}$\margincite{bs.6}
\enddemo
\bigskip

\proclaim{\stag{bs.7} Claim}  1) Assume
\mr
\item "{$(a)$}"  $\bar D \in$ {\bf IUF}
\sn
\item "{$(b)$}"  $\bar D \in N \prec ({\Cal H}(\chi),\in)$
\sn
\item "{$(c)$}"  $\rho \in {}^\omega 2 \backslash N$
\sn
\item "{$(d)$}"  $\underset\tilde {}\to \eta = {\underset\tilde {}\to
\eta_{{\Bbb Q}_{\bar D}}}$.
\ermn
\ub{Then}
\mr
\item "{$(\alpha)$}"  $\Vdash_{{\Bbb Q}_{\bar D}}$ ``if $f \in N[
{\underset\tilde {}\to G_{{\Bbb Q}_{\bar D}}}]$ and $f \in \dsize
\prod_{n < \omega} {}^{\underset\tilde {}\to \eta(n)} 2$ then
$(\forall^*n)(\neg f(n) \triangleleft \rho)"$, \nl
 moreover
\sn
\item "{$(\beta)$}"  $\Vdash_{{\Bbb Q}_{\bar D}} ``\text{if }f \in
N[G_{{\Bbb Q}_D}],f$ a function with domain $\omega,f(n) \subseteq
{}^{\underset\tilde {}\to \eta(n)} 2$ and $|f(n)| < \underset\tilde
{}\to \eta(n-1)$ when $n \ge 0$, then $(\forall^*n)(\rho \restriction \underset\tilde {}\to \eta(n)
\notin f(n))$ stipulating $\underset\tilde {}\to \eta(-1)=1$.
\endroster
\endproclaim
\bigskip

\remark{\stag{bs.7a} Remark}  1)  Instead $\bar D \in N \prec ({\Cal
H}(\chi),\in)$ it is enough to have assumptions like
\cite{Sh:630}. \nl
2) This may be relevant to trying to get universes satisfying the Borel conjecture. 
\endremark
\bigskip

\demo{Proof}  Let $\underset\tilde {}\to f \in N$ be such that
$\Vdash_{{\Bbb Q}_{\bar D}} ``\underset\tilde {}\to f$ is a function
with domain $\omega$ such that $|f(n)| \le \underset\tilde {}\to
\eta(n-1)$ and $\emptyset \ne \underset\tilde {}\to f(n) \subseteq
2^{\underset\tilde {}\to \eta(n)}"$ and we shall prove that

$$
\Vdash_{{\Bbb Q}_{\bar D}} ``\rho \restriction \underset\tilde {}\to
\eta(n) \notin f(n) \text{ for every } n < \omega \text{ large
enough}".
$$
\mn
This clearly suffices.  For each 
$n > 0$ and $\nu \in {}^{n+1} \omega$ we can find
$q_\nu \in \Bbb Q_{\bar D}$ and $\rho^m_\nu$ for $m < \nu(n-1)$, 
such that tr$(q_\nu) = \nu$
and $q_\nu \Vdash_{{\Bbb Q}_{\bar D}} ``\underset\tilde {}\to f(n) =
\{\rho^m_\nu:m < \underset\tilde {}\to \eta(n-1)\}"$.  Note
that $q_\nu \Vdash ``k < \ell g(\nu) \Rightarrow
\underset\tilde {}\to \eta(n) = \nu(k)$ in particular for $k=n-1"$ and $\Vdash ``\underset\tilde
{}\to f(n) \subseteq 2^{\underset\tilde {}\to \eta(n)},1 \le |f(n)|
\le \underset\tilde {}\to \eta(n-1))"$ hence $\rho^\ell_\nu \in
{}^{\nu(n)}2$.  
As $\underset\tilde {}\to f \in N$ \wilog \,
$\langle(q_\nu,\rho^m_\nu):m < \nu(n-1),n < \ell g(\nu)$ and
$\nu \in {}^{n+1} \omega,n < \omega$ belongs to $N$.  Now
for each $\nu \in {}^{\omega >} \omega,m < \nu(\ell g(\nu)-1) < k$ clearly
$\rho^m_{\nu \char 94 <k>} \in {}^k 2$, so for every $\ell <
\omega$ for some $\rho^m_{\nu,\ell} \in {}^\ell 2$ we have $\{k <
\omega:\rho^m_{\nu \char 94 <k>} \restriction \ell = \rho^m_{\nu,\ell}$
and $k > \ell\}$ belongs to $D_\nu$, and clearly $\rho^m_{\nu,\ell}
\triangleleft \rho^m_{\nu,\ell +1}$ and let $\rho^m_{\nu,*} =
\dbcu_{\ell < \omega} \rho^m_{\nu,\ell}$ so $\rho^m_{\nu,*} \in N \cap
{}^\omega 2$ hence $\rho^m_\nu \ne \rho$ so for some $\ell(\nu,m) <
\omega$ we have $\rho^m_\nu \restriction \ell(\nu,m) \ne \rho
\restriction \ell(\nu,m)$ hence $\{k:(\exists m)\rho^m_{\nu \char 94
<k>} \triangleleft \rho\} = \emptyset$ mod $D$.

Let $p \in \Bbb Q_{\bar D}$, let us define

$$
\align
q = \{\nu:&\nu \trianglelefteq tr(p) \text{ or } tr(p) \triangleleft
\nu \in p \text{ and if } k \in [\ell g(tr(p)),\ell g(\nu)) \\
  &\text{ and } m < \nu(k-1) \text{ then } \neg \rho^m_{\nu \restriction
(k+1)} \triangleleft \rho)\}
\endalign
$$
\mn
is a condition above $p$ forcing $\{n:\underset\tilde {}\to f(n)
\triangleleft \rho\}$ is bounded by $\ell g(tr(p))$, so we are
done. \nl
${{}}$  \hfill$\square_{\scite{bs.7}}$\margincite{bs.7}
\enddemo
\bigskip

\proclaim{\stag{bs.8} Claim}  Assume
\mr
\item "{$(a)$}"  $\bar D =$ {\bf IUF}
\sn
\item "{$(b)$}"  $D$ is a non principal ultrafilter on $\omega$
\sn
\item "{$(c)$}"  for some $n^*$, if $(*)$ below holds,
\ub{then} $D \nsubseteq \cup\{D'_{\nu,m}:\nu \in {}^{\omega >} \omega
\backslash {}^{n^* \ge} \omega$ and $m < \nu(\ell g(\nu)-1)\}$
{\roster
\itemitem{ $(*)$ }   for $\nu \in {}^{\omega >} \omega
\backslash {}^{n^* \ge} \omega$ we have $D'_{\nu,m} \le_{RK} D_\nu$ by a finite to
one function, and, of course $D'_{\nu,m}$ a non principal ultrafilter
on $\omega$.
\endroster}
\ermn
\ub{Then} in $\bold V^{{\Bbb Q}_{\bar D}}$ we have:
\mr
\item "{$(*)$}"  if $w_n \subseteq [\underset\tilde {}\to \eta(n),
\underset\tilde {}\to \eta(n+1)),|w_n| \le
\underset\tilde {}\to \eta(n)$, \ub{then} $\cup\{w_n:n < \omega\}$
is disjoint to some member of $D$.
\endroster
\endproclaim
\bigskip

\remark{Remark}  This is relevant to trying to get no $Q$-point.
\endremark
\bigskip

\demo{Proof}  Without loss of generality $n^* = 0$ (just fix
$\underset\tilde {}\to \eta \restriction n^*$) and assume
$\Vdash_{{\Bbb Q}_{\bar D}} ``{\underset\tilde {}\to {\bar w}} =
\langle {\underset\tilde {}\to w_n}:n < \omega \rangle,
{\underset\tilde {}\to w_n} \subseteq
\underset\tilde {}\to \eta(n+1),|w_n| \le \underset\tilde {}\to \eta(n)"$.  
Without loss of generality $\Vdash_{{\Bbb Q}_{\bar D}} ``w_n \ne \emptyset"$.
For $\nu \in {}^{n+2} \omega$ let $q_\nu \in
{\Bbb Q}_{\bar D}$ be such that tr$(q_\nu) = \nu$ and $q_\nu \Vdash
``{\underset\tilde {}\to w_n} = \{t^m_\nu:m < \nu(n)\}"$.  For $m < \nu(\ell g(\nu)-1)$, let 
$D'_{\nu,m} = \{A \subseteq \omega:\{k:t^m_{\nu^n
\char 94 <k>} \in A\} \in D_\nu\}$; clearly $D'_{\nu,m}$ is a non
principal ultrafilter on $\omega$ which is $\le_{RK} D_\nu$ as
$\nu(\ell g(\nu)-1) < t^m_{\nu \char 94 <k>} < k$. \nl
The rest should be clear.. \hfill$\square_{\scite{bs.8}}$\margincite{bs.8}
\enddemo
\bigskip

Results here are used in the next section; formally we have to specialize them
as ${\Bbb Q}_0$ is just $j$ random reals forcing. 

For preservation, including cardinals not collapsed we use \S2 or \S3 (really
more explicit version).  
\demo{\stag{c.1} Hypothesis}  
\roster
\item "{$(a)$}"  $\bold V \models CH$
\sn
\item "{$(b)$}"  ${\Cal F}^*$ is the standard trunk controller
\sn
\item "{$(c)$}"  ${\frak K}(0)$ is a family whose elements we denote by
$\bar R$ and Lim$(\bar R)$ is a c.c.c. forcing notion possibly
with extra structure
\sn
\item "{$(d)$}"  $<_{{\frak K}(0)}$ is a partial order on ${\frak K}(0)$
such that $\bar R' \le_{{\frak K}(0)} \bar R'' \Rightarrow \text{ Lim}
(\bar R') \lessdot \text{ Lim}(\bar R'')$; recall that ``$\kappa$-closed" 
means every increasing sequence of length $< \kappa$ has an upper bound.  
\ermn
We say ${\frak K}(0)$ is $\theta$-exactly closed if for $\langle \bar
R^\ell:i < \theta \rangle,\le_{{\frak K}(0)}$-increasing there is
$\bar R,i < \theta \Rightarrow \bar R^i \le_{{\frak K}(0)} \bar R$ and
Lim$(\bar R) = \dbcu_{i < \theta} \text{Lim}(\bar R_i)$.
\enddemo
\bigskip

\definition{\stag{c.2} Definition}  1) For an ordinal $\alpha > 0$
(assuming ${\Cal F}^*$ is based on some $\alpha^* \ge \alpha$) let
${\frak K}_\alpha$ be the family of $\bar{\Bbb Q}$ such that:
\mr
\item "{$(a)$}"  $\bar{\Bbb Q}$ is an 
${\Cal F}^*$-iteration of length $\alpha$
\sn
\item "{$(b)$}"  $\Bbb Q_0$ is a c.c.c. forcing notion from ${\frak K}(0)$,
i.e. it is Lim$(\bar R),\bar R \in {\frak K}(0)$, in principle Lim$(\bar R)$
may not determine $\bar R$ uniquely but we shall ignore this  \nl
[neater is to demand that it satisfies a strong version but not really
necessarily]
\sn
\item "{$(c)$}"  if $0 < \beta < \alpha$ \ub{then} 
${\underset\tilde {}\to {\Bbb Q}_\beta}$ is 
${\underset\tilde {}\to {\Bbb Q}_{\bar{\underset\tilde {}\to D_\beta}}}$ where
$\Vdash_{P_\beta} ``\bar{\underset\tilde {}\to D_\beta} \in$ {\bf IUF}"
(on ${\underset\tilde {}\to {\Bbb Q}_{\bar D}}$ see \scite{br.2}, on
{\bf IUF}, see \scite{br.1}).
\ermn
2) Let ${\frak K} = \dbcu_\alpha {\frak K}_\alpha$ and ${\frak K}_{< \alpha} =
\cup\{{\frak K}_\beta:\beta \le \alpha\}$ and ${\frak K}_{\le \alpha} =
{\frak K}_{< \alpha +1}$. \nl
We use $\Bbb P_\alpha = \text{ Lim}_{\Cal F}(\bar Q \restriction \alpha)$, so
e.g. $\Bbb P^1_\alpha = \text{ Lim}_{\Cal F}
(\bar{\Bbb Q}^1 \restriction \beta)$. 
\enddefinition
\bigskip

\definition{\stag{c.3} Definition}  1) For $\bar{\Bbb Q}_1,\bar{\Bbb Q}_2 
\in {\frak K}$ let $\bar{\Bbb Q}_1 \le_{\frak K} \bar{\Bbb Q}_2$ if:
\mr
\item "{$(a)$}"  $\ell g(\bar{\Bbb Q}_1) \le \ell g(\bar{\Bbb Q}_2)$
\sn
\item "{$(b)$}"  for $\beta < \ell g(\bar{\Bbb Q}_1)$ we have 
$\Bbb P_{1,\beta} \lessdot \Bbb P_{2,\beta}$, i.e.
Lim$_{\Cal F}
(\bar{\Bbb Q}_1 \restriction \beta) \lessdot \text{ Lim}_{\Cal F}(\bar{\Bbb Q}_2
\restriction \beta)$
\sn
\item "{$(c)$}"  for $\beta < \ell g(\bar{\Bbb Q}_1),\beta \ne 0$ 
and $\eta \in {}^{\omega >}
\omega$ we have $\Vdash_{\text{Lim}_{\Cal F}(\bar{\Bbb Q}_2
\restriction \beta} ``{\underset\tilde {}\to D_{1,\beta,\eta}}
\subseteq {\underset\tilde {}\to D_{2,\beta,\eta}}"$
\sn
\item "{$(d)$}"  if $\ell g(\bar{\Bbb Q}_1) = \beta < \ell g(\bar{\Bbb Q}_2)$ \ub{then}
Lim$_{\Cal F}(\bar{\Bbb Q}_1 \restriction \beta) \lessdot \text{ Lim}_{\Cal F}
(\bar{\Bbb Q}_2 \restriction \beta)$.
\endroster
\enddefinition
\bigskip

\proclaim{\stag{c.3p} Claim}  Assume $\bar{\Bbb Q}_1 \le \bar{\Bbb
Q}_2$ are from ${\frak K}_\alpha$. \nl
1) If $\alpha$ is not a limit ordinal \ub{then}

$$
\text{Lim}_{\Cal F}(\bar{\Bbb Q}_1) \lessdot \text{ Lim}_{\Cal
F}(\bar{\Bbb Q}_2).
$$
\mn
2) If $\alpha$ is a limit ordinal, ${\Cal G} \subseteq {\Cal G}^{\bold
V},\beta < \alpha,\underset\tilde {}\to \nu$ is a $\Bbb
P_{2,\beta}$-name such that, for every $\gamma \in [\beta,\alpha)$ we
have $\Vdash_{{\Bbb P}_{2,\gamma}} ``\underset\tilde {}\to \nu$ is
$\bold T_{\Cal G}$-continuous over $\bold V^{{\Bbb P}_{1,\gamma}}"$
and $\Bbb P_{1,\alpha} = \text{ Lim}_{\Cal F}(\bar{\Bbb Q}_1) \lessdot
\Bbb P_{2,\alpha} = \text{ Lim}(\bar{\Bbb Q}_2)$, \ub{then} 
$\Vdash_{{\Bbb P}_{2,\alpha}} ``\underset\tilde {}\to \nu$ is
$\bold T_{\Cal G}$-continuous over $\bold V^{{\Bbb P}_{1,\alpha}}"$.
\endproclaim
\bigskip

\demo{Proof}  1) By \scite{it.6}. \nl
2) By diagonalizing.  \hfill$\square_{\scite{c.3p}}$\margincite{c.3p}
\enddemo
\bn
The following is intended to help mainly in chain conditions, but at
present not used.
\proclaim{\stag{c.7} Claim}   Let $\bar{\Bbb Q} \in {\frak K}_\alpha$, \nl
1) We define cr$(\bar{\Bbb Q})$ as the family of objects 
${\frak p}$ consisting of:
\mr
\item "{$(a)$}"  $w[{\frak p}] \in [\alpha]^{\le \aleph_0}$ and let
$\gamma({\frak p}) = \cup\{\beta +1:\beta \in w({\frak p})\}$
\sn
\item "{$(b)$}"  val$({\frak p}) \in {\Cal F}$, Dom$(p) = w[{\frak p}]$
\sn
\item "{$(c)$}"  for $\beta \in w[{\frak p}] \cup \{\gamma[{\frak p}]\}
\backslash \{0\}$, a countable subset ${\Cal P}_\beta[{\frak p}]$ of
$\{q:q \in \Bbb P_\beta$, Dom$(q) = w({\frak p}) \cap \beta$ and 
${\Cal F} \models$ {\rm val}$({\frak p}) \le^{\Cal F}_{\text{apr}}$
val$(q)\}$, so the minimal such $\beta$ we get a subset of $\Bbb Q_0$
\sn
\item "{$(d)$}"  for $\beta \in w[{\frak p}] \cup \{\gamma[{\frak p}]\}$, a
countable family $\bold \tau_\beta[{\frak p}]$ of $\Bbb P_\beta$-names
$\underset\tilde {}\to \tau$ of a member of $\{$true,false$\}$ and for each
$\underset\tilde {}\to \tau \in \tau_\beta[{\frak p}]$ we have a set
${\Cal I}_{\underset\tilde {}\to \tau}[{\frak p}] \subseteq 
{\Cal P}_\beta[{\frak p}]$ such that 
$q \in {\Cal I}_{\underset\tilde {}\to \tau}[{\frak p}] \Rightarrow q$
forces a value to $\underset\tilde {}\to \tau$
\sn
\item "{$(e)$}"  for $\beta \in w[{\frak p}] \cup \{\gamma[{\frak
p}]\}$ and $\underset\tilde {}\to \tau \in \tau_\beta[{\frak p}]$ the set
${\Cal I}_{\underset\tilde {}\to \tau}[{\frak p}]$ is a predense subset of
$\{q \in {\Bbb P}_\beta:\text{val}(q) \in {\Cal P}_\beta[{\frak p}]\}$
\sn
\item "{$(f)$}"  if $\beta \in w[{\frak p}] \cup \{\gamma[{\frak p}]\}
\backslash \{0\}$ and $q \in {\Cal P}_\beta[{\frak p}]$ and $\gamma \in
w[{\frak p}] \cap \beta$ and $\eta \in {}^{\omega >} \omega$ \ub{then} each
``truth value$(\eta \in q(\gamma))$" belongs to $\tau_\gamma[{\frak
p}]$.  [So ${\frak p}$ involves a countable subset of $\Bbb Q_0$.] 
\ermn
2) We say ${\frak p}_1,{\frak p}_2$ (which are in $cr(\bar{\Bbb Q})$,
or more generally ${\frak p}_\ell \in cr(\bar{\Bbb Q}^\ell)$ with the
obvious changes) are strongly isomorphic as
witnessed by the function $h$ if:
\mr
\item "{$(a)$}"  $w[{\frak p}_1] = w[{\frak p}_2]$
\sn
\item "{$(b)$}"  $h \restriction {\Cal P}_\beta[{\frak p}_1]$ is a 1-to-1
mapping from ${\Cal P}_\beta[{\frak p}_1]$ onto ${\Cal P}_\beta[{\frak p}_2]$
\sn
\item "{$(c)$}"  $h \restriction \bold \tau_\beta[{\frak p}]$ is a one-to-one
mapping from $\bold \tau_\beta[{\frak p}_1]$ onto $\bold \tau_\beta
[{\frak p}_2]$
\sn
\item "{$(d)$}"  for $q \in {\Cal P}_\beta[{\frak p}_1]$ we have
val$(h(q)) = \text{ val}(q)$
\sn
\item "{$(e)$}"  if $\underset\tilde {}\to \tau \in \bold \tau_\beta
[{\frak p}_1]$ then $h$ maps 
${\Cal I}_{\underset\tilde {}\to \tau}[{\frak p}_1]$ onto ${\Cal
I}_{h({\underset\tilde {}\to \tau})}[{\frak p}_2]$ that is
${\Cal I}_{h(\underset\tilde {}\to \tau)}
[{\frak p}_2] = \{h(q):q \in {\Cal I}_{\underset\tilde {}\to \tau}
[{\frak p}_1]\}$
\sn
\item "{$(f)$}"  if $\bold t$ is a truth value,
$\underset\tilde {}\to \tau \in \bold \tau_\beta
[{\frak p}_1]$ and $q \in {\Cal I}_\tau[{\frak p}_1]$ and $q \Vdash ``
\underset\tilde {}\to \tau = \bold t"$ then $h(q) \Vdash
``h(\underset\tilde {}\to \tau) = \bold t"$.
\ermn
3) We omit the ``strongly" if we replace in part (2) clause (a) by:
\mr
\item "{$(a)'$}"  $h$ is an order preserving map from $w[{\frak p}_1]$
onto $w[{\frak p}_2]$ and such that $0 \in [{\frak p}_1]
\Leftrightarrow 0 \in w[{\frak p}_2]$.
\endroster
\endproclaim
\bigskip

\definition{\stag{c.8} Definition}  1) By induction on $\alpha \ge 1$ we 
define ${\frak K}^+_\alpha$ as the family of $\bar{\Bbb Q}_1 \in 
{\frak K}_\alpha$ such that if $\bar{\Bbb Q}_2 \in {\frak K}_\alpha$ and 
$\bar{\Bbb Q}_1 \le_{\frak K} \bar{\Bbb Q}_2$, \ub{then} for every
${\frak p}_2 \in \text{ cr}(\bar{\Bbb Q}_2)$ there is ${\frak p}_1 \in
\text{cr}(\bar{\Bbb Q}_1)$, strongly isomorphic to ${\frak p}_2$ say
as witnessed by the identity function such that
\mr
\item "{$(*)$}"  if $\beta \in w[{\frak p}] \cup \{\gamma[{\frak p}]\}
\backslash \{0\}$ is minimal and $r_1 \in {\Cal P}_\beta[{\frak p}]$
and $r_2 = h(r_1)$ then $r_2,r_1$ are compatible in $\Bbb Q_0$ whenever
$r_1 \le r \in \Bbb Q_0$.
\ermn
2) Let ${\frak K}^+ = \dbcu_\alpha {\frak K}^+_\alpha$.
\enddefinition
\bigskip

\remark{\stag{c.9} Remark}  Alternatively ${\frak K}^+_\alpha$ is the family
of $\bar{\Bbb Q}_1 \in {\frak K}^\pm_\alpha = \{\bar{\Bbb Q} \in 
{\frak K}_\alpha:(\forall \beta)(1 \le \beta < \alpha \rightarrow 
\bar{\Bbb Q} \restriction \beta \in {\frak K}^+_\beta)\}$ such that 
$\bar{\Bbb Q}_1 \le^*_0 \bar{\Bbb Q}_2 \in {\frak K}^\pm_\alpha 
\Rightarrow \text{ Lim}_{\Cal F}(\bar{\Bbb Q}_1) \lessdot 
\text{ Lim}_F(\bar{\Bbb Q}_2)$.
\endremark
\bn
\ub{\stag{c.10} Observation}:  1) If $\bar{\Bbb Q}_1 \le_{\frak K} \Bbb
Q_2$ are both from ${\frak K}_{\alpha +1}$ then Lim$_{\Cal
F}(\bar{\Bbb Q}_1) \lessdot \text{ Lim}_{\Cal F}(\bar{\Bbb Q}_2)$. \nl
2) $\le_{\frak K}$ is a partial order on ${\frak K}$.
\nl
3) If $\bar{\Bbb Q}_1,\bar{\Bbb Q}_2 \in {\frak K}^+_\alpha$ and 
$\bar{\Bbb Q}_1 \le_{\frak K} \bar{\Bbb Q}_2$ \ub{then} Lim$_{\Cal F}
(\bar{\Bbb Q}_1) \lessdot \text{ Lim}_{\Cal F}(\bar{\Bbb Q}_2)$. \nl
4)  If $\bar{\Bbb Q}_1 \le_{\frak K} \bar{\Bbb Q}_2 \in {\frak K}_\alpha$ 
and $\Bbb Q_{1,0} = \Bbb Q_{2,0}$ and
$\bar{\Bbb Q}_1 \in {\frak K}^+_\alpha$ then $\bar{\Bbb Q}_2 
\in {\frak K}^+_\alpha$ (in fact $\bar{\Bbb Q}_2 = \bar{\Bbb Q}_1$). \nl
5) If $\bar{\Bbb Q}_i \in {\frak K}^+_\alpha$ for $i < \delta$ is
$\le_{\frak K}$-increasing, $\delta < \kappa$ and ${\frak K}(0)$ is
$\kappa$-closed \ub{then} there is $\bar{\Bbb Q} \in {\frak K}^+_\alpha$ 
such that $i < \delta \Rightarrow \bar{\Bbb Q}_i \le_{\frak K} \bar{\Bbb Q}$. \nl
6) If in (5) if $K(0)$ is cf$(\delta)$-exactly closed, cf$(\delta) >
\aleph_0$, then we can add Lim$_{\Cal F}(\bar{\Bbb Q}) = \dbcu_{i < \delta}
\text{ Lim}_{\Cal F}(\bar{\Bbb Q}_i)$. \nl
7) If $\bar{\Bbb Q} \in {\frak K}_\alpha,(\forall \beta < \kappa)(|\beta|
^{\aleph_0} < \kappa =  \text{ cf}(\kappa))$ and 
${\frak K}(0)$ is $\kappa$-closed, 
each member of ${\frak K}_\alpha(0)$ is of cardinality $<
\kappa,\alpha < \kappa,
\aleph_0 < \theta = \text{ cf}(\theta) <
\kappa,{\frak K}(0)$ is exactly $\theta$-closed and $|{\Cal F}| <
\kappa$ \ub{then} there is $\bar{\Bbb Q}' \in {\frak K}^+_\alpha$ 
such that $\bar{\Bbb Q} \le_{\frak K}
\bar{\Bbb Q}'$ (we can normally bound the cardinality of 
Lim$_{\Cal F}(\bar{\Bbb Q}')$.
\bigskip

\demo{Proof}  By induction on $\alpha$, quite straightforward.

Anyhow we do not use this.  \hfill$\square_{\scite{c.10}}$\margincite{c.10}
\enddemo 
\newpage

\head {\S6 \\
On a relative of Borel conjecture with large ${\frak b}$}
\endhead  \resetall 
\bigskip

\demo{\stag{bt.5} Hypothesis}
\mr
\item "{$(a)$}"  $\bold V \models CH$
\sn 
\item "{$(b)$}"  cf$(\lambda) = \lambda,(\forall \alpha <
\lambda)(|\alpha|^{\aleph_0} < \lambda),S \subseteq \lambda$ is stationary,
$(\forall \delta \in S)(\text{cf}(\delta) > \aleph_0)$
\sn
\item "{$(c)$}"  ${\Cal F}$ is a rich enough standard trunk controller
(e.g. full).
\endroster
\enddemo
\bn
We now specify the ${\frak K}$ from \S5.
\definition{\stag{bt.8} Definition}  1) Let ${\frak K}(0)$ be the family of
$\{$Random$_A:A \subseteq \lambda\}$ where 
$d(A) = \{\omega \alpha +n:\alpha \in A,n < \omega\}$ and Random$_A$ is the 
family of Borel subsets of ${}^{d(A)}2$ of positive Lebesgue measure.  
Let ${\underset\tilde {}\to \nu_\alpha} = \cup\{f:f$ a finite function from
$[\omega \alpha,\omega \alpha + \omega)$ to $\{0,1\}$ such that $[f] =
\{g \in {}^A 2:f \subseteq g\}$ belongs to the generic$\}$.
Let $A({\Bbb Q}) = A$ if ${\Bbb Q} = \text{ Random}_A$.  
Let $\text{Random}_A \le_{{\frak K}(0)}
\text{ Random}_B$ if $A \subseteq B$ hence $\text{Random}_A \lessdot
\text{Random}_B$. 
(So $\le_{pr}$ will be just equality, $\le_{apr}$ will be the usual order).
\nl
2) For $\alpha \ge 1$, let ${\frak K}_\alpha$ be defined as in \scite{c.2}.
\nl
3) We define for any ordinal $\alpha$ and $\ell < 2$ 
the class ${\frak K}'_{\ell,\alpha} \subseteq
{\frak K}_\alpha$ as the class of $\bar{\Bbb Q}$ such that:
\mr
\item "{$(a)$}"  $\bar{\Bbb Q}$ is an ${\Cal F}$-iteration
\sn
\item "{$(b)$}"  $\ell g(\bar{\Bbb Q}) = \alpha$
\sn
\item "{$(c)$}"  ${\Bbb Q}_0 \in {\frak K}(0)$ and $A[\Bbb Q_0] \in
[\lambda]^{< \lambda}$ if $\alpha < \lambda \and \ell=0$ and $A[\Bbb Q_0] =
\lambda$ if $\ell =1$ 
\sn
\item "{$(d)$}"  if $0 < \beta < \alpha$ then 
${\underset\tilde {}\to {\Bbb Q}_\beta} = {\Bbb Q}
(\bar{\underset\tilde {}\to D_\beta})$ where $\Vdash_{{\Bbb P}_\beta} ``
{\underset\tilde {}\to {\bar D}_\beta} = \langle 
{\underset\tilde {}\to D_{\beta,\eta}}:\eta \in {}^{\omega >} \omega
\rangle \in$ {\bf IUF}"
\sn
\item "{$(e)$}"  $\bar{\Bbb Q} \restriction \gamma \in {\frak K}^+_\gamma$  for
every $\gamma < \alpha$ where ${\frak K}^+_\gamma$ is
defined in \scite{c.8} for our particular case.
\ermn
3A) If we omit $\ell$, we mean $\ell =0$ when $\alpha <
\lambda$ and we mean $\ell=1$ when $\alpha \ge \lambda$.  We let ${\frak K}'_\ell =
\cup\{{\frak K}'_{\ell,\alpha}:\alpha$ an ordinal$\}$ and ${\frak K}'
= \cup\{{\frak K}'_\alpha:\alpha$ an ordinal$\}$.
\nl
4) For $\ell = 0,1$, we define a partial order 
$\le_{{\frak K}'_\ell}$ on ${\frak K}'_\ell$ by:
\nl
$\bar{\Bbb Q}^1 \le_{{\frak K}'_\ell} \bar{\Bbb Q}^2$ \ub{iff}
$\bar{\Bbb Q}^1 \le_{{\frak K}_\ell} \bar{\Bbb Q}^2$ and $\ell
g(\bar{\Bbb Q}^1) < \ell g(\bar{\Bbb Q}^2) 
\Rightarrow [\text{there is } 
\bar{\Bbb Q}' \in {\frak K}^+_{\ell g(\bar{\Bbb Q}')}$ such that
$\bar{\Bbb Q}^1 \le_{{\frak K}_\alpha} \bar{\Bbb Q}' \le_{{\frak K}_\alpha} \bar{\Bbb Q}^2]$ and
\mr
\item "{$(*)$}"  if $\gamma$ is the minimal member of $A(\bar{\Bbb Q}^2) \backslash
A(\bar{\Bbb Q}^1)$ \ub{then} (so if $\bar{\Bbb Q}^1 = \bar{\Bbb Q}^2$ this
holds vacuously) \nl
$\Vdash_{\text{Lim}_{\Cal F}(\bar{\Bbb Q}^2)} ``
{\underset\tilde {}\to \nu_\gamma}$ is ${\Cal G}^{\bold V}$-continuous over
$\bold V^{\text{Lim}(\bar{\Bbb Q}^1)}"$ (see \scite{pr.5}(1),
(3)). 
\ermn
4A) We similarly define the partial order $\le_{{\frak K}'}$ on
${\frak K}'_\alpha$ and let $\le_{{\frak K}'_{\ell,\alpha}} =
\le_{{\frak K}'_\ell} \restriction {\frak K}'_\alpha$. \nl
5) Let ${\frak K}''_\alpha$ be the family of $\bar{\Bbb Q} \in {\frak K}'_\alpha$
such that:
\mr
\item "{$(e)$}"  $A(\bar{\Bbb Q}) = \lambda$ if $\alpha \ge
\lambda,A(\bar{\Bbb Q}) \in [\lambda]^{< \lambda}$ if $\alpha < \lambda$
\sn
\item "{$(f)$}"  if $\alpha \ge \lambda,
\beta \in (0,\alpha),\eta \in {}^{\omega >} \omega$ and
in $\bold V^{{\Bbb P}_\beta}$ letting $D_{\beta,\eta}$ be the interpretation of the
$\Bbb P_\beta$-name ${\underset\tilde {}\to D_{\beta,\eta}}$, 
we have: if $g \in {\Cal G}^{\bold V}$ and $T,T_n \in \bold T^{\bold
V[{\Bbb P}_\beta]}_g,T = \text{ lim}\langle T_n:n < \omega \rangle$,
\ub{then} for a club of $j \in S$, we 
have ${\underset\tilde {}\to \nu_\alpha} \in
\text{ lim}(T) \Rightarrow (\exists^\infty n){\underset\tilde {}\to \nu_\alpha} \in 
\text{lim}(T_n)$.
\endroster
\enddefinition
\bigskip

\proclaim{\stag{6.2abt8} Claim}  1) The two place relations
$\le_{{\frak K}'_\alpha},\le_{{\frak K}'_{1,\alpha}},\le_{{\frak
K}'_\alpha}$ are partial orders. \nl
2) The two place relation $\le_{{\frak K}'_0},\le_{{\frak
K}'_1},\le_{{\frak K}'}$ are partial orders. \nl
3) Assume
\mr
\item "{$(a)$}"  $\delta$ is a limit ordinal
\sn
\item "{$(b)$}"   $\bar{\Bbb Q}^1_1,\bar{\Bbb Q}^2_2 \in {\frak
K}_{0,\delta}$ and $\bar{\Bbb Q}^1 \le_{{\frak K}'_{0,\delta}}
\bar{\Bbb Q}^2$
\sn
\item "{$(c)$}"  $\Bbb P^1_\delta = \text{ Lim}_{\Cal F}(\bar{\Bbb
Q}^1_\delta) \lessdot \Bbb P^2_\delta = \text{ Lim}_{\Cal F}(\bar{\Bbb
Q}^2_\delta)$
\sn
\item "{$(d)$}"  $\alpha < \delta \Rightarrow \bar{\Bbb Q}^1
\restriction \alpha \le_{{\frak K}'_{\ell,\alpha}} \bar{\Bbb Q}^2
\restriction \alpha$.
\ermn
\ub{Then} $\bar{\Bbb Q}^1 <_{{\frak K}'_\ell} \bar{\Bbb Q}^2$. \nl
4) If $\bar{\Bbb Q}^1 \le_{{\frak K}'_0} \bar{\Bbb Q}^2$ and $\ell g(\bar{\Bbb
Q}^1) = \alpha + 1 < \ell g(\bar{\Bbb Q}^2)$ \ub{then} Lim$_{\Cal
F}(\bar{\Bbb Q}^1) \lessdot \text{ Lim}_{\Cal F}(\bar{\Bbb Q}^2 \restriction \alpha)$.
\endproclaim
\bigskip

\demo{Proof} 1) Easy. \nl
2) Easy. \nl
3) By \scite{c.3p}. \nl
4) By \scite{bs.6}(2).  \hfill$\square_{\scite{6.2abt8}}$\margincite{6.2abt8}
\enddemo
\bn
\ub{\stag{bt.9} Observation}:  1)  If $\bar{\Bbb Q} \in {\frak K}''_{1,\beta}$
and $\lambda \le \beta$, \ub{then} in $V^{{\Bbb P}_\beta}$ 
we have:  if ${\Cal B}^*$ is a Borel
subset of Lebesque measure $1,{\Cal B}^* = \dbcu_n \text{
lim}(T_n),\{T_n:n < \omega\} \subseteq \bold T_g,g \in {\Cal G}^{\bold V}$ 
\ub{then} for a club of $j \in S$ we have
cf$(j) > \aleph_0 \Rightarrow {\underset\tilde {}\to \nu_j} \in {\Cal B}^*$. \nl
2) If $\bar{\Bbb Q} \in {\frak K}'$ and $\beta \le \ell g(\bar{\Bbb Q})$, 
\ub{then} the ${\Cal F}$-forcing notion ${\Bbb P}_\beta$ is essentially 
${\Cal F}$-psc and has 2-pure decidability over ${\Bbb Q}_0$. \nl
3) Moreover, (in (2)) $\bar{\Bbb Q}$ is semi-simple, hence if 
$p \in \Bbb P_\beta,p \Vdash ``\underset\tilde {}\to \tau \in
{}^\omega \text{Ord}"$ \ub{then} 
for some $q$ we have $p \le_{\text{pr}} q$ and
for each $n,{\Cal I}_n = \{r:q \le_{\text{apr}} r$ or just
val$^{{\Bbb P}_\beta}(q) \le_{\text{apr}}$ val$^{{\Bbb P}_\beta}(r)$ and $r$
forces a value to $\underset\tilde {}\to \tau(n)\}$ is predense above $q$ or
even is $\{q':\text{val}(q) \le \text{val}(q')\}$. 
\bigskip

\demo{Proof}  1) Apply clause (e) of 
Definition \scite{bt.8} to ${\Cal B},\langle {\Cal B}_n:n < \omega \rangle$. \nl
2), 3)  By previous theorems.   \hfill$\square_{\scite{bt.9}}$\margincite{bt.9}
\enddemo
\bigskip

\proclaim{\stag{bt.11} Claim}  0) ${\frak K}'_0 \ne \emptyset$. \nl
1) If $\alpha_* \le \alpha < \lambda,\bar{\Bbb Q} \in {\frak K}'_{\alpha_*},
j \in \lambda \backslash A(\bar{\Bbb Q})$, \ub{then} there is
$\bar{\Bbb Q}' \in {\frak K}'_\alpha$ such that 
$\bar{\Bbb Q} \le_{{\frak K}'} \bar{\Bbb Q}' \in {\frak K}^+_\gamma$
and $A(\bar{\Bbb Q}') \supseteq A(\bar{\Bbb Q}) \cup \{j\}$. \nl
2) If $\bar{\Bbb Q}^\zeta \in {\frak K}'_{\le \alpha}$ for $\zeta < \delta$ and
$\varepsilon < \zeta < \delta \Rightarrow \bar{\Bbb Q}^\varepsilon
\le_{{\frak K}'_\alpha} \bar{\Bbb Q}^\zeta$ and (cf$(\delta) >
\aleph_0 \vee \alpha < \lambda$), \ub{then} for some
$\bar{\Bbb Q}^\delta \in {\frak K}'_\alpha$ we have $\zeta < \delta
\Rightarrow \bar{\Bbb Q}^\zeta \le_{{\frak K}'_\alpha} \bar{\Bbb Q}^\delta$.
\nl
3) For $\alpha < \lambda$ and $A \in [\lambda]^{< \lambda}$ there is
$\bar{\Bbb Q} \in {\frak K}'_\alpha$ with $A(\bar Q) \supseteq A$. \nl
4) Part (2) holds for ${\frak K}''_\alpha$, too, if cf$(\delta) > \aleph_0$.
\endproclaim
\bigskip

\remark{Remark}  Note that we have to prove that 
``${\underset\tilde {}\to \nu_\gamma}$ is 
${\Cal G}^{\bold V}$-continuous" is preserved.
\endremark
\bigskip

\demo{Proof}  We prove by induction on $\alpha$ (all parts simultaneously). \nl
0)  Trivial \nl
1) We choose $\bar{\Bbb Q}^\beta$ for $\beta \in [1,\alpha]$
by induction on $\beta$ such that:
\mr
\widestnumber\item{$(iii)$}
\item "{$(i)$}"  $\bar{\Bbb Q}^\beta \in {\frak K}'_\beta$
\sn
\item "{$(ii)$}"  $\bar{\Bbb Q} \restriction \beta \le_{\frak K}
\bar{\Bbb Q}^\beta$
\sn
\item "{$(iii)$}"  $\gamma < \beta \Rightarrow \bar{\Bbb Q}^\gamma =
\bar{\Bbb Q}^\beta \restriction \alpha$
\sn
\item "{$(iv)$}"  $\bar{\Bbb Q}^\beta \in {\frak K}^+_\beta$.  
\endroster
\enddemo
\bn
\ub{Case 1}:  $\beta = 1$.

We choose $\Bbb Q^\beta_0 = \text{ Random}_{A(\bar{\Bbb Q}) \cup \{j\}}$.
\bn
\ub{Case 2}:  $\beta = \beta_* +1 > 1$.

By \scite{bs.6}(1) + (2) we can choose $\bar{\Bbb Q}'_\beta \in
{\frak K}_\beta$ such that $\bar{\Bbb Q}'_\beta \restriction \beta_* =
\bar{\Bbb Q}_{\beta_*}$ and $\bar{\Bbb Q} \restriction \beta
\le_{{\frak K}_\beta} \bar{\Bbb Q}'$.  But maybe $\bar{\Bbb Q}'_\beta
\notin {\frak K}^+_\beta$.  So we try to choose by induction on $i <
\lambda,\bar{\Bbb Q}'_{\beta,i}$ such that:
\mr
\widestnumber\item{$(iii)$}
\item "{$(i)$}"  $\bar{\Bbb Q}'_{\beta,i} \restriction \beta_* \in
{\frak K}^+_{\beta_*}$ and $\bar{\Bbb Q}'_{\beta,0} = \bar{\Bbb Q}'_\beta$
\sn
\item "{$(ii)$}"  $\bar{\Bbb Q}'_{\beta,i}$ is $\le_{{\frak
K}_\beta}$-increasing
\sn
\item "{$(iii)$}"  for each $i,\bar{\Bbb Q}'_{\beta,i+1}$ exemplifying
$\bar{\Bbb Q}'_{\beta,i} \notin {\frak K}^+_\beta$.
\ermn
For $i=0$ no problem, for $i$ limit use the part (2) for $\beta_*$ by
the induction hypothesis.  For $i$ successor if we cannot continue, we
have succeeded.  Having carried the induction by the L.S.-argument (as
$(\forall \alpha < \lambda)(|\alpha|^{\aleph_0} < \lambda = \text{
cf}(\lambda))$ we can find $\delta < \lambda$, cf$(\delta) > \aleph_0$
such that every ``strong isomorphic type" occurs, contradicting (iii).
\bn
\ub{Case 3}:  $\beta$ limit. 

Mostly also easy and similar to part (2) but we need to preserve the
${\Cal G}^{\bold V}$-continuous.
\nl
2) We choose $\bar{\Bbb Q}^{\delta,\beta}$ by induction on 
$\beta \in [1,\alpha]$ such that:
\mr
\widestnumber\item{$(iii)$}
\item "{$(i)$}"  $\bar{\Bbb Q}^{\delta,\beta} \in {\frak K}'_\beta$
\sn
\item "{$(ii)$}"  $\gamma < \beta \Rightarrow \bar{\Bbb Q}^{\delta,\beta}
\restriction \gamma \le_{{\frak K}'_\gamma} \bar{\Bbb Q}^{\delta,\gamma}$
\sn
\item "{$(iii)$}"  $\zeta < \delta \Rightarrow \bar{\Bbb Q}^\zeta
\restriction \beta \le_{{\frak K}'_\alpha} \bar{\Bbb Q}^{\delta,\beta}$.
\endroster
\bn
\ub{Case 1}:  $\beta = 1$.

Let $\Bbb Q^{\delta,\beta}_0 = \text{ Random}_{\cup\{A(\bar Q^\zeta):
\zeta < \delta\}}$.
\bn
\ub{Case 2}:  $\beta = \gamma +1$.

Easy.
\bn
\ub{Case 3}:  $\beta$ limit.  

Take limit. [We should be careful about ${\Cal G}^{\bold
V}$-continuous]. \nl
2), 3) should be clear.   \hfill$\square_{\scite{bt.11}}$\margincite{bt.11}
\bn
\ub{\stag{bt.12} Conclusion}:  1) For any ordinal $\alpha < \lambda$ we can find
$\langle \bar{\Bbb Q}^\zeta:\zeta \le \lambda \rangle$ such that
\mr
\item "{$(a)$}"  $\bar{\Bbb Q}^\zeta \in {\frak K}''_\alpha$
\sn
\item "{$(b)$}"  $\varepsilon < \zeta \Rightarrow \bar{\Bbb Q}^\varepsilon
\le_{{\frak K}'_\alpha} \bar{\Bbb Q}^\zeta$
\sn
\item "{$(c)$}"  $A(\bar{\Bbb Q}^\zeta) \supseteq \zeta$ and $\zeta <
\lambda \Rightarrow |A(\bar{\Bbb Q}^\zeta)| < \lambda$
\sn
\item "{$(d)$}"  if $\zeta < \lambda$ is a limit ordinal of
uncountable cofinality then Lim$(\bar{\Bbb Q}^\zeta) = \dbcu_{\varepsilon <
\zeta}$ Lim$(\bar{\Bbb Q}^\varepsilon)$.
\ermn
2) Let $\Bbb P$ be Lim$(\bar{\Bbb Q}^\lambda)$.  Then
\mr
\item "{$(a)$}"  $\Bbb P$ is a proper forcing notion of cardinality
$\lambda$ satisfying the $\aleph_2$-c.c.c. (so cardinal arithmetic in
$\bold V^{\Bbb P}$ should be clear)
\sn
\item "{$(b)$}"  if cf$(\alpha) > \aleph_0$ then $\Vdash_{\Bbb P}
``{\frak b} = \text{ cf}(\alpha) = {\frak d}"$
\sn
\item "{$(c)$}"  $\Vdash_{\Bbb P} ``\text{there is a set }
\{\nu_\zeta:\zeta < \lambda\} \subseteq {}^\omega 2$ which is not in the
${\Cal G}^{\bold V}$-ideal"
\sn
\item "{$(d)$}"  the continuum in $\bold V^{\Bbb P}$ is $\lambda$.
\ermn
3) For limit $\zeta < \lambda$ of uncountable cofinality, letting
$\Bbb P_\zeta = \text{ Lim}(\bar{\Bbb Q}^\zeta)$, we have
\mr
\item "{$(a)$}"  $\Bbb P_\zeta$ is a proper forcing notion of cardinality
$(|\alpha| + |\zeta|)^{\aleph_0}$ satisfying the $\aleph_2$-c.c.c.
\sn
\item "{$(b)$}"   if cf$(\alpha) > \aleph_0$ then $\Vdash_{\Bbb P}
``{\frak b} = \text{ cf}(\alpha) = {\frak d}"$
\sn
\item "{$(c)$}"  if $A \subseteq \zeta = \sup(A)$ then \nl
$\Vdash_{\Bbb P}$ ``the set $\{\nu_i:i \in A\}$ is not in the ${\Cal
G}^{\bold V}$-ideal
\sn
\item "{$(d)$}"  the continuum in $\bold V^{\Bbb P}$ is $((|\alpha| +
|\zeta|)^{\aleph_0})^{\bold V}$.
\endroster
\bn
\ub{Discussion}:  1) Is there a Cohen in $\bold V^{\Bbb P}$ over $\bold V$?
By the way we construct in general, yes, as possibly $\Bbb P^0_2 \lessdot
\Bbb P$ and $\Bbb P^2_0$ is $\Bbb Q_{{\bar D}^0_1}$ which may add Cohen.
To replace $\langle {\underset\tilde {}\to \nu_i}:i < \lambda \rangle$ by
say a Sierpinski set in $\bold V$ we do not know. \nl
2) Similarly, the Borel conjecture may fail.   
\bigskip

\demo{Proof of \scite{bt.12}}:  Easy by quoting.
\enddemo
\bigskip

\proclaim{\stag{bt.13} Claim}  1) There is $\bar{\Bbb Q} \in {\frak
K}''_{\lambda^+}$. \nl
2) If $\bar{\Bbb Q} \in {\frak K}''_{\alpha},\alpha \ge \lambda$ and
$\Bbb P = \text{ Lim}(\bar{\Bbb Q})$ \ub{then}
\mr
\item "{$(a)$}"  $\Bbb P$ is a proper forcing notion of cardinality
$|\alpha|^{\aleph_0}$ satisfying the $\aleph_2$-c.c.c.
\sn
\item "{$(b)$}"  if cf$(\alpha) \ge \aleph_0$ then $\Vdash_{\Bbb P}
``{\frak b} = \text{ cf}(\alpha) = {\frak d}"$
\sn
\item "{$(c)$}"  $\Vdash_{\Bbb P} ``\{{\underset\tilde {}\to \nu_i}:i
< \lambda\}$ is not in the ${\Cal G}^{\bold V}$-null ideal".
\endroster
\endproclaim
\bigskip

\demo{Proof}  1) Using $\diamondsuit_\lambda$, as in
\cite[IV]{Sh:f}. \nl
2) Like \scite{bt.12}.  \hfill$\square_{\scite{bt.13}}$\margincite{bt.13}
\enddemo
\newpage

\head {\S7 Continuing \cite{Sh:592}} \endhead  \resetall \sectno=7
\bn
\ub{\stag{c.s8.0} Context}:  As in \S6.

At present we can deal with this for ``${\Cal G}^{{\bold V}_0}$-continuous"
instead Random.  To do it fully we need to make the ultrafilter 
${\underset\tilde {}\to D_{\alpha,\eta}}$ Ramsey but we do not know to
guarantee this. 
\proclaim{\stag{c7.2} Theorem}  Assume
\roster
\item "{$(*)(i)$}"  $\kappa \le \theta < \mu < \lambda = \lambda^{< \mu}
= 2^\kappa$
\sn
\item "{$(ii)$}"  $\kappa$ regular and $(\forall \alpha < \kappa)(|\alpha|
^{\aleph_0} < \kappa)$
\sn
\item "{$(iii)$}"  $\theta = \text{ cf}(\theta)$ and $(\forall \alpha < \mu)
(|\alpha|^{\aleph_0} < \mu)$
\sn
\item "{$(iv)$}"  $\mu$ is a limit cardinal
\sn
\item "{$(v)$}"  ${\Cal G} = {\Cal G}^{\bold V}$.
\ermn
\ub{Then} for some forcing notion ${\Bbb P}$ we have:
\mr
\item "{$(\alpha)$}"  ${\Bbb P}$ is an $\aleph_2$-c.c. proper 
forcing notion of cardinality $\lambda$
\sn
\item "{$(\beta)$}"  in $\bold V^{\Bbb P}$ we have 
cov(${\Cal G}$-continuous ideal = Null$_{\Cal G}$) = $\mu$
\sn
\item "{$(\gamma)$}"  in $\bold V^{\Bbb P}$ we have ${\frak b} = {\frak d}= \theta$.
\endroster
\endproclaim
\bigskip

\remark{\stag{c7.3} Remark}  1) We rely on \cite{Sh:592}, if instead we rely
on \cite{Sh:619}, then we can weaken the assumptions on the cardinals. \nl
2) By observation in $\bold V^{\Bbb P}$ we have: 
the covering by closed null sets number is also $\mu$ 
so can be $\aleph_\omega$, i.e. ``${\frak d} < \mu$".
\endremark
\bigskip

\demo{Proof of \scite{c7.2}}  Let ${\frak K}(0)$ be the family of 
$\bar{\Bbb Q} \in {\frak K}^3$ from Definition 2.11 of 
\cite{Sh:592} of length $< \lambda$ ordered naturally: $\bar{\Bbb Q}' \le 
\bar{\Bbb Q}''$ iff $\bar{\Bbb Q}' = \bar{\Bbb Q}'' \restriction \ell g
(\bar{\Bbb Q}')$, of course, $\bar{\Bbb Q}$ stands
for the forcing Lim$(\bar{\Bbb Q})$, \scite{c.1}(d).  
Clearly $\bar{\Bbb Q}$ is FS iteration, this fixes the choice in \scite{c.1}.
Sometimes we replace the ordinal $< \lambda$ by a set of ordinals, with
obvious meaning.  To avoid confusion we use $\bar R$ for such iterations 
and if ${\Bbb Q}_0$ is such a forcing, i.e. Lim$(\bar R)$ we let $\bar R = 
\bar R_{\Bbb Q}$ and
$\langle {\underset\tilde {}\to \eta_\zeta}:\zeta < \ell g(\bar R) \rangle$
for the generic in \scite{c.2} - \scite{c.7}.  Clearly ${\frak K}(0)$ is
cf$(\lambda)$-closed, but as $\lambda = \lambda^{< \mu}$ necessarily
cf$(\lambda) > \mu > \theta$.  
For $\bar{\Bbb Q} \in {\frak K}_\alpha$ let 
$\langle {\underset\tilde {}\to \nu_i}[\bar{\Bbb Q}]:0 < i < \alpha
\rangle$ denotes the sequence of generic reals,
${\underset\tilde {}\to \nu_i}[\bar{\Bbb Q}]$ for
${\underset\tilde {}\to {\Bbb Q}_i}$.  So we can find
$\bar{\Bbb Q}_0 \in {\frak K}'_\theta$.  Now by induction on 
$\zeta < \lambda$ we define $\bar{\Bbb Q}_\zeta$ such that:
\mr
\item "{$(a)$}"  $\bar{\Bbb Q}_\zeta \in {\frak K}^+_\theta \cap
{\frak R}'_\theta$
\sn
\item "{$(b)$}"  $\varepsilon < \zeta \Rightarrow \bar{\Bbb Q}_\varepsilon 
\le_{{\frak K}'} \bar{\Bbb Q}_\zeta$ (hence Lim$_{\Cal F}
(\bar{\Bbb Q}_\varepsilon) \lessdot
\text{ Lim}_{\Cal F}(\bar{\Bbb Q}_\zeta))$
\sn
\item "{$(c)$}"  letting $\bar R_{\bar{\Bbb Q}_\zeta} \in {\frak K}(0)$ be such that
$\Bbb Q_{\zeta,0} = \text{ Lim}(\bar R_{\bar{\Bbb Q}_\zeta})$ 
we have $\ell g(\bar R_\zeta) = \ell g(\bar R_0) + \xi_\zeta$
\sn
\item "{$(d)$}"  \ub{if}
{\roster
\itemitem{ $(i)$ }  $\varepsilon < \lambda,\theta \le \chi = \chi^{\aleph_0}
< \mu$
\sn
\itemitem{ $(ii)$ }  $\bar{\Bbb Q}' \in {\frak K}^+_\theta,
\bar{\Bbb Q}' \le_{{\frak K}'} \bar{\Bbb Q}_\varepsilon$
\sn
\itemitem{ $(iii)$ }  $A \subseteq \ell g(\bar R_{Q_{\varepsilon,0}})$ is
of cardinality $\le \chi$
\endroster}
\ub{then} for some $\zeta \in [\varepsilon,\lambda)$ we have
$\Vdash_{\text{Lim}_{\Cal F}(\bar{\Bbb Q}_{\zeta +1})}
``{\underset\tilde {}\to \eta_{\ell g(\bar R_{\bar{\Bbb Q}_{\zeta,0}})}}$, the partial
random real of the $\ell g(\bar R_{\bar{\Bbb Q}_{\zeta,0}})$-iterand
of the iteration $\Bbb Q_{\zeta,0} \in {\frak K}(0)$ is ${\Cal
G}^{\bold V}$-continuous over
$\bold V[\langle {\underset\tilde {}\to \eta_\zeta}:\zeta \in A \rangle \char 94
\langle {\underset\tilde {}\to \nu_i}:i < \theta \rangle]"$. 
\ermn
There is no problem for $\zeta = 0$ and $\zeta$ limit.  For $\zeta =
\varepsilon +1$, let $\gamma = \ell g(\bar R_{{\Bbb Q}_{\varepsilon,0}}),
\bar R = \bar R_{\bar{\Bbb Q}_{\varepsilon,0}}$.
By bookkeeping we are given $\xi \le \varepsilon$ and $A_\varepsilon \subseteq
\ell g(\bar R_{{\Bbb Q}_{\xi,0}})$ of cardinality $\le \chi$.  By 
the Lowenheim-Skolem argument, choose 
$A^*_\varepsilon \subseteq \gamma$ of cardinality $\le \chi$ including
$A_\varepsilon$, closed enough (in particular, as required in \cite{Sh:592},
see 2.16 (1),(2) and $\bar R_{\bar{\Bbb Q}_{\varepsilon,0}} 
\restriction A^*_\varepsilon \in {\frak K}(0))$ and there is 
$\bar{\Bbb Q}_\varepsilon \le_{{\frak K}'} \bar{\Bbb Q}_\varepsilon$ such that 
$\bar R_{{\Bbb Q}_{\varepsilon,0}} =
\bar R_{\bar{\Bbb Q}_{\varepsilon,0}} \restriction A^*_\varepsilon$.  By
the bookkeeping we can ensure every $A$ will appear.  Let
${\underset\tilde {}\to R_\gamma} = 
\text{ Random}^{\bold V[{\underset\tilde {}\to \eta_\beta}:\beta \in 
A^*_\varepsilon]}$ and let $\bar R'_\varepsilon$ be
$\bar R_{\bar{\Bbb Q}'_{\varepsilon,0}}$ when we add
${\underset\tilde {}\to R_\gamma}$, i.e.
$\bar R_{\bar{\Bbb Q}'_{\varepsilon,0}} \le_{\frak K} \ell g
(\bar R''_\varepsilon) = \ell g(\bar R_{\bar{\Bbb Q}'_{\varepsilon,0}}) \cup
\{\gamma\}$ (abusing notion) $(R''_\varepsilon)_\gamma = R_\gamma$
($\bar R''_\varepsilon$ exists by \cite{Sh:592}).  
By \scite{c.7} we can find $\bar{\Bbb Q}'' \in {\frak K}^+ \cap
K'_\theta$ and $R''_\varepsilon <_{{\frak K}(0)} \Bbb Q''_0$ which 
satisfies $\Vdash_{\text{Lim}_{\Cal F}(\bar{\Bbb Q}'')}
``{\underset\tilde {}\to \eta_\gamma}$ is random over
$\bold V^{\text{Lim}(\bar R \restriction A)}"$.
\nl
By renaming, \wilog \, $A(\Bbb Q''_0) \cap A(\bar{\Bbb Q}^\varepsilon)
= B^*_\varepsilon$.
\mn
Now we can define $\bar{\Bbb Q}_\zeta$ by amalgamation,
i.e. \scite{c7.4} below.  Let
$\bar{\Bbb Q}_\lambda$ by $\dbcu_{\zeta < \lambda} \bar{\Bbb Q}_\zeta$ and
$\Bbb P = \text{ Lim}_{\Cal F}(\bar{\Bbb Q}_\lambda) = \dbcu_{\zeta < \lambda}
\text{ Lim}_{\Cal F}(\bar{\Bbb Q}_\zeta)$.  
It is as required: $\Vdash_{\Bbb P} ``{\frak b}
= {\frak d} = \theta"$ easily, and $\Vdash_{\Bbb P}$ ``cov(null) $\ge \mu$"
by clause (d) and the bookkeeping concerning the $A_\varepsilon$'s.
Lastly $\Vdash_{\Bbb P}$ ``cov(null) $\le \mu$" because $\Vdash_{{\Bbb Q}_{\lambda,0}}$
``cov(null) $\le \mu$" by \cite{Sh:592} and properties of
Lim$_{\Cal F}(\bar{\Bbb Q}_{\lambda,\theta})/{\Bbb Q}_{\lambda,0}$.
\hfill$\square_{\scite{c7.2}}$\margincite{c7.2}
\enddemo
\bigskip

\proclaim{\stag{c7.4} Claim}  Assume
\mr
\item "{$(a)$}"  ${\Bbb Q}_{\ell,0} \in {\frak K}(0)$ for $\ell = 0,1,2,3$ and
${\Bbb Q}_{0,0} \lessdot {\Bbb Q}_{\ell,0} \lessdot {\Bbb Q}_{3,0}$ moreover
${\Bbb Q}_{3,0} = {\Bbb Q}_{1,0} \ast_{{\Bbb Q}_{0,0}} {\Bbb Q}_{2,1}$
\sn
\item "{$(b)$}"  $\bar{\Bbb Q}_\ell \in {\frak K}^+_\alpha$ for $\ell=0,1,2$
\sn
\item "{$(c)$}"  $\bar{\Bbb Q}_0 \le_{{\frak K}'} \bar{\Bbb Q}_1$ and $\bar{\Bbb Q}_0 
\le_{{\frak K}'} \bar{\Bbb Q}_2$.
\ermn
\ub{Then} we can find $\bar{\Bbb Q}_3 \in {\frak K}^+_\alpha$ such that
$\bar{\Bbb Q}_\ell \le_{{\frak K}'} \bar{\Bbb Q}_3$ for $\ell < 3$.
\endproclaim
\bigskip

\remark{Remark}  1) How do we get such $\Bbb Q_{\ell,0} \in {\frak
K}(0)$?  By \cite[Lemma 2.16]{Sh:592}. \nl
2) We can replace $\Bbb Q_{3,0}$ by $\bar Q_{3,\alpha}$ as the proof.
\endremark
\bigskip

\demo{Proof}  We choose $\bar{\Bbb Q}_3 \restriction \beta$ by induction on 
$\beta \in [1,\alpha]$, for $\beta = 1$ there is nothing to do.  For
$\beta$ limit just use $\bar{\Bbb Q}_\ell \restriction \beta \in {\frak
K}^+_\alpha$.
For $\beta = \gamma +1$ use \scite{c.11} below;
we could have demanded something on how $\bar{\Bbb Q}_0 \le_{{\frak
K}'} \bar{\Bbb Q}_2$ (i.e.
choosing $\bar A_\varepsilon$ in the proof of \scite{c7.2} but not needed).
\enddemo
\bigskip

\proclaim{\stag{c.11} Claim}  Assume
\mr
\item "{$(a)$}"  $Q_\ell$ is a forcing notion for $\ell \le 3$
\sn
\item "{$(b)$}"  $Q_0 \lessdot Q_\ell \lessdot Q_3$
\sn
\item "{$(c)$}"  $Q_3 = \Bbb Q_1 \ast_{{\Bbb Q}_0} {\Bbb Q}_1$
\sn
\item "{$(d)$}"  for $\ell = 0,1,2$ we have ${\underset\tilde {}\to D_\ell}$
is a $Q_\ell$-name of an ultrafilter on $\omega$
\sn
\item "{$(e)$}"  for $\ell =1,2$ we have $\Vdash_{Q_\ell} ``
{\underset\tilde {}\to D_0} \subseteq {\underset\tilde {}\to
D_\ell}"$. 
\ermn
\ub{Then} we can find a $\Bbb Q_3$-name $\underset\tilde {}\to D$ such that
$\Vdash_{{\Bbb Q}_3} ``\underset\tilde {}\to D$ is an ultrafilter on $\omega$
extending ${\underset\tilde {}\to D_1} \cup {\underset\tilde {}\to D_2}$.
\endproclaim
\bigskip

\demo{Proof}  As in \cite[\S3]{Sh:326}.  \hfill$\square_{\scite{c.11},\scite{c.10}}$
\enddemo
\bigskip

\proclaim{\stag{c.12} Claim}  Let $\bar{\Bbb Q} \in {\frak K}_\alpha$ and for
$\beta \in (0,\alpha)$ let ${\underset\tilde {}\to \nu_\beta}$ be the
generic real of ${\Bbb Q}_\alpha$.  \ub{Then} 
${\underset\tilde {}\to G_{\text{Lim}
_{\Cal F}(\bar{\Bbb Q})}}$ can be computed from $\langle
{\underset\tilde {}\to G_{{\Bbb Q}_0}} \rangle \char  94 \langle
{\underset\tilde {}\to \nu_\beta}:\beta \in (0,\alpha) \rangle$.
\endproclaim
\bigskip

\demo{Proof}  As usual (or see \cite{Sh:592}).
\enddemo
\newpage

\head {\S8 On $\eta$ is ${\Cal L}$-big over $M$} \endhead  \resetall \sectno=8
\bigskip

\definition{\stag{p.1} Definition}  1) Let $\bold T = \{{\Cal T}:{\Cal
T} \subseteq {}^{\omega >} {\Cal H}(\aleph_0),{\Cal T} \ne
\emptyset,{\Cal T}$ closed under initial segments, no $\triangleleft$-maximal
member and ${\Cal T}_n = \{\eta \in {\Cal T}:\ell g(\eta) = n\}$
finite for $n < \omega\}$. \nl
2) For ${\Cal T}_1,{\Cal T}_2 \in \bold T$ let $\bold R_{{\Cal
T}_1,{\Cal T}_2} = \{R:R$ a closed subset of lim$({\Cal T}_1) \times
\text{ lim}({\Cal T}_2)\}$.  Similarly for $R_{\Cal T}$.

We write $\eta R_n \nu$ instead of
$(\eta,\nu) \in R_n$.  We always assume that ${\Cal T}_1,{\Cal T}_2$ can be
reconstructed from $\bar R \in \bold R_{{\Cal T}_1,{\Cal T}_2}$ and
write ${\Cal T}_1[R],{\Cal T}_2[R]$; 
similarly $R \in \bold R_{\Cal T}$.  Let $\bold R_* = \cup\{\bold
R_{{\Cal T}_1,{\Cal T}_2}:{\Cal T}_1,{\Cal T}_2 \in \bold T\}$. \nl
3) If $R$ is a closed subset of lim$({\Cal T}_1) \times \text{
lim}({\Cal T}_2)$ and $k < \omega$ then let $R^{<k>} = \{(\eta
\restriction k,\nu \restriction k):(\nu,\eta) \in R\}$.  Similarly for
$R \subseteq \text{ lim}({\Cal T})$. \nl
4) For every ${\Cal Y} \subseteq \bold Y =: \{(f,{\Cal
T}):{\Cal T} \in\bold T,f \in \dsize \prod_{n < \omega} {\Cal P}({\Cal
T}_n)\}$ let $D_{\Cal Y} = \{A \subseteq \omega$: for some $k,m < \omega$ and
$(f_\ell,{\Cal T}_\ell) \in {\Cal Y}$ and $\nu_\ell \in \text{
lim}(T_\ell)$ for $\ell < k$ we have $A \supseteq \{n:n \ge m$ and
$\ell < k \Rightarrow \nu_\ell(n) \in f_\ell(n)\}$.  We say ${\Cal Y}$ is
nontrivial if $\emptyset \notin D_{\Cal Y}$.  Let $J_{\Cal Y}$ be the dual
ideal.
\enddefinition
\bigskip

\definition{\stag{p.2} Definition}  1) We say $D$ is $(f,{\Cal T})$-narrow \ub{if}:
\mr
\widestnumber\item{$(iii)$}
\item "{$(i)$}" $f \in \dsize \prod_{n < \omega} {\Cal P}({\Cal T}_n)$
\sn
\item "{$(ii)$}"  $D$ is a filter on $\omega$ containing the co-finite
subsets
\sn
\item "{$(iii)$}"  for every $\nu \in {}^\omega 2$ the set $\{k <
\omega:\nu \restriction n \in f(n)\}$ belongs to $D$.
\ermn
2)  For ${\Cal Y} \subseteq \bold Y$, we say 
$D$ is ${\Cal Y}$-narrow if $D$ is $(f,{\Cal T})$-narrow 
for every $(f,{\Cal T}) \in {\Cal Y}$. \nl
3)  $\bold Z = \{(\eta,R):\eta \in \text{ lim}({\Cal T}_2[R]),
R \in \bold R_*\},\bold Z_M = \{(\eta,R) \in \bold Z:R
\in M$ and $\eta \in \text{ lim}({\Cal T}_2[R])\}$. \nl
4)  We say that $D$ is $(\eta,R)$-\ub{big} over $M$ if: it is
$\{(\eta,R)\}$-big over $M$ (see below). \nl
5) We say that $(\eta,D)$ is ${\Cal L}$-big over $M$ if
\mr
\widestnumber\item{$(iii)$}
\item "{$(i)$}"  $M$ is a set or a class (usually an inner model), $D$
is a filter on $\omega$ containing the co-bounded subsets of $\omega$
\sn 
\item "{$(ii)$}"  ${\Cal L} \subseteq \bold Z_M,R \in \bold R^M_{{\Cal
T}_1{\Cal T}_2},{\Cal T}_1,{\Cal T}_2,{\Cal T} \in \bold T^M$
\sn
\item "{$(iii)$}"  $\eta \in \text{ lim}({\Cal T}_2[R])$ when ${\Cal
L} = \{(\eta,R)\}$ for \footnote{note that if $D$ is an ultrafilter
then the case $m^*=1$ suffices}
every $m^* < \omega,\langle \nu_{m,n}:n < \omega \rangle \in M,\nu_m
\in M$ for $m < m^*,(\eta_m,R_m) \in {\Cal L}$ such that
$\{\nu_{m,n},\nu_m:n < \omega\} \subseteq \text{ lim}({\Cal
T}_1[R_m])$, if $m < m^* \Rightarrow \nu_m = \text{ lim}_D \langle
\nu_{m,n}:n < \omega \rangle$ \ub{then} $\{n:m < m^* \Rightarrow \nu_m R_m
\eta_m \equiv \nu_{m,n} R_m \eta_m\} \ne \emptyset \text{ mod } \in D$.
\ermn
6) We say $\eta$ is $R$-big over $M$ if the filter of co-finite
subsets of $\omega$ is $(\eta,R)$-big.  We say that $\eta$ is ${\Cal
L}$-big over $M$ if $(\eta,D)$ is with $D$ the filter of co-finite
subsets of $\omega$.
\enddefinition
\bigskip 

\proclaim{\stag{p.3} One Step Claim}  Assume (all in $\bold V_2$)
\mr
\item "{$(a)$}"  $\bold V_1 \subseteq \bold V_2$
\sn
\item "{$(b)$}"  in $\bold V_1,D_1$ is a non principal ultrafilter on $\omega$
\sn
\item "{$(c)$}"  for $\zeta < \zeta^*_1,f_\zeta \in \bold V_1,\bold
T_\zeta \in \bold T^{{\bold V}_1}$
\sn
\item "{$(d)$}"  $D_1$ is ${\Cal Y}$-narrow, ${\Cal Y} \subseteq \bold
Y^{{\bold V}_1}$ (of course, ${\Cal Y} \in \bold V_2$, but possibly
${\Cal Y} \notin \bold V_1$)
\sn
\item "{$(e)$}"  $Z \subseteq \bold Z_{\bold V_1}$
\sn
\item "{$(f)$}"  $({\Cal Y},{\Cal Z})$ is high over $\bold V_1$ which
means (in $\bold V_1$): \nl
${\Cal Y} \subseteq \bold Y^{{\bold V}_1},{\Cal Z} 
\subseteq \bold Z^{{\bold V}_2}_{{\bold V}_1}$, and $m^* <
\omega,(\eta_m,R_m) \in {\Cal Z}$ for $m < m^*,{\Cal Y}'
\subseteq {\Cal Y}$ is finite, $B \in (J_{{\Cal Y}'})^+$
and for $m < n^*,\bar \nu = \langle \nu_{m,n}:n \in B \rangle \in \bold
V_1$ and $\nu_{m,n},\nu_m \in {\text{\rm lim\/}}({\Cal T}^1[R_\ell])^{{\bold
V}_1},\nu_m = {\text{\rm lim\/}} \langle \nu_{m,n}:n \in B \rangle$ \ub{then}
$\{n \in B$: if $m < m^*$ then $\nu_{m,n} R_n \eta_m \equiv \nu R_m
\eta_m\}$ is infinite.
\ermn
\ub{Then} there is $D_2$ such that:
\mr
\item "{$(\alpha)$}"  $D_2 \, (\in \bold V_2)$  is an ultrafitler on $\omega$
\sn
\item "{$(\beta)$}"  $D_1 \subseteq D_2$
\sn
\item "{$(\gamma)$}"  $D_2$ is ${\Cal Y}$-narrow over $\bold V_1$
\sn
\item "{$(\delta)$}"  $D_2$ is ${\Cal Z}$-big over $\bold V_1$.
\endroster
\endproclaim
\bigskip

\remark{Remark}  Better if we predetermine $D_2 \cap {\Cal
P}(\omega)^{\bold V}$, good for ${\frak u} = \text{ cf}(\alpha^{\frak
q}) > \aleph_0$.
\endremark
\bigskip

\demo{Proof of \scite{p.3}}

For $(f,{\Cal T}) \in {\Cal Y}$ and $\rho \in \text{ lim}({\Cal
T})^{{\bold V}_2}$ let

$$
A^1_{f,\rho} = \{n:\rho(n) \in f(n)\}.
$$
\mn
For $(\eta,R) \in {\Cal Z}$ and $\nu_n,\nu \in \text{ lim}({\Cal
T}[R]) = \bold V_1$ for $n < \omega$ such that $\bar \nu = \langle \nu_n:n < \omega
\rangle \in \bold V_1$ and $\nu = \text{ lim}_D(\bar \nu)$ we let

$$
A^2_{\eta,R,\bar \nu,\nu} = \{n:\nu_n R \eta \equiv \nu R \eta\}.
$$
\mn
So we just need to find an ultrafilter $D_2$ on $\omega$ which extends
$D_1 \cup \{A^2_{f,\nu}:\nu \in \text{lim}({\Cal T}),(f,{\Cal T}) \in
{\Cal Y}_1\} \cup \{A^2_{\eta,R,\bar \nu,\nu}:(\eta,R) \in {\Cal Z}$
and $\bar \nu,\nu \in \bold V_1$ are as above$\}$.  For this it suffices to
prove
\mr
\item "{$(*)$}"   assume $B \in D_1,n^*_1 < \omega,n^*_2 <
\omega,A^1_{f_\ell,\rho_\ell},
A^2_{\eta_n,R_m,\bar \nu_m,\nu_m}$ well defined for $\ell <
n^*_1,m < n^*_2$ \ub{then} $B \cap \cap \{A^1_{f_\ell,\rho_\ell}:\ell < n^*_1\}
\cap \{A^2_{\eta_m,R_m,\bar \nu_m,\nu_m}:m < n^*_2\}$ is non empty
where $(f_\ell,\rho_\ell) \in {\Cal Y},(\eta_m,R_m) \in {\Cal Z}$ and
$\bar \nu_m,\nu_m$ as usual.
\ermn
As $\nu_m = \text{ lim}_D \langle \nu_{m,n}:
n < \omega \rangle$ and $B \in D_1$ we have $B_{m,k} = \{m \in
B:\nu_{m,n} \restriction k = \nu_m \restriction k\} \in D_1$ for $m
< n^*_2,k < \omega$ hence $B_k = \dbca_{m < n^*_2} B_{m,k} \in D_1$ and
clearly $B_{k+1} \subseteq B_k$, hence $B_k \ne \emptyset$ mod $D$
hence $B_k \notin J^{\bold V_1}_{\Cal Y}$, see Definition
\scite{p.1}(4)
where ${\Cal Y}^* =: \{(f_\ell,{\Cal T}_\ell):\ell < n^*_1\}$.

Clearly it suffice to prove that $B \cap \cap \{A^2_{\eta_m,R_m,\bar
\nu_m,\nu_m}:m < n^*_2\}$ is not in $J_{{\Cal Y}^*}$.  By part (2) of Claim \scite{p.4}
below there is $B^* \subseteq B$ in $\bold V_1$ such that $B^* \backslash B_k$ is
finite for $k < \omega$ and $B^* \notin J_{{\Cal Y}^*}$.

So we have:
\mr
\item "{$(*)(i)$}"  ${\Cal Y}^* \subseteq {\Cal Y}$ is finite
\sn
\item "{$(ii)$}"  $B^* \subseteq \omega$ is from $\bold V_1,B^* \notin J_{{\Cal Y}^*}$
\sn
\item "{$(iii)$}"  $\nu_m,\nu_{m,n} \in 
\text{ lim}({\Cal T}^2[R_m])$ for $m < n^*_2,n < \omega$
\sn
\item "{$(iv)$}"  $\nu_m = \text{ lim} \langle \nu_{m,n}:n \in B^*
\rangle$ for $m < n^*_2$
\sn
\item "{$(v)$}"  $\langle \nu_{m,n}:n < \omega \rangle$ and $\nu_m$
belongs to $\bold V_1$.
\ermn
By assumption (f) we are done.  \hfill$\square_{\scite{p.3}}$\margincite{p.3}
\enddemo
\bigskip

\proclaim{\stag{p.4} Claim}  1) Let ${\Cal Y} \subseteq \bold Y$, then
the following are equivalent for $B \subseteq \omega$:
\mr
\item "{$(i)$}"  $B \notin J_{\Cal Y}$
\sn
\item "{$(ii)$}"  for every $n^* < \omega,(f_\ell,{\Cal T}_\ell) \in
{\Cal Y}$ for $\ell < n^*$ and $m_0 < \omega$ there is $m_1 \in
(m_0,\omega)$ such that:
{\roster
\itemitem{ $(*)$ }  if $\nu_\ell \in ({\Cal T}_\ell)_{m_1}$
for $\ell < n^*$ \ub{then} for some $n \in B \cap [m_0,m_1)$ we have $(\forall
\ell < n^*)(\nu_\ell \restriction n \in f_\ell(n))$. 
\endroster}
\ermn
2) For ${\Cal Y} \subseteq \bold Y$, if $B_n \in J^+_{\Cal Y},B_{n+1}
\subseteq B_n$ \ub{then} there is $B \in J^+_{\Cal Y}$ such that $n <
\omega \Rightarrow B \subseteq^* B_n$. \nl
3) If $\bold V_1 \subseteq \bold V_2,{\Cal Y} \subseteq \bold
Y^{{\bold V}_1},{\Cal Y} \in V_1,A \in {\Cal P}(\omega)^{{\bold
V}_1}$, \ub{then} $A \in J^{{\bold V}_1}_{\Cal Y} \Leftrightarrow A
\in J^{{\bold V}_2}_{\Cal Y}$.
\endproclaim
\bigskip

\demo{Proof}  Easy. \nl
1) Assume clause (i), i.e. $B \notin J_{\Cal Y}$; to prove clause (ii)
assume toward contradiction that $n^*$ and $\langle (f_\ell,{\Cal
T}_\ell):\ell < n^* \rangle$ and $m_\ell < \omega$ are as there but
there is no $m_1 \in (m_0,\omega)$ such that $(*)$ there holds, so
there are $\nu^{m_1}_\ell \in ({\Cal T}_\ell)_{m_1}$ for $\ell < m^*$
such that $n \in B \cap [m_0,m_1) \Rightarrow (\exists \ell <
n^*)(\nu_\ell \notin f_\ell(n))$.  By Konig lemma there are $\nu_\ell
\in \text{ Lim}({\Cal T}_\ell)$ for $\ell < n^*$ such that $\forall m
< \omega,\exists^\infty m_2 < \omega(m < m_1 \and m_0 < m_1 \and \dsize
\bigwedge_{\ell < n^*} \nu^{m_2}_\ell \restriction m = \nu_\ell
\restriction m)$.

Now for each $\ell$, by Definition \scite{p.1}(4) as $(f_\ell,{\Cal
T}_\ell) \in {\Cal Y}$, the set $A_\ell =: \{m < \omega:\nu_\ell
\restriction m \in f_\ell(m)\} \in D_{\Cal Y}$ hence $A = \dbca_{\ell
< n^*} A_\ell \in D_Y$, but we assume $B \ne \emptyset$ mod $D_{\Cal
Y}$ hence $A \cap B \ne \emptyset$ mod $D_{\Cal Y}$ so there is
$m_\ell,m_0 < m \in A \cap Y$.  Let $m_1 > m$ be such that $\ell < n^*
\Rightarrow \nu^{m_1}_\ell \restriction m = \nu'_\ell \restriction m$
and this $m_1$ contradicts the choice of $\langle \nu^{m_1}_\ell:\ell
< n^* \rangle$.  So $(i) \Rightarrow (ii)$ indeed.  The other
direction is even easier. \nl
2) Just use clause (ii) of part (i) as the definition.  This is
straight. \nl
3) Follows using clause (ii) of part (1).
\enddemo
\bigskip

\proclaim{\stag{p.5} The limit Claim}  Assume:
\mr
\item "{$(a)$}"  $\delta$ a limit ordinal
\sn
\item "{$(b)$}"  $\langle \bold V_\zeta:\zeta < \delta \rangle$ is an
increasing sequence of inner models
\sn
\item "{$(c)$}"  ${\Cal Y}_\zeta \subseteq \bold Y^{{\bold V}_\zeta}$
is increasing with $\zeta$
\sn
\item "{$(d)$}"  $D_\zeta$ is a filter on ${\Cal P}(\omega)^{{\bold
V}_\zeta}$, increasing with $\zeta$
\sn
\item "{$(e)$}"  $D_\zeta$ is disjoint to $J_{\Cal Y}$ for every
finite ${\Cal Y} \subseteq {\Cal Y}_\zeta$.
\ermn
\ub{Then} $\cup \{D_\zeta:\zeta < \delta\}$ can be extended to a
uniform ultrafilter on $\omega$ disjoint to $J_{\cup\{{\Cal
Y}_\zeta:\zeta < \delta\}}$.
\endproclaim
\bigskip

\demo{Proof} Easy.
\enddemo
\bigskip

\definition{\stag{p.6} Definition}  1) Assume ${\Cal T} \in \bold T,h
\in {}^\omega \omega,\infty = \text{ lim}\langle h(n):n < \omega
\rangle$ and $f \in \dsize \prod_{n < \omega} {\Cal P}({\Cal
T}_{h(n)})$.  We say that ``$D$ on $\omega$ is $(f,h,{\Cal
T})$-narrow" if
\mr
\item "{$(a)$}"  $D$ is a filter on $\omega$ containing the co-bounded
subsets
\sn
\item "{$(b)$}"  for every $\nu \in \text{ lim}({\Cal T})$, the set
$\{n:\nu \restriction h(n) \in f(n)\}$ belongs to $D$.
\ermn
2) We say $(f',{\Cal T}')$ is the translation of $(f,h,{\Cal T})$ if:

$$
{\Cal T}' = \{\langle \eta \restriction h(m):m < n \rangle:n <
\omega,\eta \in \text{ lim}({\Cal T})\},
$$

$$
f'(n) = \{\langle \eta \restriction h(m):m \le n \rangle:\eta \in
\text{ lim}({\Cal T}) \text{ and } \eta \restriction h(n) \in f(n)\}.
$$
\enddefinition
\bigskip

\remark{\stag{p.6a} Remark}  We may like to have ${\Cal Y}' \subseteq
\bold Y^{{\bold V}_2}$ is this needed?  Helpful.
\endremark
\bigskip

\definition{\stag{p.7} Definition}  1) Let for a class $M,{\frak Z}_M$
be the set of $(\eta,\bar R)$ such that:
\mr
\item "{$(a)$}"  $\bar R \in {}^\omega(R^M_*),{\Cal T}_1[R_n] = {\Cal
T}_1[R_0]$ and ${\Cal T}_2[R_n] = {\Cal T}_2[R_0]$ and $\bar R$ we let
${\Cal T}_\ell[\bar R] = {\Cal T}_\ell[R_0]$ for $\ell = 1,2$
\sn
\item "{$(b)$}"  $\eta \in \text{ lim}({\Cal T}_2)$
\sn
\item "{$(c)$}"  $\eta$ does $\bar R$-cover $M$, which means $(\forall \nu
\in \text{ lim}({\Cal T}_1[\bar R])^M)(\exists n < \omega)[\nu R_n
\eta]$ [return: content with one $R$].
\endroster
\enddefinition
\bigskip

\proclaim{\stag{p.8} Claim}  1) Assume (with $\bold V = \bold V_2$)
\mr
\item "{$(a)$}"  $\bold V_1 \subseteq \bold V_2 = \bold V$
\sn
\item "{$(b)$}"  ${\Cal L} \subseteq \bold Z_{{\bold V}_1}$
\sn
\item "{$(c)$}"  $\bar D_2 = \langle D_{2,\eta}:\eta \in {}^{\omega >}
\omega \rangle,D_{2,\eta}$ a non principal ultrafilter on $\omega$
(all in $\bold V_2$)
\sn
\item "{$(d)$}"  $\langle D_{1,\eta}:\eta \in {}^{\omega >} \omega
\rangle \in \bold V_1$ where $D_{1,\eta} = D_{2,\eta} \cap {\Cal
P}(\omega)^{{\bold V}_1}$
\sn
\item "{$(e)$}"  $D_{2,\eta}$ is ${\Cal L}$-big ultrafilter over
$\bold V_1$ for every $\eta \in {}^{\omega >} \omega$.
\ermn
\ub{Then} $\Vdash_{{\Bbb Q}_{\bar D_2}} ``\eta$ is $R$-big over
$\bold V_1[{\underset\tilde {}\to \eta_{Q_{\bar D_2}}}]"$. \nl
2) Assume (a), (c), (d) above and
\mr
\item "{$(b)'$}"  $(\eta,\bar R) \in {\frak Z}_{{\bold V}_1}$
\sn
\item "{$(e)'$}"  $D_{2,\eta}$ is $(\eta,R_n)$-big ultrafilter when
$\eta \in {}^{n} \omega$.
\ermn
\ub{Then} $\Vdash_{{\Bbb Q}_{\bar D_2}} ``(\eta,\bar R) \in {\frak
Z}_{{\bold V}_1[{\underset\tilde {}\to \eta}[Q_{\bar D_2}]]}"$.
\endproclaim
\bigskip

\demo{Proof}  [Saharon \ub{revised}: copied from \scite{bs.6}(2).
\nl
1) By \sciteu{p.7a} this is a special case of part (2).
So assume that $p \in \Bbb Q_{{\bar D}_1},m^* < \omega$ and for
each $m < m^*,(\rho_m,R_m) \in {\Cal L}$ and ${\underset\tilde {}\to \nu^m},
\langle {\underset\tilde {}\to \nu^m_n}:n < \omega \rangle \in \bold V$
are $\Bbb Q^{\bold V}_{\bar D}$-names hence $\Bbb Q_{{\bar
D}_1}$-names such that
\mr
\item "{$(*)_1$}"  $p \Vdash_{{\Bbb Q}_{{\bar D}_1}} ``{\underset\tilde {}\to \nu^m},
{\underset\tilde {}\to \nu^m_n} \in \lim({\Cal T}_1[R_m])$
and
${\underset\tilde {}\to \nu^m} = \text{ lim}\langle {\underset\tilde {}\to \nu_m}:
n < \omega \rangle"$.
\ermn
By the definition and what we need to prove, \wilog
\mr
\item "{$(*)_2$}"  $p \Vdash ``{\underset\tilde {}\to \nu^m}
\restriction n = {\underset\tilde {}\to \nu^m_n} \restriction n"$. 
\ermn
We shall find $p' \ge p$ in $\Bbb Q_{{\bar D}_1}$ such that $p'
\Vdash ``(\nu^m_n R_m \rho_m) \equiv (\nu^m R_m \rho_m)$ for every $m < \omega$; 
for some $n < \omega"$, this suffice (see \scite{pr.6}(2)); work in $\bold V$.  Let
$q_0 = ({}^{\omega >}\omega)$, so $q_0 \in \Bbb Q_{\bar D}$, now
we find $\langle \nu^m_\eta,\nu^m_{n,\eta}:\eta \in q_0,n < \omega \rangle$ of
course in $\bold V$ such that:
\mr
\item "{$(*)_3(i)$}"  $\nu^m_\eta,\nu^m_{n,\eta} \in \lim({\Cal T}_1[R_m]),m < m^*$
\sn
\item "{$(ii)$}"  for every $\eta \in q_0$ and $k < \omega$ we can find
$q^m_{\eta,k},q^m_{n,\eta,k} \in \Bbb Q_{\bar D}$ such that: \nl
$q^{[\eta]}_0 \le_{\text{pr}} q^m_{\eta,k},q^{[\eta]}_0 \le_{\text{pr}}
q^m_{n,\eta,k}$ \nl
$q^m_{\eta,k} \Vdash_{{\Bbb Q}_{\bar D}} ``{\underset\tilde {}\to \nu^m} \restriction
k = \nu^m_\eta \restriction k"$
\nl
$q^m_{n,\eta,k} \Vdash_{{\Bbb Q}_{\bar D}} ``{\underset\tilde {}\to
\nu^m_n} \restriction k = \nu^m_{n,\eta} \restriction k"$.
\ermn
Now clearly
\mr
\widestnumber\item{$(*)_4(iii)$}
\item "{$(*)_4$(i)}"   $\nu^m_\eta = \text{ lim}_{D_\eta} \langle 
\nu^m_{\eta \char 94 <k>}:k < \omega \rangle$
\sn
\item "{${}(ii)$}"  $\nu^m_{n,\eta} = \text{ lim}_{D_\eta}\langle
\nu^m_{n,\eta \char 94 <k> \rangle}:k < \omega \rangle$.
\ermn
Next note that
\mr
\item "{$(*)_5$}"   $\nu^m_\eta = \text{ lim}\langle \nu^m_{n,\eta}:
n < \omega \rangle$. \nl
[Why by $(*)_2$ let $u_\eta = \{m < m^*:\nu^m_\eta R_m \rho_m \text{ holds}\}$.
\ermn
Now as each $R_m$ is closed (see Definition \scite{p.1}(2)) there is
$k_\eta < \omega$ such that
\mr
\item "{$(*)_6$}"  if $m < m^*,m \notin u_\eta$ and $\nu^m_\eta
\restriction k \triangleleft \nu \in \text{ lim}({\Cal
T}_1[R_m]),\rho_m \restriction k \triangleleft \rho \in \lim({\Cal
T}_2[R])$ then $\neg(\nu R_m \rho)$.
\ermn
By $(*)_4(i) + (ii) + (*)_6$ we have
\mr
\item "{$(*)_7(i)$}"  if $\neg(\nu^m_\eta R_m \rho_m)$ implies $\{k <
\omega:\nu^m_{\eta \char 94 \langle k \rangle} R_m \rho_m\} \notin
D_\eta$
\sn
\item "{$(ii)$}"  if $\neg(\nu^m_{n,\eta} R_m \rho_m)$ implies $\{k <
\omega:\nu^m_{n,\eta \char 94 \langle k \rangle} R_m \rho_m\} \notin
D_\eta$.
\ermn
By the assumption on $D_{1,\eta}$ and $(*)_4(i) + (ii)$ 
we have
\mr
\item "{$(*)_8(i)$}"  $\nu^m_\eta R_m$  iff $\{k < \omega:\nu^m_{\eta
\char 94 \langle k \rangle} R_m \rho_m\} \in D_\eta$
\sn
\item "{$(ii)$}"  $\nu^m_{n,m} R_m \rho_m$  iff $\{k:\nu^m_{n,\eta
\char 94 \langle k \rangle} R_m \rho_m \text{ holds}\} \in D_\eta$.
\ermn
By $(*)_5$ applied to $\eta = \text{ tr}(p)$, we can find $n(*) < \omega$ such that
$(\forall m < m^*)[(\nu^m_\eta R_m \rho_m) \equiv (\nu^m_{n(*),\eta}
R_m \rho_m)]$.  Next let

$$
p^* =: \{\nu \in p:\text{if } \ell g(r(p)) \le 
\ell \le \ell g(\nu) \text{ and } m < m^* \text{ then }
\nu^m_{\nu \restriction \ell} R_m \rho_m \equiv \nu^m_{n(*),\nu
\restriction \ell} R_m \rho_m\}.
$$
\mn
Now $p \le_{\text{pr}} q \in \Bbb Q_{\bar D}$ by $(*)_8$.  Lastly, let
$q^* =: \{\nu \in p^*:\text{if } \ell < \ell g(\nu)$, \ub{then}
$\nu \in q^m_{\nu \restriction \ell,k_\eta}$ and $\nu \in
q^m_{n(*),\nu \restriction \ell,k_\eta}\}$.

Does $q^* \Vdash_{{\Bbb Q}_{{\bar D}_1}} ``({\underset\tilde {}\to
\nu^m} R_m \rho_m) \equiv ({\underset\tilde {}\to \nu^m_{n(*)}} R_m
\rho_m)"$?  If not, then for some $q^{**}$ we have $q^* \le q^{**}$ and
$q^{**} \Vdash_{{\Bbb Q}_{{\bar D}_1}} ``({\underset\tilde {}\to
\nu^m} R_m \rho_m) \equiv \neg({\underset\tilde {}\to \nu^m_{n(*)}}
R_m \rho_m)"$; moreover, without loss of generality for some truth value $\bold t,
q^{**} \Vdash_{{\Bbb Q}_{\bar D_2}} ``({\underset\tilde {}\to \nu^m}
R_m \rho_m) \equiv \bold t$ and $({\underset\tilde {}\to \nu^n_{n(*)}}
R_m \rho_m) \equiv \neg \bold t"$ and for some $k^* < \omega, q^{**}
\Vdash_{{\Bbb Q}_{\bar D_2}} ``\bold t =$ false $\Rightarrow (\forall
\nu,\rho)
[{\underset\tilde {}\to \nu^m} \restriction k^* \trianglelefteq \nu
\in \lim({\Cal T}_1[R_m]) \and \rho_m \restriction k^* \triangleleft
\rho \in \lim({\Cal T}_2[R_m]) \rightarrow (\nu,\rho) \notin R_m$ and
$\neg \bold t =$ false $\Rightarrow (\forall
\nu,\rho)[{\underset\tilde {}\to \nu^m_{n(*)}} \restriction k^*
\triangleleft \nu \in \lim({\Cal T}_1[R_m]) \and \rho_m \restriction
k^* \triangleleft \rho \in \lim({\Cal T}_2[R_m]) \rightarrow
(\nu,\rho) \notin R_m]$.
But $q^{**},q^m_{\text{tr}(q^{**}),k^*},
q^m_{n(*),\text{tr}(q^{**}),k^*}$ for $m < m^*$ 
are compatible having the same trunk, so let
$q'$ be a common upper bound with tr$(q'') = \text{ tr}(q^{**})$ and we get a
contradiction. \nl
2) Return!  \hfill$\square_{\scite{p.8}}$\margincite{p.8}
\enddemo
\bigskip

\definition{\stag{p.9} Definition}  1) For $g \in {\Cal G}$ let

$$
{\Cal T}[g] = {}^{\omega >} 2
$$

$$
{\Cal T}_2[g] = \{\langle T \cap {}^{g(\ell)} 2:\ell < n\rangle:T \in
\bold T_g,n < \omega \},
$$
\mn
so $\eta \in \text{ lim}({\Cal T}_2[g])$ can be identified with

$$
T = T[\eta] \in \bold T_g:\eta = \eta_T = \left \langle \langle T \cap {}^{\ell
\ge}2:\ell < n \rangle:n < \omega \right \rangle
$$

$$
R_g = \{(\eta,\nu):\nu \in {}^{\omega >}2,\eta \in \text{ lim}({\Cal
T}_2[g]) \text{ and } \nu \in \text{ lim}(T[\eta])\}.
$$
\mn
2) Let $\bar w^* = \langle w^*_k:k < \omega \rangle$ list with no
repetition $\cup\{{\Cal P}({}^n 2) \backslash \{\emptyset\}:n <
\omega\}$ such that $w^*_\ell \in {}^{\text{les}[w^*]}2,\ell_1 < \ell_2
\Rightarrow {}^{\text{les}}[w^*_{\ell_1}] \le [w^*_{\ell_2}]$ and $\ell _1 <
\ell_2 \and \text{ les}[w^*_{\ell_1}] = \text{ les}[w^*_{\ell_2}] \Rightarrow$ the
$<_{\text{lex}}$-first $\rho$ such that $\rho \in w^*_{\ell_1} \equiv
\rho \notin w^*_{\ell_2}$ satisfies $\rho \in w^*_\ell$. Let it be
defined as ${\Cal T}^0 = {}^{\omega >}2$ and $h^*(\ell) =
n[w^*_\ell]$. 
\hfill$\square_{\scite{p.9}}$\margincite{p.9}
\enddefinition
\bn
The following shows that the ``${\Cal G}$-continuous" treated in
\S4-\S6 fit our present framework.
\proclaim{\stag{p.10} Claim}  1) Assume $\bold V_1 \subseteq \bold
V_2,g \in {\Cal G}^{{\bold V}_1},\bold r \in ({}^\omega 2)^{{\bold
V}_2}$, \ub{then} we have: $\bold r$ is $\{g\}$-continuous over $\bold V_1$ iff
$\bold r$ is $R_g$-big over $\bold V_1$. \nl
2) Assume:
\mr
\item "{$(a)$}"  $\bold V_1 \subseteq \bold V_2 = \bold V$
\sn
\item "{$(b)$}"  ${\Cal Z} \subseteq \{(\bold r,R_g):g \in {\Cal G}^{{\bold
V}_1},\bold r \in ({}^\omega 2)^{{\bold V}_1}$ and $\bold r$ is
$R_g$-big over $\bold V_1\}$.
\ermn
\ub{Then} ${\Cal Y} = \emptyset$ and ${\Cal Z}$ are as required in
\scite{p.3}, i.e. $({\Cal Y},{\Cal L})$ is high over $\bold V_1$
(i.e. clause (f) there). \nl
3) Assume (a), (b) as in (2) and
\mr
\item "{$(c)$}"  ${\Cal Y} = \{(\bar w,h^*,{\Cal Y}^0_*)\}$.
\ermn
\ub{Then} $({\Cal Y},{\Cal Z})$ are as required in \scite{p.3}.
\endproclaim
\bigskip

\demo{Proof}   (Or use ${\Cal Z}$ with $R^1$, see
later). \nl
1) Compare the definition \scite{pr.5}(1) + (3) and
\scite{p.2}(3). \nl
2) So assume $m^* < \omega$ and $(\eta_m,R_m) \in {\Cal L}$ for $m <
m^*$ and $B \subseteq \omega$ is infinite and $\nu_m,\nu_{m,n} \in
\lim({\Cal T}_1[R_m]),\langle \nu_{m,n}:n < \omega \rangle \in \bold
V_1$ and $\nu_m = \lim \langle \nu_{m,n}:n \in B \rangle$ for $m <
m^*$.  We should prove that $\{n \in B:\text{if } m < m^*$ then $\nu_m
R_m \eta \equiv \nu_{m,n} R_m \eta\}$ is infinite. \nl
3) Left to the reader.
\enddemo
\newpage

\head {Glossary} \endhead  \resetall 
\bn
\ub{\S1 Trunk Controller}
\mn
\scite{it.1} (Definition)  Trunk controller, standard (trunk control),
$\aleph_1$-complete based, fully based, ${\Cal F}^{[\beta]}$, trivial, transparent
\mn
\scite{it.1a} (Definition)  (A trunk controller ${\Cal F}$ is) simple
(= purely regressively $\aleph_2$-c.c. on $S^2_1$), semi-simple 
(= $\aleph_2$-c.c. for pure extensions), (semi) simply based (all are
simple (except the first is semi simple)
\mn
\scite{it.2} (Definition)  ${\Cal F}$-forcing
\mn
\scite{it.2a} (Definition)  (An ${\Cal F}$-forcing $\Bbb Q$ is) clear (help put
together extensions) basic, straight (help put together $p_\ell
\le_{apr} q_\ell,p_q \le_{pr} p_2$, (used in \scite{mr.1}(1)) and transparent
\mn
\scite{it.4} (Definition)  ${\Cal F}$-iteration
\mn
\scite{it.5a} (Claim) Lim$_{\Cal F}(\bar{\Bbb Q})$ is a ${\Cal
F}$-forcing
\mn
\scite{it.6} (Claim)  in a ${\Cal F}$-iteration, $\beta < \gamma
\Rightarrow \Bbb P_\beta \lessdot P_\gamma$ naturally
\mn
\scite{it.7} (Claim) Preservation of clear (+ variant) and straight
\mn
\scite{it.7a} (Claim) Simplicity of ${\Cal F} +$ clearly of the ${\Cal
F}$-forcing $\Bbb Q$ impure (pure) $\aleph_2$-c.c. (+ variants)
\mn
\scite{it.8} (Claim) Existence of ${\Cal F}$-iteration
\mn
\scite{it.9} (Claim)  Associativity (of ${\Cal F}$-iterations)
\mn
\scite{it.10} (Discussion)
\mn
\scite{it.11} (Definition)  Pure decidability
\bn
\ub{\S2 Being ${\Cal F}$-pseudo c.c.c. (${\Cal F}$-psc) is preserved
by ${\Cal F}$-iterations}
\mn
\scite{ct.1} (Definition)  $\Bbb Q$ is ${\Cal F}$-psc, $({\Cal F},{\Cal P})$-psc,
clear, straight; consider $({\Cal F},{\Cal P})$
\mn
\scite{ct.2} (Definition)  $\Bbb Q$ is ${\Cal F}$-psc iteration as witnessed by
$\bar{\bold H}$, is essentially ... (except $Q_0$), semi-simple ${\Cal
F}$-psc strong
\mn
\scite{ct.2a} (Definition)  $\Bbb Q$ is strong ${\Cal F}$, psc, $\bar{\Bbb Q}$ is
strong
\mn
\scite{ct.3a} (Claim) Sufficient conditions for psc
\mn
\scite{ct.2b} (Definition)  ${\Cal F}$ is psc, strongly psc; cont,
Knaster, explicite, semi
\mn
\scite{ct.2c} (Claim) Implications
\mn
\scite{ct.2d} (Claim)  Basic fact on the explicit version
\mn
\scite{ct.3} (Claim) F-psc implies pure
$(\infty,\aleph_1)$-decidability for $\aleph_0$ anmes; the strong
version and the explicit version implies purely proper
\mn
\scite{ct.3b} (Remark)  On ``straight" and on stationary ${\Cal S}
\subseteq [\lambda]^{\aleph_0}$
\mn
\scite{ct.3c} (Remark) On $\aleph_1$-completeness
\mn
\scite{ct.4} (Lemma) Preserving psc under ${\Cal F}$-iteration
\mn
\scite{ct.4a} (Claim)  Preserving the explicit version under ${\Cal F}$-iteration
\mn
\scite{ct.6} (Definition) $({\Cal H},<_{\Cal H})$ is a c.c.c. witness
\mn
\scite{pr.2}  $\bold T_g$
\mn
\scite{pr.3}  $T = \text{ lim}_D\langle T_n:n < \omega \rangle,{\Cal
G} \subseteq {\Cal G}^{\bold V},{\Cal B} = ms - \text{lim}_D \langle {\Cal
B}_n:n < \omega \rangle$
\mn
\scite{p.5}  $\rho$ is ${\Cal G}$-continuous over $N$ (for $D$),
Null$_{{\Cal G},D}$
\bn
\ub{\S3 Nicer pure properness and pure decidability}
\mn
\scite{mr.1} (Claim)  Sufficient conditions for pure decidability
\mn
\scite{mr.2} (Claim) Purely properness
\bn
\ub{\S4 Averages by an ultrafitler and restricted non null set}
\mn
\scite{pr.1} (Definition)  For $\bold V_1 = \bold V[\bold r],\bold r$
random over $\bold V_1$, we consider extending an ultrafilter $D$ on
$\omega$ from $\bold V$ to an ultrafilter $D_1$ on $\omega$ from
$\bold V$ relevant to the randomness of $\bold r$.
\mn
\scite{pr.2} (Definition)  We define $\bold T_g$ as the set of $T
\subseteq {}^{\omega >} 2$ whose convergence to their Lebesgue measure
is bounded by $g$.
\mn
\scite{pr.3} (Definition)  We define $T = \text{ lim}_D \langle T_n:n
< \omega \rangle,{\Cal G}^v$ and ${\Cal B} = \text{ ms-lim}_D \langle
{\Cal B}_n:n < \omega \rangle$.
\mn
\scite{pr.4} (Claim) We note obvious things on lim$_D \langle T_n:n <
\omega \rangle$.
\mn
\scite{pr.5} (Definition)  We define ``$\rho$ is ${\Cal G}$-continuous
over $N$ for $D$", the ideal Null$_{{\Cal G},D}$ and Null$_{\Cal G}$.
\mn
\scite{pr.6} (Observation) If $\bold V \subseteq \bold V_1,({}^\omega
2)^{\bold V}$ not in $(\text{Null}_{\Cal G})^{{\bold V}_1}$ \ub{then}
$\bold V_1$ has no Cohen over $\bold V$ and the ultrafilter over
\scite{pr.5}.
\mn
\scite{pr.7} (Conclusion) We can extend an ultrafilter $D \in \bold V$
to an ultrafilter $D_1$ in $\bold V_1$ preserving ``$\bold r$ is
${\Cal G}$-continuous over $\bold V$".
\bn
\ub{\S5  On $\Bbb Q_{\bar D}$ and iterations}
\mn
\scite{br.1} (Definition)  We define $\bar D \in \bold{IF},\bar D \in
\bold{IUF} (\bar D = \langle D_\eta:\eta \in {}^{\omega >} \omega
\rangle,D_\eta$ a filter or ultrafilter on $\omega$ (non principle).
\mn
\scite{br.2} (Definition)  $\Bbb Q_{\bar D}$, a forcing notion, for
$\bar D \in \bold{IF},\underset\tilde {}\to \eta = \underset\tilde
{}\to \eta(\Bbb Q_{\bar D})$, the generic.
\mn
\scite{br.3} (Fact)  For $\bar D \in \bold{IUF},\Bbb Q_{\bar D}$ is
straight, clear, simple, $\sigma$-centered purely proper, ${\Cal
F}$-psc forcing when ${}^{\omega >} \omega \subseteq {\Cal F}$, with
$\underset\tilde {}\to \eta(\Bbb D_{\bar D})$ a generic real.
\mn
\scite{br.4} (Claim)  On 2-pure decidability fronts and absoluteness
for $\bar{\Bbb Q}_D$.
\mn
\scite{bs.5} (Claim)  New $f \in {}^\omega \omega$ run away from old
on Rang$[\underset\tilde {}\to \eta(\Bbb Q_{\bar D})]$
\mn
\scite{bs.6} (Claim)  (1) On $\Bbb Q_{\bar D} \lessdot \Bbb P * 
\Bbb Q_{\underset\tilde {}\to {\bar D}}$, when $D_\eta \le
{\underset\tilde {}\to D'_\eta}$. \nl
(2) Preserving $\underset\tilde {}\to G$-continuity.
\mn
\scite{bs.7} (Claim) On new $f \in \dsize \prod_n {}^{{\underset\tilde
{}\to \eta}(\Bbb Q_{\bar D})(n)} 2$ running away from old $\rho \in
{}^\omega 2$
\mn
\scite{bs.8} (Claim)  For $D \cup f$ on $\omega$, when does $\Bbb
Q_{\bar D}$ satisfy: in $\bold V^{{\Bbb Q}_{\bar D}},w_n
\subseteq
[{\underset\tilde {}\to \eta_{\Bbb Q_{\bar D}}}(n),
{\underset\tilde {}\to \eta_{\Bbb Q_{\bar D}}}(n)],|w_n| \le
\eta_{\Bbb Q_{\bar D}}(n)$ then $\cup\{w_n:n < \omega\}$ is disjoint to some member
of $D$.
\mn 
\scite{c.1} (Hypothesis)  $(CH + {\Cal F}^* + {\frak K}(0)$, with
Lim$(\bar R)$ - c.c.c.
\mn
\scite{c.2} (Definition)  Of ${\frak K}_\alpha,{\frak K}$.
\mn
\scite{c.3} (Definition)  $\bar{\Bbb Q}_1 \le_{\frak K} \bar{\Bbb
Q}_2$.
\mn
\scite{c.7} (Definition)  $ir(\Bbb Q)$, set of ${\frak p}$ giving an
autonomous description of a condition in the iteration
\mn
\scite{c.8} (Definition)  We define ${\frak K}^+_\alpha$ for the
content above.
\mn
\scite{c.9} (Remark)  
\mn
\scite{c.10} (Observation)  Collect the properties, not used.
\bn
{\S6 On a relative of Borel conjecture with large ${\frak b}$}
\mn
\scite{bt.5} (Hypothesis)
\mn
\scite{bt.8} (Definition) (1) We fix ${\frak K}(0)$,the candidates for
first forcing in the iterations as adding $\lambda$-randoms and define
$A(Q_0)$. \nl
(2) ${\frak K}_\alpha$ in this content. \nl
(3) ${\frak K}'_\alpha,{\frak K}'_{\ell,\alpha}$; mainly $\bar{\Bbb Q}
\restriction \beta \in {\frak K}^+_\beta$ for $\beta <
\alpha,\bar{\Bbb Q} \in {\frak K}'_\alpha$ but for $\alpha \ge
\lambda,Q_0$ has all random and they look internally a Sierpinski such
that if $\alpha < \lambda$ has $< \lambda$ randoms. \nl
(4) $\le_{{\frak K}'_{\ell,\alpha}},\le_{{\frak K}'_\alpha}$ (the
first new random is ${\Cal G}^{\bold V}$-continuous over the (mainly
smaller forcing). \nl
(5) ${\frak K}''_\alpha$ (mainly $\subseteq {\frak K}^+_\alpha$).
\mn
\scite{bt.9} (Observation) (1) The Sierpinski-ness of the randoms. \nl
(2) Essentially ${\Cal F}$-psc with 2-pure decidability over
$Q_0$. \nl
(3) Semi-simple + .
\mn
\scite{bt.11} (Claim)  Existence of extensions and appear bound of
increasing sequences for ${\frak K}'_{\le \alpha},{\frak K}''_{\le
\alpha}$.
\mn
\scite{bt.12} (Conclusion)  We get $\bar{\Bbb Q} \in {\frak K}_\alpha$
with $\alpha < \lambda$, using $\lambda$ or less of the randoms
manipulating ${\frak b},{\frak d}$ covering number for ?
\mn
\scite{bt.13} (Claim)  Similar to \scite{bt.12}, for $\alpha =
\lambda^+$.
\bn
{\S7 Continuing \cite{Sh:592}}
\mn
\scite{c7.2} (Theorem)  We find a forcing as in \cite{Sh:592} replacing
the null ideal by Null$_{\Cal G}$ but with ${\frak b} = {\frak d}$ is
quite small, e.g. in $\bold V^{\Bbb P}$, cov$(\text{Null}_{\Cal G}) =
\aleph_\omega,{\frak b} = {\frak d} = \aleph_2$.  
\mn
\scite{c7.3} (Remark)  Connection to \cite{Sh:619}.
\mn
\scite{c7.4} (Claim)   Amalgamation in ${\frak K}^+_\alpha$ or above
an amalgamation in ${\frak K}(0)$.
\mn
\scite{c.11} (Claim)  The fact needed for the induction step in
\scite{c.10} putting two ultrafilters together over amalgamated
forcings.
\mn
\scite{c.12} (Claim)  The generic real for $\bar{\Bbb Q} \in {\frak
K}'_{\ell,\alpha}$ are enough.
\bn
{\S8 On ``$\eta$ is ${\Cal L}$-big over $M$" 
\mn
\scite{p.1} (Definition)  (1) We define $\bold T$, the set of finitary
trees $\subseteq {}^{\omega >} {\Cal H}(\aleph_0)$. \nl
(2) $\bold R_{{\Cal T}_1,{\Cal T}_2}$ is the set of closed subsets of
Lim$({\Cal T}_2)$; also $\bold R_{\Cal T},\bold R_\alpha$ and $
\subseteq \lim({\Cal T}_1[R]) \times \lim({\Cal T}_2[R])$. \nl
(3) $R^{<k>}$. \nl
(4) $\bold Y$ is a set of $(f,{\Cal T}),D_Y$ for ${\Cal Y} \subseteq
\bold Y$. 
\mn
\scite{p.2} (Definition) (1) $D$ is $(f,{\Cal T})$-narrow for $f \in
\pi\{{\Cal P}({\Cal T}_n):n < \omega\}$. \nl
(2) $D$ is ${\Cal Y}$-narrow. \nl
(3) $\bold Z_M$ set of $(\eta,R)$ with $R \in M,\eta \in \lim({\Cal
T}_2[R])$. \nl
(4) $D$ is $(\eta,R)$-big over $M$. \nl
(5) $D$ is ${\Cal L}$-big over $M$. \nl
(6) $\eta$ is $R$-big or ${\Cal L}$-big over $M$.
\mn
\scite{p.3} (Claim)  Sufficient condition for extending the
ultrafitler $D_1 \in \bold V_1$ to $D_2 \in \bold V_2$ which is ${\Cal
L}$-big, ${\Cal Y}$-narrow over $\bold V_1$.
\mn
\scite{p.4} (Claim)  Equivalent condition to $B \notin J_{\Cal Y}$,
the narrowness ideal for ${\Cal Y} \subseteq \bold Y$.
\mn
\scite{p.5} (Claim)  Limit of ${\Cal Y}_\zeta$-narrow filters which is
$\dbcu_\zeta Y_\zeta$-narrow.
\mn
\scite{p.6} (Definition)  $(f,h,{\Cal T})$-narrow (return?) [releveling ${\Cal T}$]
\mn
\scite{p.7} (Definition)  ${\frak Z}_\mu$ set of $(\eta,\bar R)$ (?)
[return?]
\mn
\scite{p.8} (Claim)  (1) Extending $\bar D$ to preserve such that
$\eta \in \bold V_2$ is $R$-big on $\bold V_1[\eta_{\Bbb Q_{\bar
D}}]$. \nl
(2) Similarly for ${\frak Z}$.
\mn
\scite{p.9} (Definition) (1) Tree of subsets for $\bold T_g$.
 \nl
(2) More notation on trees.
\mn
\scite{p.10} (Claim) (1) ${\Cal G}$-continuous and $R_g$-bigness
equivalent. \nl
(2) Sufficient conditions for \scite{p.3} for ${\Cal Y} = G$. \nl
(3) Similarly ${\Cal Y} \ne \emptyset$ using \scite{p.9}(2). \nl

\newpage
    
REFERENCES.  
\bibliographystyle{lit-plain}
\bibliography{lista,listb,listx,listf,liste}

\shlhetal

\enddocument